\documentclass[onefignum,onetabnum]{siamart250211}
\usepackage{arydshln}
\usepackage{mathtools}

\usepackage{amsmath}
\usepackage{tikz}

\usepackage{array}
\usepackage{float}
\usepackage{amsfonts}
\usepackage{graphicx}
\usepackage{multirow}
\usepackage{epstopdf}
\usepackage{color}
\usepackage{algorithmic}
\usepackage[title]{appendix}

\usepackage[labelformat=simple]{subcaption}

\usepackage{amssymb}

\usepackage{xcolor, colortbl}
\definecolor{Gray}{gray}{0.90}
\definecolor{LightCyan}{rgb}{0.65,1,1}
\definecolor{Cyan}{rgb}{0.95,1,1}
\definecolor{GRAY}{gray}{0.75}

\usepackage[T1]{fontenc}
\usepackage[utf8]{inputenc}
\usepackage[font=small,labelfont=bf]{caption}
\usepackage[numbers,sort&compress]{natbib}

\newcommand{\bsub}{\begin{subequations}}
\newcommand{\esub}{\end{subequations}$\!$}

\renewcommand\thesection{\arabic{section}}
\renewcommand\thesubsection{\thesection.\arabic{subsection}}
\renewcommand\thesubsubsection{\thesubsection.\arabic{subsubsection}}
\renewcommand\theequation{\thesection.\arabic{equation}}

\renewcommand{\theequation}{\arabic{section}.\arabic{equation}}

\usepackage{wrapfig}
\usepackage{mathtools}

\usepackage[format=plain,
            labelfont={it},
            textfont=it]{caption}

\newtheorem{prop}{Proposition}

\newcommand{\R}{\mathbb{R}}
\newcommand{\zero}{\v{0}}
\renewcommand{\v}[1]{\mathbf{#1}}

\newcommand{\x}{\v{x}}
\newcommand{\q}{\v{q}}

\newcommand{\y}{\v{y}}

\newcommand{\eps}{\varepsilon}

\title{Synchronized Memory-Dependent Intracellular Oscillations for a
  Cell-Bulk ODE-PDE Model in $\R^3$} \author{Merlin Pelz\thanks{School
    of Mathematics, University of Minnesota, Minneapolis, MN, USA
    (corresponding author), \texttt{email:mpelz@umn.edu}} \and Michael
  J. Ward\thanks{Department of Mathematics, University of British
    Columbia, Vancouver, BC, Canada, \texttt{email:ward@math.ubc.ca}}}

\headers{Memory-Dependent Dynamics for a Cell-Bulk Model in 3-D}
  {M. Pelz, M. J. Ward}

\ifpdf
\hypersetup{
  pdftitle={Memory-Dependent Dynamics for a Cell-Bulk Model in 3-D},
  pdfauthor={M. Pelz, M. J. Ward}}
\fi

\begin{document}

\maketitle

\begin{abstract}
We develop a hybrid asymptotic-numerical theory for a cell-bulk
ODE-PDE model in $\R^3$ to analyze the long-range diffusive coupling
of spatially localized, dynamically active cells. In the limit of
small cell radius $\eps\ll 1$, we derive a novel integro-differential
ODE system for the intracellular dynamics, in which the 3-D bulk
diffusion field is effectively replaced by a time-memory kernel whose
strength depends on the inter-cell distances. To compute long-time
solutions efficiently, we implement a fast time-marching scheme for
this integro-differential system that is inspired by the {\em
  sum-of-exponentials} method. This scheme is applied to Sel'kov
intracellular kinetics, which is a conceptual two-component model for
a conditional oscillator. Our analysis reveals that, in contrast to
the more strongly coupled 2-D case, each cell in $\R^3$ experiences
only a weak ${\mathcal O}(\eps)$ coupling effect from the other cells
in the population. Consequently, individual cells that exhibit
oscillatory intracellular dynamics can only synchronize over very long
time scales of order ${\mathcal O}(\eps^{-1})$ when they are separated
by ${\mathcal O}(1)$ distances. Moreover, cells whose kinetic
parameters are ${\mathcal O}(\eps)$ close to their marginal stability
threshold will be influenced by the weak ${\mathcal O}(\eps)$
cell-cell coupling from the entire population. For this latter
situation, in order to resolve the effect of these subtle interactions
on collective intracellular dynamics, we develop a novel matrix
perturbation theory for the linear stability of steady-state solutions
of the cell-bulk model. Finally, we validate our analytical
predictions of the perturbed marginal stability thresholds of
steady-state solutions using the time-marching numerical scheme,
illustrating various types of synchronized intracellular dynamics that
can emerge over long time scales.
\end{abstract}

\section{Introduction}\label{sec:intro}

 Collective behavior mediated by a diffusing chemical that
 intermittently binds with spatially localized, dynamically active
 units in a multi-dimensional spatial domain arises across diverse
 applications in biology and chemical physics.  In particular,
 microbial cell-to-cell communication and quorum-sensing phenomena
 often rely on the secretion and feedback of extracellular signaling
 molecules called autoinducers (cf.  \cite{QS_whiteley} and
 \cite{ingalls}). In this bulk-mediated process, cells release signals
 into a shared environment to form a collective autoinducer field,
 allowing individual cells to adjust their intracellular dynamics
 based on the total population. Key examples of this mechanism include
 cyclic adenosine monophosphate (cAMP) triggering aggregation in the
 amoeba {\em Dictyostelium discoideum} (cf.~\cite{gregor2010},
 \cite{nandy1998}, \cite{goldbeter1997biochemical}), acetaldehyde
 (Ace) coordinating glycolytic oscillations in yeast colonies
 (cf.~\cite{dano}, \cite{dano2}, \cite{de2007dynamical}), and acylated
 homoserine lactones (AHLs) inducing bioluminescence behavior in the
 marine bacterium {\em Aliivibrio fischeri}
 (cf.~\cite{taga2003chemical}).  To study cell-cell interactions via a
 diffusive medium in a more controlled setting, early efforts in
 modeling living systems through synthetic biology demonstrated
 diffusion-mediated synchronization among engineered genetic clocks
 \cite{danino2010}.  Recent progress in synthetic biology has focused
 on engineering chemical communication interfaces in completely
 artificial cells that are coupled exclusively by a bulk medium
 \cite{mukwaya2021} (see \cite{karoui2022chemical} for a recent
 survey).

  Collective behavior of spatially localized dynamically active units
  coupled by a bulk diffusion field also occurs in chemical
  physics. In particular, chemically coated small spherical beads
  immersed in a Belousov-Zhabotinsky (BZ) reaction mixture trigger
  synchronous chemical oscillations when the number of particles
  exceeds a quorum-sensing threshold (cf.~\cite{tinsley1},
  \cite{tinsley2}, \cite{taylor2}). In other applications, collective
  oscillations can emerge in chemical emulsions where BZ reagents
  encapsulated in microfluidic droplets communicate via the diffusion
  of reaction intermediates through a bulk medium that consists of a
  continuous oil phase (cf.~\cite{torbensen}, \cite{toiya}). 

  The common feature in each of these applications is that, through
  the secretion and absorption of a diffusible species, a spatially
  localized dynamically active ``unit'' can either initiate or
  synchronize its internal dynamics with other spatially distant units
  through the mediation of a bulk diffusion field.  Cell-bulk ODE-PDE
  models provide an optimal theoretical framework for investigating
  collective dynamics driven by bulk diffusion fields, established in
  3-D by \cite{muller2006, muller2013} and in 2-D by \cite{gou2d,
    smjw_diff}. From a mathematical viewpoint, they are
  compartmental-reaction diffusion (CRD) systems, and they can model
  how intracellular oscillations within spatially segregated,
  dynamical units are mediated by the simultaneous production and
  absorption of diffusing chemical species. This formulation directly
  aligns with Kuramoto's foundational analysis of oscillator
  synchronization via diffusive coupling \cite{kuramoto1995}, and
  contributes to a largely open research frontier highlighted in
  \cite{show_epstein} for theoretical frameworks that couple discrete
  dynamical elements.

  Most prior investigations of CRD systems have primarily
  focused on pattern formation properties near bifurcation points,
  such as characterizing the onset of collective dynamics via Hopf
  bifurcations of steady-state solutions. In 1-D, for reaction
  compartments located at the domain endpoints, the marginal conditions
  for the onset of either in-phase or anti-phase oscillations for
  various reaction kinetics with either one or two bulk diffusing
  species \cite{gomezmarin,gou2015, gou2016, gou2017, pelz2026a} or
  with a time-delayed kinetics \cite{xu} have been
  analyzed. Amplitude equations characterizing branching properties
  for these oscillatory instabilities, obtained from a weakly
  nonlinear analysis, have been derived in \cite{paquin_1d}. More
  recently, in \cite{pelz2026b} a rigorous analysis of bifurcating
  in-phase and anti-phase oscillations from a symmetric steady-state
  solution for a scalar bulk diffusing species with dynamically active
  Wentzell boundary conditions was undertaken. With two bulk diffusing
  species, it was shown in \cite{turing1d} that large boundary
  reaction rate ratios can trigger symmetry-breaking bifurcations to
  stable asymmetric patterns. Parallel 2-D
  investigations with one bulk diffusing species in bounded domains containing
  small spatially segregated dynamically active cells have been used
  to model either quorum sensing behavior or the onset of oscillatory
  instabilities through a Hopf bifurcation \cite{gou2d, smjw_diff,
    ridgway, q_survey,smjw_quorum}.

  In contrast, there has been more limited work on analyzing
  large-scale spatiotemporal patterns of cell-bulk ODE-PDE models
  away from bifurcation points. For a cell-bulk model in $\R^2$, an
  integro-differential system of ODEs characterizing intracellular
  dynamics of a population of dynamically active units was derived in
  \cite{pelzward2D} by extending the strong localized perturbation
  theory surveyed in \cite{ward2018spots} to a time-dependent
  setting. By using a novel numerical method, the synchronizaion of
  intracellular oscillations as quantified by the Kuramoto-order
  parameter was illustrated for various spatial configurations of
  cells. The 2-D cell-bulk framework of \cite{pelzward2D} was used in
  \cite{pelz26} to model a specific application in neuroscience.
  
\begin{figure}[htbp]
  \centering
      \begin{subfigure}[b]{0.48\textwidth}
      \includegraphics[width=\textwidth]{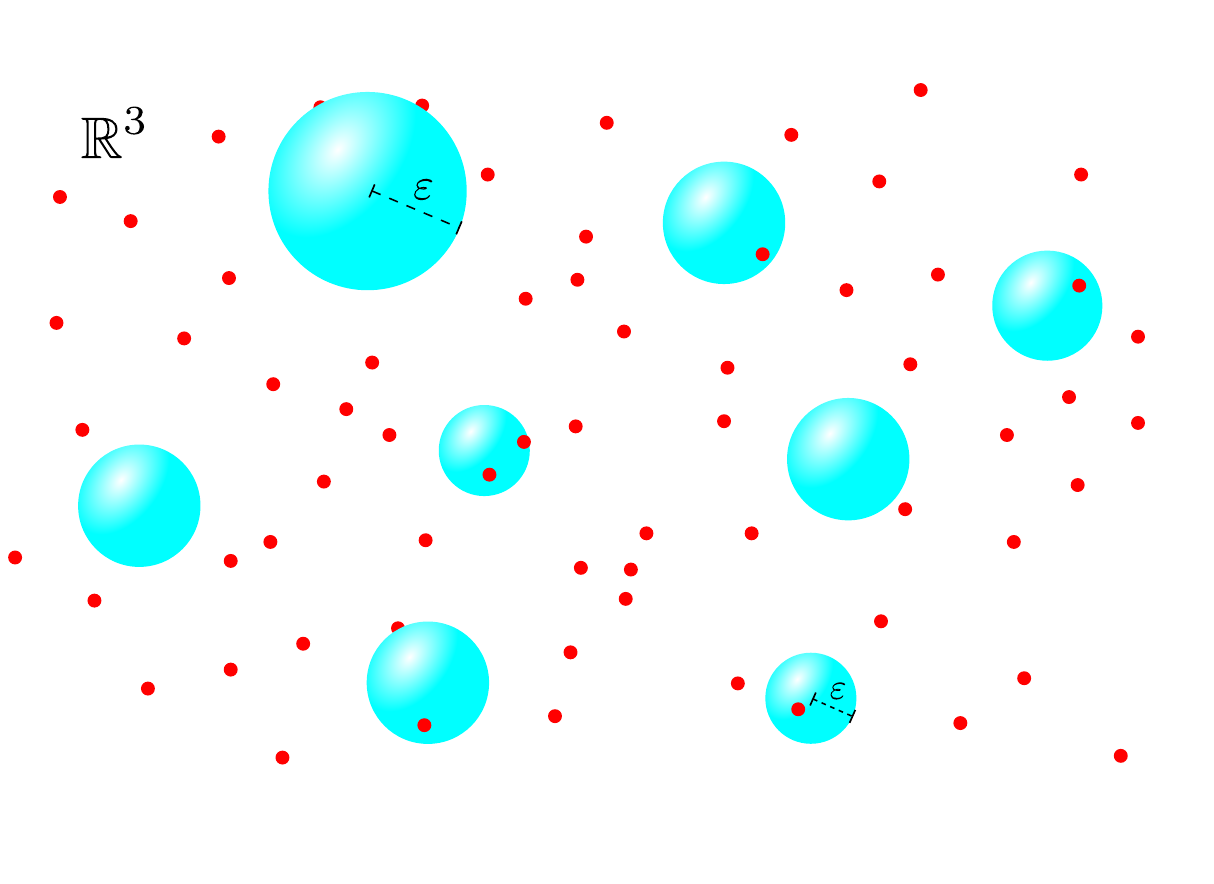}
        \label{fig:onebulk_3d}
    \end{subfigure}
      \begin{subfigure}[b]{0.48\textwidth}
      \includegraphics[width=\textwidth]{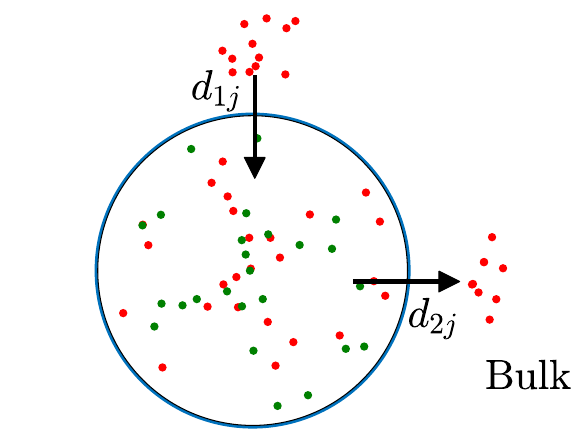}
        \label{fig:onebulk_zoom}
    \end{subfigure}
    \caption{Left: Schematic showing dynamically active
        signaling compartments (in cyan) in 3-D coupled by a bulk
        diffusion field (red dots) for the dimensionless cell-bulk
        model \eqref{DimLess_bulk}.  Right: A zoomed-in illustration
        of the compartments showing two intracellular species (green
        and red dots). Only the red species, labeled by $u_{1j}$ (see
        \eqref{DimLess_bulkb}), is secreted into and absorbed from the
        bulk, as regulated by the influx and efflux permeabilities
        $d_{1j}$ and $d_{2j}$.}
    \label{fig:onebulk_scheme}
\end{figure}

The goal of this paper is to study large-scale spatial-temporal
pattern formation for a cell-bulk model with one bulk species in $\R^3$
and with $N$ dynamically active cells, as inspired by \cite{muller2006}
and \cite{muller2013}. In this model, the dynamically active cells are
assumed to be a collection of non-overlapping spheres of a common
radius $R_0$ and are labeled as $\Omega_j\equiv \lbrace{\v{X} \, \vert
  \, |\v{X}-\v{X}_j|\leq R_0\rbrace}$ for
$j\in\lbrace{1,\ldots,N\rbrace}$.  In the bulk, or extracellular,
region $\R^{3} \setminus \cup_{j=1}^{N}\,\Omega_j$, the bulk
concentration $\mathcal{U}(\v{X},T)$ with constant bulk diffusivity
$D_B>0$ and degradation $k_B>0$ is assumed to satisfy
\begin{subequations}\label{Dim_bulk} 
\begin{align}
 \partial_T \mathcal{U} = & D_B \, \Delta \,
     \mathcal{U} - k_B \, \mathcal{U}\,, \quad T>0\,,\quad 
    \v{X} \in \R^3\,\setminus \cup_{j=1}^{N}\,\Omega_j\,, \label{Dim_bulka}\\
    D_B \, \partial_{n_{\v{X}}} \, \mathcal{U}  &=
  \beta_{1j} \, \mathcal{U} - \beta_{2j} \, \mu_{1j}\,, \quad
   \v{X} \in \partial \Omega_j\,,\qquad j \in \lbrace{1, \ldots, N\rbrace}\,,
                                                  \label{Dim_bulkb}
\end{align}
In this model, the influx and efflux into and out of the
$j^{\text{th}}$ cell are regulated by the permeabilities
$\beta_{1j}>0$ and $\beta_{2j}>0$. In \eqref{Dim_bulkb},
$\partial_{n_{\v{X}}}$ denotes the outer normal derivative pointing into
the bulk region, while $\mu_{1j}$ is an intracellular species specified
below.

Within each cell we assume that $m$ intracellular species
$\v{\mu}_j = (\mu_{1j},\ldots,\mu_{mj})^T$ interact with reaction
kinetics $\v{F}_j$, and we neglect any spatial gradients within the
cells. Moreover, we assume that the cell boundary is permeable to only
one species, labeled by $\mu_{1j}$. The intracellular dynamics in the
$j^{\mathrm{th}}$ cell are coupled to the bulk diffusion field
\eqref{Dim_bulk} through an integration of the diffusive flux across the
boundary of the cell as
\begin{equation}\label{Dim_Intra}
\begin{split}
  \frac{d\v{\mu}_j}{dT} & = k_R\, \mu_c \, \v{F}_j
  \left( \v{\mu}_j/\mu_c  \right)  + \v{e}_1 \int_{\partial \Omega_j}\,
  \left( \beta_{1j} \, \mathcal{U} - \beta_{2j} \, \mu_{1j} \right) \,
  \text{d}S_{\v{X}} \,, \quad j\in \lbrace{1,\ldots,N\rbrace}\,,
\end{split}
\end{equation}
\end{subequations}
where $\v{e}_1 \equiv (1,0,\ldots,0)^T$. In this model, $k_R > 0$ is a
representative dimensional intracellular reaction rate, $\mu_c>0$ is a
typical magnitude of the vector $\v{\mu}_j$, $\beta_{2j}$ is the
efflux rate of $\mu_{1j}$ per unit area of the cell surface, while
$\beta_{1j}$ is the influx rate times the sphere volume per area of the
cell surface. This influx rate regulates the feedback into the
$j^{\mathrm{th}}$ cell from the bulk diffusion field that is generated
by the entire collection of cells.

We will analyze the cell-bulk model in the asymptotic limit where the
common radius $R_0$ of the signaling compartments is small relative to
the minimum inter-cell separation distance
$L\equiv \min{|\v{X}_i-\v{X}_j|}$ for $i\neq j$. As such, we introduce
a small parameter by $\eps \equiv R_0/L \ll 1$. In Appendix
\ref{app:nondim} we non-dimensionalize the cell-bulk model
\eqref{Dim_bulk} to obtain that the dimensionless bulk field
$U(\v{x},t)$ satisfies
\begin{subequations}\label{DimLess_bulk}
\begin{align}
 \partial_t U =&\,  D \, \Delta U -  \sigma \, U\,, \quad
   t > 0\,,\quad \x \in \R^3 \setminus \cup_{j=1}^{N}\,
   \Omega_{\varepsilon_j}\,; \quad U(\v{x},0)=0 \,, \label{DimLess_bulka} \\
  \varepsilon  D\, \partial_n U  &= d_{1j} \, U -
                         \frac{d_{2j}}{\eps} \, u_{1j}\,,
      \quad \x \in \partial \Omega_{\varepsilon_j}\,,
                               \qquad j \in \lbrace{1,\ldots, N\rbrace}\,.
    \label{DimLess_bulkb}
\end{align}
This bulk diffusion field is coupled to the dimensionless
intracellular species $\v{u}_j\equiv(u_{1j},\ldots,u_{mj})^T$
within the $j^{\text{th}}$ cell by
\begin{align}\label{DimLess_Intra}
  \frac{\text{d} \v{u}_j}{\text{d}t} & =  \,
 \v{F}_j \left( \v{u}_j  \right)  + \frac{\v{e}_1}{\eps}
    \int_{\partial \Omega_{\varepsilon_j}}\, \left( d_{1j} \, U -
    \frac{d_{2j}}{\eps} \, u_{1j}\right) \,\,
              \text{d}S_{\v{x}} \,,\qquad j \in \lbrace{1, \ldots, N \rbrace}\,.
\end{align}
\end{subequations}
In \eqref{DimLess_bulk}, the small signaling cells are labeled by
$\Omega_{\varepsilon_j}\equiv \lbrace{\v{x} \, \vert \, \,
  |\v{x}-\v{x}_j| \leq \varepsilon\rbrace}$, and we assume that
they are spatially segregated in the sense that
$\mathrm{dist}(\v{x}_j,\v{x}_k)={\mathcal O}(1)$ for $j\neq k$. The
dimensionless parameters in the cell-bulk model \eqref{DimLess_bulk}
are
\begin{equation}\label{dim:param}
  D \equiv \frac{D_B}{k_R L^2} \,, \quad
  d_{1j} \equiv \varepsilon \frac{\beta_{1j}}{k_R L} = {\mathcal O}(1) \,, \quad
  d_{2j} \equiv \varepsilon^2 \frac{\beta_{2j} L^2}{k_R} =
  {\mathcal O}(1) \,, \quad
   \sigma \equiv \frac{k_B}{k_R} \,,
 \end{equation}
 where $D$ and $\sigma$ are referred to as the effective bulk
 diffusivity and effective bulk degradation parameter,
 respectively. In \eqref{dim:param}, the permeability parameters
 $\beta_{1j}$ and $\beta_{2j}$ are assumed to have a different scaling
 in $\eps$ in order to ensure that the small cells contribute an
 ${\mathcal O}(1)$ diffusive flux into the bulk medium (see Appendix
 \ref{app:nondim}). A schematic figure of the dimensionless cell-bulk
 model \eqref{DimLess_bulk} is shown in Fig.~\ref{fig:onebulk_scheme}
 for the case of two intracellular species.  Our goal for
 \eqref{DimLess_bulk} is to develop a hybrid asymptotic-numerical
 method to investigate how the bulk diffusion field, generated by the
 entire population of cells, initiates and synchronizes large-scale
 intracellular oscillations.

In the limit $\eps\to 0$ of common small cell radius, in \S
\ref{2:ode_memory_trans} we use a strong localized perturbation
approach, based on the method of matched asymptotic expansions in a
time-dependent setting, to derive a new integro-differential ODE
system characterizing how the intracellular dynamics within each cell
are coupled by the bulk diffusion field. This integro-differential
system, summarized in Proposition \ref{prop:3ode}, shows that the
intracellular dynamics within each cell is only weakly coupled through
${\mathcal O}(\eps)$ terms to the entire cell population. The
${\mathcal O}(\eps)$ coupling terms are nonlocal in time and arise
from the combined effect of the rapid ${\mathcal O}({1/r})$ spatial
decay of the free-space Green's function at a radial distance $r$ from
the center of a cell, and the memory-dependent Green's kernel in
$\R^3$.  In modeling scenarios where bulk-diffusion in 3-D is believed
to provide the mechanism for intercellular coupling of spatially
segregated dynamically active units with permeable boundaries, the
weakly coupled memory dependent system formally derived in
\eqref{2:reduced} provides an alternative to the more phenomenological
Kuramoto-type models discussed in \cite{viana2020} and in \S 4.1 of
the review article \cite{rodriguez}.  Our limiting
integro-differential ODE system is an extension of that derived in
\cite{muller2013} for one intracellular species and where the bulk
diffusion field has no degradation.

In \S \ref{sec:ss} we characterize the set of steady-state limiting
values for the integro-differential ODE system. In Appendix
\ref{app:strong} we show using strong localized perturbation theory
that this set corresponds to all possible steady-state solutions of
the cell-bulk model \eqref{DimLess_bulk}. In \S \ref{sec:stab} we use
strong localized perturbation theory to derive a nonlinear matrix
eigenvalue problem whose spectrum characterizes the linear stability
properties of steady-state solutions to \eqref{DimLess_bulk}. This
eigenvalue problem shows that the linear stability of a steady-state
solution within a particular cell is only weakly influenced through
${\mathcal O}(\eps)$ perturbations by all of the other cells in the
population. For arbitrary intracellular kinetics, in \S
\ref{sec:eig_pert} we use matrix perturbation theory to derive some
general results for how discrete eigenvalues associated the
steady-state of a single cell are affected by the ${\mathcal
  O}(\eps)$ perturbation from the entire cell population.

In \S \ref{sec:selkov} we provide detailed linear stability results
regarding the possibility of Hopf bifurcations for the two-component
kinetics used in \cite{MNS} to model chemical oscillations. These
kinetics are a modification of the original Sel’kov kinetics
\cite{Selkov} modeling glycolysis oscillations.  By implementing the
matrix perturbation framework of \S \ref{sec:eig_pert} for this
conditional oscillator, we determine the marginal threshold conditions
for oscillations for various spatial configurations of cells. More
specifically, we will consider two non-identical cells and spatial
configurations of multiple identical cells where a certain Green's
matrix has a constant row sum, such as when cells are centered at the
vertices of platonic solids on the surface of a sphere. We will also
determine marginal stability thresholds for such uniformly spaced
cells in the presence of an additional pacemaker cell located at the
center of the sphere.

A direct numerical computation of the integro-differential system in
Proposition \ref{prop:3ode} is neither efficient or practical owing to
the fact that ${N(N+1)/2}$ memory-dependent convolution integrals,
each with integrable singularities, would have to be computed
numerically up to time $t$ in order to advance the solution one time
step to time $t+\Delta t$. As
such, a different numerical approach is needed to efficiently compute
solutions to \eqref{2:reduced} over long time intervals. Using a
related approach as done for the 2-D cell-bulk model in
\cite{pelzward2D}, in \S \ref{sec:algorithm} and Appendix
\ref{app:numerics} we develop a fast time-marching numerical scheme
for the integro-differential ODE system that is inspired by the {\em
  sum-of-exponentials} method of \cite{greengard_2015} and
\cite{beylkin_2005} (see also \cite{beylkin_2010}). In this method,
the Laplace transform of the Green's kernel in \eqref{2:reduced} is
represented, with a controlled error tolerance, by a sum of
exponentials obtained from a discretization along a path in the
Laplace transform space. Since the cell-cell interactions can then be
replaced (up to calculating certain local contributions) by a sum of
convolution kernels with exponential weights, a recursive
time-stepping scheme can be developed for the integro-differential ODE
system. This eliminates the need to recalculate the system's entire
past history to advance one time-step. The details on our numerical
time-marching scheme for \eqref{2:reduced} with Sel'kov kinetics are
described in Appendix \ref{app:numerics}. In \S
\ref{sec:selkov_experiments} we use the time-marching numerical scheme
to confirm some of our analytical predictions of the perturbed
marginal stability thresholds for Sel'kov kinetics that were derived
in \S \ref{sec:selkov}, and we illustrate various types of
synchronized intracellular dynamics that can emerge over long time
scales.

Finally, in \S \ref{sec:discussion} we discuss the key differences between
the results derived herein for the synchronization of spatially
segregated dynamically active units coupled by a bulk diffusion field in
$\R^3$ with the those obtained in \cite{pelzward2D} for the corresponding
problem in $\R^2$. A few open problems in a 3-D setting are discussed.

\section{Derivation of Integro-Differential ODE System for Intracellular Kinetics}\label{2:ode_memory_trans}

In the limit $\eps\to 0$ we use the method of matched asymptotic
expansions on the cell-bulk model \eqref{DimLess_bulk} to derive an an
integro-differential ODE system for the intracellular kinetics that is
valid on time-scales $t\gg {\mathcal O}(\eps^2)$. However, in order to
provide initial conditions for this limiting system we first must
analyze the transient solution near each cell that occurs on a short
${\mathcal O}(\eps^2)$ time-scale.

\subsection{Transient Analysis Near the Cells}\label{sec:transient}

We first analyze the transient solution near the $j^{\mathrm{th}}$
cell on the time-scale $t={\mathcal O}(\eps^2)$ when
$U(\v{x},0)=0$. As such, we introduce $\tau$ by $t=\eps^2\tau$. The
far-field behavior as $\tau\to\infty$ of this solution provides the
initial conditions for the long-time dynamics analyzed below in \S
\ref{2:ode_memory}.

Near the $j^{\mathrm{th}}$ cell, we introduce the local coordinates
$\y=\eps^{-1}(\v{x}-\v{x}_j)$ with $\rho=|\v{y}|$, and set
$\tau={t/\eps^2}$. For $t\ll 1$, we use
$u_{1j}(t)=u_{1j}(0)+ t u_{1j}^{\prime}(0)+\cdots$, which becomes
$u_{1j}(t)\approx u_{1j}(0)+\eps^2 \tau u_{1j}^{\prime}(0)$ on the
slow time-scale.  From \eqref{DimLess_bulk}, we obtain that
$V(\rho,\tau)=U(\v{x}_j+\eps\v{y}, \eps^2\tau)$ satisfies the
approximate local problem
\begin{equation}\label{3:local_trans}
\begin{split}
  \partial_{\tau} V &=D \left(\partial_{\rho\rho} V + \frac{2}{\rho}
    \partial_\rho V\right) -\eps^2\sigma V \,,
  \quad \rho>1\,, \,\,\, \tau\geq 0\,;
  \quad V(\rho,0) = 0  \,, \\
  D \partial_\rho V &= d_{1j} V - \frac{d_{2j}}{\eps}
  \left(u_{1j}(0) +\eps^2 \tau u_{1j}^{\prime}(0)\right) \,, \quad \mbox{on}
  \,\,\, \rho=1 \,.
\end{split}
\end{equation}
To determine $V(\rho,\tau)$ we take the Laplace transform of
\eqref{3:local_trans} in $\tau$ and solve the resulting boundary value
problem. Labeling
$\hat{V}(\rho,s)\equiv {\mathcal L}\left[V(\rho,\tau)\right]$ as the Laplace
transform of $V(\rho,\tau)$, we obtain that
\begin{equation}\label{3:lap_trans}
  \hat{V}(\rho,s) = \frac{d_{2j}}{\eps s\rho }\left( u_{1j}(0) +
    \frac{\eps^2 u_{1j}^{\prime}(0)}{s} \right) 
  \frac{e^{-\sqrt{{(s+\eps^2\sigma)/D}}(\rho-1)}}{\left(d_{1j} + D + \sqrt{(s+\eps^2
        \sigma) D}\right)}\,,
\end{equation}
where the branch cut is taken along 
$\mathrm{Re}(s)\leq -\eps^2\sigma$ with $\mathrm{Im}(s)=0$.

In Appendix \ref{app:laplace}, we analyze the far-field behavior of
$V(\rho,\tau)$ for $\tau\gg 1$ and $\rho$ fixed, by first
approximating $\hat{V}(\rho,s)$ in \eqref{3:lap_trans} near $s=0$,
which corresponds to the singularity of $\hat{V}$ with the largest
real part. We then invert the resulting approximation to obtain
$V_{0}(\rho,\tau)$ for $\tau\gg 1$, and replace $\tau={t/\eps^2}$.  In
this way, we conclude that on the intermediate scaling regime
${\mathcal O}(\eps^2)\ll t\ll 1$, we have
\begin{subequations}\label{3:v_final}
\begin{equation}\label{3:v_mono}
  V(\rho,t) \sim \frac{A_j}{\rho} + B_{j} \,,
\end{equation}
where
\begin{gather}
  A_j = \frac{d_{2j}}{\eps(d_{1j}+D)}
  \left( u_{1j}(0)+t u_{1j}^{\prime}(0) + \frac{{\eps d_{1j}/\sqrt{\pi D t}}}
    {\left(d_{1j}+D\right)} \left[ u_{1j}(0) + t \left(\sigma u_{1j}(0) +
  2 u_{1j}^{\prime}(0) \right) + {\mathcal O}(t^2) \right]\right) \,,
  \label{3:v_aj}\\
  B_j = - \frac{d_{2j}}{\left(d_{1j}+D\right)\sqrt{\pi D t}}
  \left[ u_{1j}(0) + t \left( \sigma u_{1j}(0) + 2 u_{1j}^{\prime}(0)\right)
    + {\mathcal O}(t^2)\right] \,. \label{3:v_bj}
\end{gather}
\end{subequations}
In \S \ref{2:ode_memory} we show that \eqref{3:v_bj} agrees precisely
with the limiting behavior as $t\to 0^{+}$ for the solution to our
integro-differential ODE system.

Finally, to determine $u_{1j}^{\prime}(0)$, as needed in
\eqref{3:v_aj} and \eqref{3:v_bj}, we set $U(\x,0)=0$ in
\eqref{DimLess_Intra} to obtain
\begin{equation}\label{3:uprime}
  \quad u_{1j}^{\prime}(0)=\v{e}_1^T \frac{ {\text{d} \v{u}_j}}{\text{d}t}
    \vert_{t=0} = \v{e}_1^{T} \v{F}_j \left( \v{u}_j(0)  \right)  -
    4\pi d_{2j} u_{1j}(0) \,.
\end{equation}

\subsection{Derivation of Integro-Differential ODE
  System}\label{2:ode_memory}

We now derive the integro-differential ODE system that is valid when
$t\gg {\mathcal O}(\eps^2)$.

For $t={\mathcal O}(1)$, in the $j^{\mathrm{th}}$ inner region we introduce
the local coordinates $\y=\eps^{-1}(\x-\x_j)$, $\rho=|\y|$, and
$V(\y,t)=U(\x_j+\eps\v{y},t)$. From \eqref{DimLess_bulka}, we obtain the
leading-order quasi-steady problem
\begin{equation}\label{2:inn_qss}
  \Delta_{\y} V = 0 \,, \quad \mbox{for} \quad \rho>1 \,; \qquad
  D V_{\rho}= d_{1j}V - \frac{d_{2j}}{\eps} u_{1j} \,, \quad \mbox{on}
  \quad \rho=1\,,
\end{equation}
which has the radially symmetric solution
\begin{subequations}\label{2:inn_sol_full}
\begin{equation}\label{2:inn_sol}
  V(y,t) =\frac{A_{j}}{\rho} + B_j \,,
\end{equation}
where $A_j=A_j(t,\eps)$ and $B_j=B_j(t,\eps)$ are to be determined. We
obtain that
\begin{equation}\label{2:ajbj_new}
  A_{j} = \frac{d_{2j} u_{1j}}{\eps(d_{1j}+D)} - \frac{d_{1j}B_j}{d_{1j}+D}
  \,.
\end{equation}
\end{subequations}

Next, we substitute \eqref{2:inn_sol} into the intracellular dynamics
\eqref{DimLess_Intra} to obtain 
\begin{equation}\label{DimLess_Intra_flux}
  \frac{\text{d} \v{u}_j}{\text{d}t} - \v{F}_j \left( \v{u}_j  \right)
  = \frac{\v{e}_1}{\eps}
      \int_{\partial \Omega_{\varepsilon_j}}\, D V_{\rho}\vert_{\rho=1} \,
    \text{d}S_{\v{x}}=
     \frac{\v{e}_1}{\eps} \int_{\partial \Omega_{\varepsilon_j}}\,
    \left[-DA_j \right] \text{d}S_{\v{x}} \,,
  \end{equation}
for each $j\in\lbrace{1,\ldots,N\rbrace}$. Then, by using
$|\partial\Omega_{\eps_j}|=4\pi \eps^2$ and \eqref{2:ajbj_new}, we
conclude that the intracellular dynamics is determined in terms of
the unknown $B_{j}(t)$ by
\begin{equation}\label{2:intra_full}
  \frac{\text{d} \v{u}_j}{\text{d}t}= \v{F}_j \left( \v{u}_j  \right)
  - \kappa_j u_{1j} \v{e}_1 + \frac{\eps d_{1j}}{d_{2j}} \kappa_j B_{j}
  \v{e}_1 \,, \quad \mbox{where} \quad \kappa_j \equiv \frac{4\pi D d_{2j}}
  {d_{1j} + D} \,.
\end{equation}

We now determine $B_{j}(t)$ by matching the inner solutions near the
cells to an outer bulk solution. Upon writing \eqref{2:inn_sol} in outer
variables, by using $\eps \rho=|\v{x}-\v{x}_j|$, we obtain that the
bulk solution $U(\v{x},t)$ satisfies
\begin{subequations}\label{2:bulk_all}
\begin{equation}\label{2:bulk_pde}
  \begin{split}
   \partial_t U = D \Delta U &- \sigma U \,, \quad \v{x}\in \R^3\backslash
    \lbrace{\v{x}_1,\ldots,\v{x}_N\rbrace} \,; \qquad U(\v{x},0)=0 \,, \\
    U \sim \frac{b_j}{4\pi D}\frac{1}{|\v{x}-\v{x}_j|} &+ B_{j}
    \,, \quad \mbox{as} \quad \v{x}\to \v{x}_j \,,
    \quad j\in \lbrace{1,\ldots,N\rbrace} \,,
  \end{split}
\end{equation}
with $U(\v{x},t)\to 0$ as $|\v{x}|\to \infty$, and where $b_j=b_j(t)$ is
defined by
\begin{equation}\label{2:bulk_cj}
  b_{j}(t) \equiv \kappa_j u_{1j}(t) -\frac{\eps \kappa_j d_{1j}}{d_{2j}} B_j(t)\,.
\end{equation}
\end{subequations}

To determine $B_{j}(t)$, we take the Laplace transform of
\eqref{2:bulk_all} to obtain that
$\hat{U}(\x,s)={\mathcal L}\left[U(\x,t \right)]$ satisfies
\begin{subequations}\label{2:vk}
\begin{align} 
  \Delta \hat{U} &- \frac{(\sigma+s)}{D} \hat{U} = -\sum_{k=1}^{N}
  \frac{\hat{b}_k}{D} \delta(\x-\x_k) \,, \quad \x\in \R^3\,; \quad
           \quad \hat{U}\to 0  \quad \mbox{as} \quad |\v{x}|\to \infty \,,
  \label{2:vk_1}  \\
  \hat{U} &\sim \frac{\hat{b}_j}{4\pi D} \frac{1}{|\x-\x_j|}
  +\hat{B}_{j} \,, \quad \mbox{as} \quad \v{x}\to \v{x}_j \,,
    \quad j\in \lbrace{1,\ldots,N\rbrace} \,. \label{2:vk_2}
\end{align}
\end{subequations}
The solution to \eqref{2:vk_1} is
\begin{equation}\label{2:vksol}
  \hat{U} = \frac{1}{4\pi D} \sum_{k=1}^{N} \frac{\hat{b}_k}{|\x-\x_k|}
  e^{-|\x-\x_k|\sqrt{(\sigma+s)/D}} \,.
\end{equation}
By expanding $\hat{U}$ as $\x\to\x_j$, and enforcing that the limiting
behavior agrees with that required in \eqref{2:vk_2}, we determine
$\hat{B}_{j}$ as
\begin{equation}\label{2:hatbj}
  \hat{B}_{j} = - \frac{\hat{b}_j \sqrt{\sigma+s}}{4\pi D^{3/2}}  +
  \frac{1}{4\pi D}\sum_{\stackrel{k=1}{k\neq j}}^{N} \frac{\hat{b}_k}{|\x_j-\x_k|}
  e^{-|\x_j-\x_k|\sqrt{(\sigma+s)/D}} \,.
\end{equation}

To determine $B_{j}(t)={\mathcal L}^{-1}\left[\hat{B}_{j}\right]$, we
first use the inverse transform
\begin{equation}\label{2:free_3d}
  {\mathcal L}^{-1} \left[ \frac{e^{-|\x_j-\x_k|\sqrt{(\sigma+s)/D}}}
    {4\pi D|\x_j-\x_k|}\right] = G(a_{jk},t) \equiv
  \frac{e^{-\sigma t}e^{-a_{jk}^2/(4t)}}{\left(4\pi D t\right)^{3/2}}  \, \quad
  \mbox{where} \quad a_{jk}\equiv \frac{|\x_j-\x_k|}{\sqrt{D}} \,.
\end{equation}
From the convolution property of the Laplace transform and by
using the leading-order term from \eqref{2:bulk_cj}, we obtain
\begin{equation}\label{2:bj_new}
  B_{j}= -\frac{\kappa_j}{4\pi D^{3/2}} {\mathcal L}^{-1}\left[ \hat{u}_{1j}
    \sqrt{s+\sigma}\right] +
  \sum_{\stackrel{k=1}{k\neq j}}^{N} \kappa_k \int_{0}^{t} u_{k1}(\tau)
  G(a_{jk},t-\tau) \, d\tau + {\mathcal O}(\eps)\,.
\end{equation}

To invert the remaining term in \eqref{2:bj_new} we add and
subtract a term to get
\begin{equation}\label{2:invert_1}
  {\mathcal L}^{-1}\left[ \hat{u}_{1j}
    \sqrt{s+\sigma}\right]= {\mathcal L}^{-1}\left[ (s\hat{u}_{1j}
    -u_{1j}(0)) \frac{\sqrt{s+\sigma}}{s} \right] +  u_{1j}(0)
  {\mathcal L}^{-1}\left[ \frac{\sqrt{s+\sigma}}{s}\right] \,.
\end{equation}
By using ${\mathcal L}^{-1}\left[s\hat{u}_{1j}-u_{1j}(0)\right]=
u_{1j}^{\prime}(t)$ together with the inverse transform
\begin{equation}\label{3:lap2}
  {\mathcal L}^{-1}\left[  \frac{\sqrt{s+\sigma}}{s}\right] = g(t)\equiv
  \frac{e^{-\sigma t}}{\sqrt{\pi t}} + \sqrt{\sigma} \, \mathrm{erf}\left(
    \sqrt{\sigma t}\right) \,,
\end{equation}
we can invert \eqref{2:invert_1} using the convolution
formula. Upon substituting the resulting expression into
\eqref{2:bj_new} we obtain $B_{j}(t)$ as needed in \eqref{2:intra_full}.
Overall, this leads to an integro-differential system for the
intracellular concentrations $\v{u}_j(t)$ that is coupled to the
time-history of the source strengths $B_j(t)$. We summarize our result
in the following formal proposition.

\begin{prop}\label[prop]{prop:3ode} For $\eps\to 0$, and with the
  initial condition
  $U(\x, 0)=0$, the solution $U(\x,t)$ and $\v{u}_j(t)$, for
  $j\in\lbrace{1,\ldots,N\rbrace}$, to the cell-bulk model
  \eqref{DimLess_bulk} is approximated for $t\gg {\mathcal O}(\eps^2)$
  by the integro-differential system
  \begin{subequations}\label{2:reduced}
  \begin{align}
    \frac{d\v{u}_j}{dt} &= \v{F}_j(\v{u}_j)- \kappa_j u_{1j} \v{e}_1 +
                          \frac{\eps d_{1j}}{d_{2j}} \kappa_j B_{j}(t)
  \v{e}_1 \,,  \label{2:reduced_1} \\
    B_{j}(t) &= -\frac{\kappa_j u_{1j}(0)}{4\pi D^{3/2}} g(t) -
    \frac{\kappa_j}{4\pi D^{3/2}} \int_{0}^{t} u_{1j}^{\prime}(\tau)
                g(t-\tau)\, d\tau \nonumber \\
    & \qquad + \sum_{\stackrel{k=1}{k\neq j}}^{N} \kappa_k
    \int_{0}^{t} u_{k1}(\tau) G(a_{jk},t-\tau) \, d\tau \,,\label{2:reduced_2}
  \end{align}
  where $g(t)$ and the free-space Green's function $G(a_{jk},t)$ are defined by
  \begin{equation}\label{2:reduced_3}
    g(t) \equiv \frac{e^{-\sigma t}}{\sqrt{\pi t}} + \sqrt{\sigma}
    \mathrm{erf}\left(\sqrt{\sigma t}\right) \,, \qquad
   G(a_{jk},t) \equiv
   \frac{e^{-\sigma t}e^{-a_{jk}^2/(4t)}}{\left(4\pi D t\right)^{3/2}}   \,.
  \end{equation}
  Here $\mathrm{erf}(z)$ is the error function, while $\kappa_j$ and
  $a_{jk}$ are defined by
  \begin{equation}\label{2:reduced_4}
    a_{jk} \equiv \frac{|\x_j-\x_k|}{\sqrt{D}} \,, \qquad \kappa_j \equiv
    \frac{4\pi D d_{2j}}{d_{1j} + D}\,.
  \end{equation}
\end{subequations}
In terms of $B_{j}(t)$, the asymptotic solution in the vicinity of the
$j^{\mathrm{th}}$ cell, with $\rho=\eps^{-1}|\v{x}-\v{x}_j|={\mathcal O}(1)$, is
\begin{equation}\label{2:reduced_local}
  U \sim\frac{d_{2j}}{\eps \rho \left(d_{1j}+D\right)} \left(
    u_{1j}(t)- \frac{\eps d_{1j}}{d_{2j}} B_{j}(t)\right) +
  B_{j}(t) \,.
\end{equation}
Moreover, the asymptotic solution in the bulk region is
\begin{equation}\label{2:reduced_global}
  U(\x,t) \sim  \sum_{j=1}^{N} \kappa_j
  \int_{0}^{t} u_{1j}(\tau) \,
  G\left( \frac{|\v{x}-\v{x}_j|}{\sqrt{D}},t-\tau\right) \, d\tau \,.
\end{equation}
\end{prop}

\subsection{Matching to the Transient Solution}\label{sec:match}

We now show that the limiting behavior as $t\to 0^{+}$, but with
$t\gg {\mathcal O}(\eps^2)$, of the local solution near the
$j^{\mathrm{th}}$ cell in \eqref{2:inn_sol_full} agrees with the
far-field behavior of the transient solution as derived in
\eqref{3:v_final}.

To do so, we determine the short-time behavior of $B_{j}(t)$ by
letting $s\to \infty$ in \eqref{2:hatbj}. Upon using the leading-order
term in \eqref{2:bulk_cj}, and neglecting the negligible inter-cell
interactions in this limit, we obtain as $s\to\infty$ that
\begin{equation}\label{3:mat_1}
  \hat{B}_{j}(s) \sim -\frac{\kappa_j}{4\pi D^{3/2}}  \hat{u}_{1j}(s)
  \sqrt{s+\sigma}  \quad \mbox{with} \quad
  \hat{u}_{1j}(s) \sim \frac{u_{1j}(0)}{s} + \frac{u_{1j}^{\prime}(0)}{s^2} \,,
\end{equation}
where $\kappa_j$ is given in \eqref{2:reduced_4}. We readily calculate for
$s\to\infty$ that
\begin{equation}\label{3:mat_2}
  \hat{B}_{j}(s) \sim -\frac{d_{2j}}{\left(d_{1j}+D\right)\sqrt{D}}
  \left(\frac{u_{1j}(0)}{s} + \frac{u_{1j}^{\prime}(0)}{s^2}\right)
  \sqrt{s} \left(1 + \frac{\sigma}{2s} + \cdots \right) \,,
\end{equation}
which reduces to
\begin{equation}\label{3:mat_3}
  \hat{B}_{j}(s) \sim -\frac{d_{2j}}{\left(d_{1j}+D\right)\sqrt{D}}
  \left[ \frac{u_{1j}(0)}{\sqrt{s}} + \frac{1}{2 s^{3/2}}\left(
      u_{1j}(0)\sigma + 2 u_{1j}^{\prime}(0) \right) \right]\,.
\end{equation}
Upon inverting \eqref{3:mat_3}, we readily recover the result \eqref{3:v_bj}
derived from the long-time behavior of the transient solution. Moreover,
substituting \eqref{3:v_bj} into \eqref{2:ajbj_new}, we also recover
\eqref{3:v_aj}.

We conclude that the limiting behavior as $t\to 0^{+}$ of the solution
to the integro-differential system matches identically with the long-time
behavior of the transient solution near each cell.

\section{Steady-State and Linear Stability Analysis}\label{sec:ss}
In this section we construct steady-state solutions to
\eqref{DimLess_bulk} and we derive a globally coupled
eigenvalue problem (GCEP) that determines the linear stability of the
steady-states.

We first study the long-time dynamics of solutions to \eqref{2:reduced}
that tend to steady-state values as $t\to \infty$. To do so, we need
the following lemma:

\begin{lemma} Suppose that $u_{1j}(t)$ is bounded for all $t>0$ and
  that $\lim_{t\to\infty} u_{1j}(t)=u_{1j\infty}$ for
  $j\in\lbrace{1,\ldots,N\rbrace}$.  Then, with $\kappa_j$ as given in
  \eqref{2:reduced_4}, we have
\begin{equation}\label{4:long_time}
  \lim_{t\to\infty} B_{j}(t)= B_{j\infty} \equiv -\frac{\kappa_j\sqrt{\sigma}}
  {4\pi D^{3/2}} u_{1j\infty}  +
  \sum_{\stackrel{k=1}{k\neq j}}^{N} \frac{\kappa_k u_{k1\infty}}
  {4\pi D|\x_j-\x_k|} e^{-|\x_j-\x_k|\sqrt{\sigma/D}} \,.
\end{equation}

\end{lemma}

\begin{proof}
   To determine the long-time behavior of $B_{j}(t)$, we use the
   Tauberian theorem. From \eqref{2:free_3d} and
   \eqref{3:lap2} the Laplace transform of \eqref{2:reduced_2} is
   given by \eqref{2:hatbj} as
   \begin{equation}
     \hat{B}_{j}(s) = -\frac{\kappa_j \hat{u}_{1j}(s)}{4\pi D^{3/2}}
     \sqrt{s+\sigma} 
     + \sum_{\stackrel{k=1}{k\neq j}}^{N} \frac{\kappa_k\hat{u}_{k1}(s)}
      {4\pi D|\x_j-\x_k|} e^{-|\x_j-\x_k|\sqrt{(s+\sigma)/D}} \,, \label{4:lt_b}
  \end{equation}
  where the branch cut is chosen along $\mathrm{Re}(s)\leq -\sigma<0$
  and $\mathrm{Im}(s)=0$. Since $u_{1j}(t)$ is assumed to be bounded for
  $t>0$ with the limiting behavior $u_{1j}(t)\to u_{1j\infty}$ as
  $t\to \infty$, it follows that $\hat{u}_{1j}(s)$ is analytic in
  $\mathrm{Re}(s)>0$ and has the singular behavior
  $\hat{u}_{1j}(s)\sim {u_{1j\infty}/s}$ as $s\to 0$.  As a result,
  the singularity of $\hat{B}_{j}(s)$ with the largest real part is at
  $s=0$. To determine the long-time behavior of $B_{j}(t)$, we expand
  \eqref{4:lt_b} near this dominant singularity to obtain
  \begin{equation}\label{4:lt_b1}
    \hat{B}_{j}(s) \sim -\frac{\kappa_j \sqrt{\sigma}}{4\pi D^{3/2}}
    \frac{u_{1j\infty}}{s} 
     + \sum_{\stackrel{k=1}{k\neq j}}^{N} \frac{\kappa_k u_{k1\infty}}
     {4\pi s D|\x_j-\x_k|} e^{-|\x_j-\x_k|\sqrt{\sigma/D}}\,, \quad
     \mbox{as} \quad s\to 0 \,.
   \end{equation}
   The Tauberian result $B_{j\infty}=\lim_{s\to 0} s\hat{B}_{j}(s)$
   yields \eqref{4:long_time}. 
 \end{proof}

 From this lemma, we obtain the following result that characterizes
 those solutions of \eqref{2:reduced} tending to limiting values as
 $t\to \infty$.

 \begin{prop}\label{prop:limit} For \eqref{2:reduced}, suppose that
   $\v{u}_{j}(t)$ is bounded for $t>0$ and satisfies
   $\v{u}_{j}\to \v{u}_{j\infty}$ as $t\to \infty$ for each
   $j\in\lbrace{1,\ldots,N\rbrace}$. Then, for each
   $j\in \lbrace{1,\ldots,N\rbrace}$, we have that
   $\lim_{t\to\infty} B_{j}=B_{j\infty}$ and $\v{u}_{j\infty}$ satisfy
   the nonlinear algebraic system (NAS)
   \begin{subequations}\label{4:limit_sys}
     \begin{align}
      & \v{F}_j(\v{u}_{j\infty}) -\kappa_j u_{1j\infty} \v{e}_1 +
       \eps \frac{d_{2j}}{d_{1j}} \kappa_j B_{j\infty} \v{e}_1 = 0 \,,\\
      & B_{j\infty} = -\frac{\kappa_j}{D} R_{0}\, u_{1j\infty} -
       \sum_{\stackrel{k=1}{k\neq j}}^{N} \frac{\kappa_k}{D} u_{k1\infty}\,
        G_{0jk} \,, 
      \end{align}
    \end{subequations}
    where $u_{1j\infty}=\v{e}_1^T\v{u}_{j\infty}$. In
    \eqref{4:limit_sys}, $G_{0jk}$ and $R_0$ are defined by
    \begin{equation}\label{4:lim_sys_green}
      R_0 \equiv \frac{1}{4\pi} \sqrt{\frac{\sigma}{D}} \,, \qquad
      G_{0jk} \equiv -\frac{1}{4\pi|\x_j-\x_k|} e^{-|\x_j-\x_k|\sqrt{\sigma/D}}\,.
     \end{equation}
     We identify $G_{0jk}\equiv G_{0}(\x_j,\x_k)$, where $G_{0}(\x,\x_k)$
     is the steady-state free-space Green's function for $\x\in \R^3$,
     with regular part $R_{0}$, satisfying
     \begin{equation}\label{4:ge}
       \Delta G_0- \frac{\sigma}{D} G_0 =\delta(\x-\x_k) \,; \qquad
       G_0= -\frac{1}{4\pi|\x-\x_k|}+R_0 + o(1) \quad \mbox{as}\quad
       \x\to\x_k\,.
     \end{equation}
 \end{prop}

 In Appendix \ref{app:strong} we establish that the limiting NAS
 \eqref{4:limit_sys} is the steady-state system that results from an
 asymptotic analysis of \eqref{DimLess_bulk} for the time-independent
 case.
  
\subsection{Linear Stability Analysis}\label{sec:stab}

Next, we derive a globally coupled eigenvalue problem (GCEP) for the
linearization around the steady-state solution of
\eqref{DimLess_bulk}. To determine the linear stability of the
steady-states, we let $\v{u}_{js}=\left(u_{1js},\ldots,u_{1Ns}\right)^T$ and
introduce the eigenperturbation
\begin{equation}\label{5:pert}
  U=U_{s}(\v{x}) + e^{\lambda t} \Phi(\v{x}) \,, \quad
  \v{u}_j = \v{u}_{js} + e^{\lambda t} \v{\zeta}_j \,, \quad
  j\in \lbrace{1,\ldots,N\rbrace}\,,
\end{equation}
into \eqref{DimLess_bulk}. By linearizing the resulting system we obtain the
eigenvalue problem
\begin{subequations}\label{5:eig_prob}
\begin{align}
 &  \Delta \Phi - \frac{(\lambda+\sigma)}{D} \Phi =0\,, \quad
   \x \in \R^3 \setminus \cup_{j=1}^{N}\,\Omega_{\varepsilon_j}\,; \quad
   \Phi(\v{x})\to 0 \,\,\,  \mbox{as}\,\,\, |\v{x}|\to \infty \,,
      \label{5:eig_prob_1} \\
  &  \varepsilon  D\, \partial_n \Phi  = d_{1j} \, \Phi - \frac{d_{2j}}{\eps}
    \, \zeta_{1j}\,,
    \quad \x \in \partial \Omega_{\varepsilon_j}\,,\qquad j \in
    \lbrace{1, \ldots, N\rbrace}\,,  \label{5:eig_prob_2}
\end{align}
which is coupled to the linearized intracellular dynamics within the
$j^{\text{th}}$ cell by
\begin{align}\label{5:eig_prob_3}
  \left(\lambda I - J_j\right) \v{\zeta}_j & =  
   \frac{\v{e}_1}{\eps}\int_{\partial \Omega_{\varepsilon_j}}\, \left( d_{1j} \,
     \Phi - \frac{d_{2j}}{\eps} \, \zeta_{1j}\right) \, \text{d}S_{\x}\,,
                \quad j \in \lbrace{1,\ldots, N\rbrace} \,.
\end{align}
\end{subequations}
Here $I$ is the $m\times m$ identity matrix and
$J_j\equiv {\partial \v{F}_j/\partial \v{u}_j}$ is the $m\times m$
dimensional Jacobian matrix of the intracellular kinetics evaluated at the
steady-state $\v{u}_j=\v{u}_{js}$.

For $\eps\to 0$, we now analyze \eqref{5:eig_prob} using strong
localized perturbation theory \cite{ward2018spots}. In the inner
region near the $j^{\mathrm{th}}$ cell we let $\rho=|\y|$ with
$\y= \varepsilon^{-1}(\x-\x_j)$ and to leading-order we get that
\begin{equation}\label{5:inn_eig}
  \Delta_{\y} \Phi = 0 \,, \quad \mbox{for} \quad \rho>1 \,; \qquad
  D \partial_\rho \Phi= d_{1j}\Phi - \frac{d_{2j}}{\eps} \zeta_{1j} \,, \quad
  \mbox{on} \quad \rho=1\,,
\end{equation}
which has the radially symmetric solution
\begin{equation}\label{5:inn_sol}
  \Phi =\frac{c_j}{\eps \rho} -\frac{c_j}{\eps}\left(1 + \frac{D}{d_{1j}}
    \right) + \frac{d_{2j}}{\eps d_{1j}} \zeta_{1j} \,,
\end{equation}
where $c_j$ is to be found. In terms of $c_j$, we get
$\eps^{-1}\int_{\Omega_{\eps_j}} \left(d_{1j}\Phi - d_{2j} \eps^{-1}
  \zeta_{1j}\right)\, \text{d}S_{\v{x}}= -4\pi D c_j$ so that from
\eqref{5:eig_prob_3} we conclude that
\begin{equation}\label{5:intra_1}
  \left(\lambda I - J_j\right) \v{\zeta}_j = -4\pi D c_j \v{e}_1 \,,
  \quad j\in \lbrace{1,\ldots, N\rbrace}\,.
\end{equation}

Next, we write \eqref{5:inn_sol} in terms of the outer variable, which
provides an asymptotic matching condition for the eigenperturbation
$\Phi$ in the outer bulk region. In this way, we obtain that $\Phi$ satisfies
\begin{subequations}\label{5:outer}
  \begin{align}
 & \Delta \Phi - \frac{(\lambda+\sigma)}{D} \Phi = 0\,, \quad
   \x \in \R^3 \setminus \cup_{j=1}^{N}\,\Omega_{\varepsilon_j}\,;
   \quad \Phi(\v{x})\to 0 \,\,\,  \mbox{as}\,\,\, |\v{x}|\to \infty \,,
      \label{5:outer_1} \\
    &   \Phi \sim
      \frac{c_j}{|\x-\x_j|} -\frac{c_j}{\eps}\left(1 + \frac{D}{d_{1j}}
      \right) + \frac{d_{2j}}{\eps d_{1j}} \zeta_{1j} \,,
      \quad \mbox{as} \quad \v{x}\to \v{x}_j \,, \quad
      j\in \lbrace{1,\ldots,N\rbrace} \,.
      \label{5:outer_2}
  \end{align}
\end{subequations}
The pre-specification of the regular part of this singularity condition
\eqref{5:outer_2} will provide a linear homogeneous system for
$\v{c}=\left(c_1,\ldots,c_N\right)^{T}$. 

We represent the solution to \eqref{5:outer} as
\begin{equation}\label{5:eigbulk_sol}
  \Phi(\v{x}) = -4\pi \sum_{k=1}^{N} c_k G_{\lambda}(\x,\x_k)\,,
\end{equation}
where $G_{\lambda}(\x,\x_k)$ is the eigenvalue-dependent Green's function
for $\x\in \R^3$ satisfying
\begin{equation}\label{4:ge_lambda}
  \Delta G_\lambda- \frac{(\sigma+\lambda)}{D} G_\lambda =
  \delta(\x-\x_k) \,; \quad
  G_\lambda\sim -\frac{1}{4\pi|\x-\x_k|}+R_\lambda+o(1)
  \quad \mbox{as}\quad \x\to\x_k\,.
\end{equation}
The solution $G_{\lambda}$ to \eqref{4:ge_lambda} and its regular part
$R_{\lambda}$ are given explicitly by
\begin{equation}\label{4:ge_lambda_ex}
  G_{\lambda}(\x,\x_k) = -\frac{e^{-|\x-\x_k|\sqrt{(\sigma+\lambda)/D}}}
  {4\pi|\x-\x_k|}\,, \qquad R_{\lambda} =
  \frac{1}{4\pi} \sqrt{ \frac{\sigma+\lambda}{D}} \,.
\end{equation}
We identify that when $\lambda=0$ we recover $G_{0}(\x,\x_k)$, which satisfies
\eqref{4:ge}, which was key to the construction of steady-state solutions.

Upon enforcing that the limiting behavior as $\x\to\x_j$ for $\Phi$ in
\eqref{5:eigbulk_sol} agrees with that required in the matching
condition \eqref{5:outer_2}, we obtain the linear system
\begin{equation}\label{5:lin_cj}
  c_j\left(1 + \frac{D}{d_{1j}}\right) - 4\pi \eps c_j R_\lambda -
  4\pi \eps \sum_{\stackrel{k=1}{k\neq j}}^{N} c_{k} G_{\lambda}(\x_j,\x_k)
  = \frac{d_{2j}}{d_{1j}} \zeta_{1j}\,, \quad j\in \lbrace{1,\ldots,N\rbrace} \,,
\end{equation}
where $\zeta_{1j}$ is given in terms of $c_j$ by \eqref{5:intra_1}.
When $\lambda$ is not an eigenvalue of $J_j$, we use \eqref{5:intra_1}
to calculate
$\zeta_{1j}=-4\pi D c_j \v{e}_1^T \left(\lambda I-J_j\right)^{-1}
\v{e}_1$. Finally, by substituting this expression into \eqref{5:lin_cj} and
writing the resulting homogeneous linear system for
$\v{c}=(c_1,\ldots,c_N)^T$ in matrix form, we conclude that the
discrete eigenvalues $\lambda$ of the linearization around a
steady-state solution of \eqref{DimLess_bulk} are obtained from a
nonlinear matrix eigenvalue problem, which we refer to as the
GCEP. Our result is as follows:
  
\begin{prop}\label{prop:stab} For $\eps\to 0$, the discrete
  eigenvalues $\lambda$ associated with the linearization around a
  steady-state solution to \eqref{DimLess_bulk}, for which
  $\mathrm{det}(\lambda I-J_j)\neq 0$ for all
  $j\in\lbrace{1,\ldots,N\rbrace}$, are the set of values
  \begin{subequations}\label{5:gcep}
    \begin{equation}\label{TransDent}
   \Lambda({\mathcal M}) \equiv \lbrace{ \lambda \, \vert \,
    \mathrm{det}\mathcal{M}(\lambda)=0 \rbrace} \,.
\end{equation}
Here the $N\times N$ dimensional matrix ${\mathcal M}(\lambda)$ is
defined by
\begin{equation}\label{5:gcep_2}
  \mathcal{M}(\lambda) \equiv I + D P_1 -4\pi \eps \mathcal{G}_{\lambda}  +
  4\pi D P_2 \mathcal{H}(\lambda) \,.
\end{equation}
\end{subequations}
In \eqref{5:gcep_2}, $P_1$, $P_2$ and ${\mathcal H}$ are the
diagonal matrices
\begin{subequations}
\begin{align}
  P_1 &\equiv \mathrm{diag}\Big{(} \frac{1}{d_{11}}, \ldots,
        \frac{1}{d_{1N}} \Big{)}
        \,, \qquad
        P_2\equiv \mathrm{diag}\left(\frac{d_{21}}{d_{11}},\ldots,
        \frac{d_{2N}}{d_{1N}} \right) \,, \label{5:p1p2} \\
   \mathcal{H}(\lambda) &\equiv \mathrm{diag}\left( \mathcal{H}_1,\ldots,
  \mathcal{H}_N\right) \quad \mbox{where} \quad {\mathcal H}_j(\lambda) \equiv
  \v{e}_1^T \left(\lambda I - J_j  \right)^{-1}\v{e}_1 \,.\label{5:kjac}
\end{align}
\end{subequations}
In addition, ${\mathcal G}_{\lambda}$ is the eigenvalue-dependent
Green's matrix with matrix entries
\begin{equation}\label{GreenMat_lambda}
  ({\mathcal G}_{\lambda})_{ij} = ({\mathcal G}_{\lambda})_{ji} \equiv
  -\frac{e^{-|\x_i-\x_j|\sqrt{(\sigma+\lambda)/D}}}{4\pi|\x_i-\x_j|} \,, \quad
      i\neq j \,; \quad
    ({\mathcal G}_{\lambda})_{jj} = R_{\lambda}= \frac{1}{4\pi}
        \sqrt{ \frac{\sigma+\lambda}{D}} \,,
      \end{equation}
      where we take the principal branch of the square root.  The
      steady-state is linearly stable if, whenever
      $\lambda\in \Lambda({\mathcal M})$, we have
      $\mathrm{Re}(\lambda)<0$.  The vector $\v{c}=(c_1,\ldots,c_N)^T$,
      satisfying ${\mathcal M}(\lambda) \v{c}=\v{0}$, is nontrivial
      only when $\lambda\in \Lambda({\mathcal M})$.
\end{prop}

If $\lambda\in \Lambda({\mathcal M})$ is complex-valued, the matrix
${\mathcal G}_{\lambda}$, and consequently the GCEP matrix
${\mathcal M}(\lambda)$, is complex-valued and symmetric but not
Hermitian. As such, the appropriate inner product is the transpose and not
the Hermitian conjugate. We now discuss the information encoded by the
complex-valued nullvector $\v{c}=(c_1,\ldots,c_N)^{T}$, which we
normalize by $\v{c}^{T}\v{c}=1$. We use the steady-state solution
$U_{s}(\x)$ and the eigenperturbation \eqref{5:pert} to calculate both
the diffusive flux into the $j^{\mathrm{th}}$ cell and the corresponding
permeable intracellular species, for each
$j\in\lbrace{1,\ldots,N\rbrace}$, as
\begin{subequations}\label{nstabform:c}
\begin{align}
  D \partial_{\rho} U \vert_{\rho=1} &\sim -\frac{D}{\eps} \left(
     \frac{d_{2j}}{d_{1j}+D} u_{1js} +
   \sum_{\lambda\in\Lambda({\mathcal M})} c_j  e^{\lambda t}\right) \,,
    \label{nstabform:c1}  \\
 u_{1j} &\sim u_{1js} -4\pi D 
  \sum_{\lambda\in\Lambda({\mathcal M})} {\mathcal H}_{j}(\lambda) c_j
          e^{\lambda t} \,, \quad \mbox{where} \quad
{\mathcal H}_j(\lambda)=\v{e}_1^T \left(I-\lambda J_j\right)^{-1}
\v{e}_1\,.   \label{nstabform:c2}
\end{align}
\end{subequations}
From \eqref{nstabform:c} we conclude that the components of the
normalized eigenvector $\v{c}$ determine the relative magnitude of the
spatial gradient of the bulk signal near the cell boundaries.  If
$\lambda\in \Lambda({\mathcal M})$ is complex-valued, the real and
imaginary parts of $c_j$ in \eqref{nstabform:c2} encode both the
relative magnitude and phase shift of intracellular oscillations
within the cell population near the steady-state.

\subsection{Perturbation Theory: Linear Stability Analysis}
\label{sec:eig_pert}

In contrast to the 2-D case studied in \cite{pelzward2D}, for our 3-D
cell-bulk model we observe from the main results \eqref{4:ss_final}
and \eqref{5:gcep} of our steady-state and linear stability analyses,
that there is only an ${\mathcal O}(\eps)$ effect due to intercellular
coupling, which is mediated by the Green's matrices that depend on the cell
locations. As such, it is worthwhile to develop a perturbation
analysis to approximate the eigenvalues from the GCEP \eqref{5:gcep}, under
the assumption that the steady-state NAS \eqref{4:ss_final} has a
solution for any $\eps\ll 1$. In the analysis below we assume that the
Jacobian of the intracellular kinetics $J_j$ is evaluated at this
steady-state that depends weakly on $\eps$.

We first consider the $\eps=0$ problem for the GCEP
\eqref{5:gcep}. For this limiting problem, we seek all values
$\lambda_0$ for which
$\mathrm{det}\left({\mathcal M}_0(\lambda_0)\right)=0$ where
${\mathcal M}_0(\lambda_0)\equiv I+DP_1 + 4\pi D P_2 {\mathcal H}(\lambda_0)$.
Since ${\mathcal M}_0$ is a diagonal matrix, the spectrum of this
limiting GCEP is
\begin{equation}\label{per:spec0}
  \Lambda_{0}({\mathcal M}_0) =
  \lbrace{ \lambda_0 \, \vert \,\,
    \v{e}_1^T (\lambda_0 I - J_j)^{-1}\v{e}_1= -\frac{1}{\kappa_j}\,,
    \,\,\, j\in \lbrace{1,\ldots,N\rbrace} \rbrace} \,,
\end{equation}
where $\kappa_j={4\pi D d_{2j}/(d_{1j}+D)}$. An equivalent formulation
of \eqref{per:spec0} can be obtained by the following
well-known matrix determinant lemma (see Lemma 1.1 \cite{ding}).

\begin{lemma} Let $A$ be an invertible $m\times m$ matrix and let
  $\v{a}\in \R^m$ and $\v{b}\in \R^m$ be two column vectors. Then,
\begin{equation}\label{mat:lemma}
   \mathrm{det}\left(A+ \v{a}\v{b}^T\right) = \left(1+ \v{b}^T A^{-1}
 \v{a} \right)\mathrm{det}(A) \,.
\end{equation}
As such, $(A+ \v{a} \v{b}^T)\v{\phi}=\v{0}$ has a nontrivial
solution if and only if $\v{b}^T A^{-1} \v{a}=-1$.
\end{lemma}

By applying this lemma to \eqref{per:spec0} we obtain the following
equivalent re-formulation of the set $\Lambda_{0}({\mathcal M}_0)$.

\begin{lemma} Suppose that $\mathrm{det}(\lambda I -J_j)\neq 0$ for each
  $j\in \lbrace{1,\ldots,N\rbrace}$. Then, we have
\begin{equation}\label{pert:mlemma}
    \mathrm{det}\left( \lambda I - J_j + \kappa_j\v{e}_1 \v{e}_1^T \right) =0\,,
\end{equation}
if and only if $\lambda$ satisfies
\begin{equation}\label{pert:orig}
  \v{e}_1^T (\lambda I - J_j)^{-1}\v{e}_1= -\frac{1}{\kappa_j}\,.
\end{equation}
Therefore, the set in \eqref{per:spec0} is equivalent to
\begin{equation}\label{per:spec0_new}
  \Lambda_{0}({\mathcal M}_0) =
 \Big{\lbrace{} \lambda_0 \, \vert \,\,
    \mathrm{det}\left( \lambda_0 I - J_j + \kappa_j\,
      \v{e}_1\v{e}_1^T \right) =0\,, \quad j\in \lbrace{1,\ldots,N \rbrace}
    \, \Big{\rbrace}}\,.
\end{equation}
\end{lemma}

\begin{proof}
  The result follows by setting $A=\lambda I-J_j$,
  $\v{a}=\v{e}_1$ and $\v{b}=\kappa_j \v{e}_1$ in \eqref{mat:lemma}.
\end{proof}

Next, for any $\lambda_0\in \Lambda_0({\mathcal M}_0)$ we will
determine a two-term asymptotic expansion of the form
$\lambda=\lambda_0+\eps\lambda_1 +\cdots$, for which we can ensure
that $\lambda\in \Lambda({\mathcal M})$. For a particular choice of
reaction kinetics, in \S \ref{sec:selkov} this will allow us to analyze
how stability thresholds of the uncoupled intracellular kinetics are
perturbed by weak intercellular coupling.

We begin by writing the nonlinear matrix eigenvalue problem
${\mathcal M}(\lambda)\v{c}=\v{0}$
in \eqref{5:gcep_2} as
\begin{subequations}\label{pert:all_main}
\begin{equation}\label{pert:main}
  \left( {\mathcal M_0}(\lambda)-4\pi \eps {\mathcal G}_{\lambda}\right)
  \v{c}=\v{0} \,, \quad \mbox{with} \quad
  {\mathcal M}_0(\lambda)=\mathrm{diag}\left(m_{1}(\lambda),\ldots,m_{N}(\lambda)
    \right)\,,
\end{equation}
where
\begin{equation}\label{pert:m0}
  \quad m_{j}(\lambda)\equiv 1+\frac{D}{d_{1j}} +
  \frac{4\pi D d_{2j}}{d_{1j}} \v{e}_1^{T}(\lambda I - J_j)^{-1}\v{e}_1\,.
\end{equation}
\end{subequations}
The perturbation expansion that we now develop depends on the
dimension of the null space of ${\mathcal M}_{0}(\lambda_0)$ when
$\lambda_0\in \Lambda_0({\mathcal M}_0)$.

We first suppose that $m_J(\lambda_0)=0$ and
$m_{J}^{\prime}(\lambda_0)\neq 0$ for some unique
$J\in\lbrace{1,\ldots,N\rbrace}$, so that $m_j(\lambda_0)\neq 0$ for
$j\neq J$. This implies that $\lambda=\lambda_0$ is a simple root of
$\mathrm{det}({\mathcal M}_{0}(\lambda))=0$ and that
${\mathcal M}_{0}(\lambda_0)\v{c}_0=0$ with $\v{c}_0=\v{e}_J$, where
$\v{e}_J$ is the unit basis vector in the $J^{\mathrm{th}}$
direction. In \eqref{pert:all_main} we expand
\begin{equation}\label{pert:exp1}
  \lambda=\lambda_0 + \eps \lambda_1 + \ldots \,, \qquad
  \v{c}=\v{e}_J + \eps \v{c}_1 + \ldots\,,
\end{equation}
where $\v{c}^{T}\v{c}=1$ implies that $\v{c}_1^{T}\v{e}_J=0$. Upon
collecting powers of $\eps$, we obtain
\begin{equation}\label{pert:exp1_a}
  {\mathcal M}_{0}(\lambda_0) \v{c}_1= \v{r} \,, \quad \mbox{where} \quad
  \v{r}\equiv \left( 4\pi {\mathcal G}_{\lambda_0} -
    {\mathcal M}_{0}^{\prime}(\lambda_0) \lambda_{1} \right)\v{e}_J\,.
\end{equation}
Since ${\mathcal M}_{0}(\lambda_0)={\mathcal M}_0^T(\lambda_0)$, the
solvability condition for \eqref{pert:exp1_a} is that
$\v{e}_J^{T}\v{r}=0$, which yields
\begin{equation}\label{part:exp1_sol}
  \lambda_{1}=\frac{4\pi \v{e}_J^{T} {\mathcal G}_{\lambda_0} \v{e}_J}
  {\v{e}_J^T {\mathcal M}_{0}^{\prime}(\lambda_0)\v{e}_J} =
  \frac{4\pi  R_{\lambda_0}}{m_{J}^{\prime}(\lambda_0)} \,, \quad
  \mbox{where} \quad m_{J}^{\prime}(\lambda_0)= -\frac{4\pi d_{2j}}{d_{1j}}
  \v{e}_1^T \left(\lambda_0 I - J_j\right)^{-2} \v{e}_1 \,.
\end{equation}
Under this condition, there is a unique solution $\v{c}_1$ to
\eqref{part:exp1_sol} for which $\v{c}_1^{T}\v{e}_J=0$. 

We remark that if $\lambda_0\in \mathbb{C}$, then we must have
$m_J(\overline{\lambda}_0)=0$. By using the explicit results in
\eqref{GreenMat_lambda}, we readily verify that
${\mathcal G}_{\overline{\lambda_0}}=\overline{{\mathcal G}_{\lambda_0}}$. As
a result, the perturbation of the conjugate eigenvalue is simply
$\lambda=\overline{\lambda_0}+\eps\overline{\lambda}_1+{\mathcal O}(\eps^2)$.

Next, we consider the case of identical cells where $d_{1j}=d_{1c}$,
$d_{2j}=d_{2c}$ and identical reaction-kinetics $\v{F}_j=\v{F}_c$, for
which $m_j(\lambda)=m_{c}(\lambda)$ for
$j\in \lbrace{1,\ldots,N\rbrace}$. Moreover, we assume that
$\lambda_0$ satisfies $m_c(\lambda_0)=0$, with
$m_{c}^{\prime}(\lambda_0)\neq 0$.  As shown below in \S
\ref{sec:selkov}, this case is relevant to the situation where the
spatial locations $\x_j$ of the $N$ cells are such that
${\mathcal G}_{\lambda}$ has a constant row sum.  For this case, we
have
\begin{equation}
  {\mathcal M}_0 = m_{c}(\lambda) I  \quad \mbox{where} \quad
  m_{c}(\lambda) \equiv  1+\frac{D}{d_{1c}} +
  \frac{4\pi d_{2c}}{d_{1c}} \v{e}_1^{T}(\lambda I - J_c)^{-1}\v{e}_1\,,
\end{equation}
where $J_c$ is the common Jacobian of the intracellular kinetics. 
Therefore, we have $\v{c}_0=\mathrm{span}(\v{e}_1,\ldots,\v{e}_N)$, and
we readily conclude that $\lambda_1$ and $\v{c}_1$ satisfy
\begin{equation}\label{part:exp2_sol}
  \left( 4\pi {\mathcal G}_{\lambda_0} - m_{c}^{\prime}(\lambda_0)
    \lambda_1 I \right)\v{c}_1=\v{0} \quad \Rightarrow \quad
    \mathrm{det}\left( 4\pi {\mathcal G}_{\lambda_0} - m_{c}^{\prime}(\lambda_0)
    \lambda_1 I\right)=0 \,.
\end{equation}
This gives $N$ different values for $\lambda_1$ (counting
multiplicity).  If $\lambda_0\in \mathbb{C}$, it follows that
$m_{c}(\overline{\lambda}_0)=0$. Since
${\mathcal G}_{\overline{\lambda_0}}= \overline{{\mathcal
    G}_{\lambda_0}}$, any conjugate eigenvalue has the perturbation
$\lambda\sim\overline{\lambda_0}+\eps\overline{\lambda}_1$.

Our key conclusion is that the spectrum of the coupled cell-bulk
problem around the steady-state is ${\mathcal O}(\eps)$ close to the
spectrum associated with the linearization of the leading-order
intracellular problems obtained by setting $\eps=0$ in
\eqref{2:reduced_1}, which is given by
\begin{equation}\label{pert:intra}
  \frac{d\v{u}_j}{dt} = \v{F}_j(\v{u}_j)- \kappa_j u_{1j} \v{e}_1 \,, \quad
 \mbox{with} \quad \kappa_j= \frac{4\pi d_{2j}D}{d_{1j}+D} \,, \quad
  j\in \lbrace{1,\ldots,N\rbrace}\,.
\end{equation}
We refer to \eqref{pert:intra} as the {\em ``single cell
    problem''}, since the cell is still embedded in the diffusion
  field in $\R^3$, but it does not receive any diffusive signal from the
  other cells in the population. For the steady-state of
\eqref{pert:intra}, the perturbation results in \eqref{part:exp1_sol}
and \eqref{part:exp2_sol} show how the eigenvalue of the linearization
is altered by ${\mathcal O}(\eps)$ owing to the cell-bulk coupling. In
particular, if the kinetic parameters are ${\mathcal O}(\eps)$ close
to an instability threshold of \eqref{pert:intra}, then the cell-bulk
coupling may effectively destabilize the steady-state. This behavior
is in contrast to the 2-D case studied in \cite{pelzward2D} where
there is a much stronger ${\mathcal O}(1)$ effect for the coupling of
the intracellular dynamics via the bulk medium.

\section{Linear Stability Theory: Sel'kov Kinetics}\label{sec:selkov}

We now implement our linear stability analysis of
steady-states given in \S \ref{sec:eig_pert} for the two-component
reaction kinetics used in \cite{MNS} to model chemical
oscillations. This conditional oscillator kinetics is a modification
of the Sel'kov kinetics \cite{Selkov} used to model glycolysis
oscillations \cite{nandy1998}. Labeling
$\v{F}_j(u_{1j},u_{2j})\equiv (F_{1j},F_{2j})^T$, the intracellular
kinetics are
\begin{equation}\label{isolated:selkov:ode}
  F_{1j} = \alpha_j u_{2j}+ u_{2j} u_{1j}^2- u_{1j}\,,
  \qquad F_{2j} =\zeta_j
  \left[\mu_j-\left(\alpha_j u_{2j}+u_{2j} u_{1j}^2\right)\right] \,,
\end{equation}
where the reaction-kinetic parameters $\alpha_j$, $\mu_j$, and $\zeta_j$
are all positive.

With \eqref{isolated:selkov:ode} the steady-state intracellular
species satisfies \eqref{4:ss_final}, which reduces to a linear algebraic
system for the permeable intracellular species $\v{u}_{1s}$ given by
\begin{subequations}\label{selkov:ss}
\begin{equation}\label{selkov:ss1}
  \left( I  +   {\mathcal K} + \frac{\eps}{D} {\mathcal P} {\mathcal K}
    {\mathcal G}_0 {\mathcal K}\right) \v{u}_{1s} = \v{\mu}\,,
\quad \mbox{with}\quad  u_{2js}=\frac{\mu_j}{\alpha_j+ u_{1js}^2} \,,
\end{equation}
where $\v{\mu}=(\mu_1,\ldots,\mu_N)^T$ and the diagonal matrices
${\mathcal K}$ and ${\mathcal P}$ are defined by
\begin{equation}\label{selkov:ss2}
  {\mathcal P} \equiv \mathrm{diag}\Big{(} \frac{d_{11}}{d_{21}}, \ldots,
        \frac{d_{1N}}{d_{2N}} \Big{)}
        \,, \qquad
        {\mathcal K}\equiv \mathrm{diag}\left( \kappa_1,\ldots,
        \kappa_N \right) \,.
\end{equation}
\end{subequations}
In terms of the solution to \eqref{selkov:ss},
$\v{B}_{s}=(B_{1s},\ldots,B_{Ns})^T$ is given by \eqref{4:b0eq}.
For ${\eps/D}\ll 1$, a two-term approximation for the solution to
\eqref{selkov:ss1} is
\begin{equation}\label{selkov:ss1_approx}
  \v{u}_{1s}= \left(I+{\mathcal K}\right)^{-1}\mu -
  \frac{\eps}{D} {\mathcal P} {\mathcal K} \left(I+{\mathcal K}\right)^{-1}
  {\mathcal G}_0 {\mathcal K} \left(I+{\mathcal K}\right)^{-1} \v{\mu}
  + {\mathcal O}\left(\left[{\eps/D}\right]^2\right)\,.
\end{equation}

To determine the set $\Lambda_0({\mathcal M}_0)$ in \eqref{per:spec0_new}
we first label $\tilde{J}_j$ by
\begin{equation}\label{selkov:tJj}
   \tilde{J}_j \equiv J_j - \kappa_j \v{e}_1\v{e}_1^{T} \,.
\end{equation}
Upon linearizing \eqref{isolated:selkov:ode} around the steady-state
we readily calculate that
\begin{equation}\label{selkov:m0}
  \Lambda_{0}({\mathcal M}_0) =
 \Big{\lbrace{} \lambda_0 \, \vert \,\,
   \lambda_0^2 - \mathrm{tr}(\tilde{J}_j) \lambda_0 +
   \mathrm{det}(\tilde{J}_j) =0\,, \quad j\in \lbrace{1,\ldots,N \rbrace}
    \, \Big{\rbrace}}\,,
\end{equation}
where, in terms of the solution $u_{1js}$ to \eqref{selkov:ss1} we have
\begin{equation}\label{selkov:det_trace}
  \begin{split}
&    \mathrm{det}(\tilde{J}_j)= \zeta_j (1+\kappa_j) \left(\alpha_j
  + u_{1js}^2\right)>0 \,,\\
& \mathrm{tr}(\tilde{J}_j)= -1-\kappa_j + \frac{2\mu_ju_{1js}}
{\alpha_j+u_{1js}^2} - \zeta_j \left(\alpha_j + u_{1js}^2\right)\,.
\end{split}
\end{equation}

\begin{figure}[!h]
\centering
\includegraphics[width=88mm]{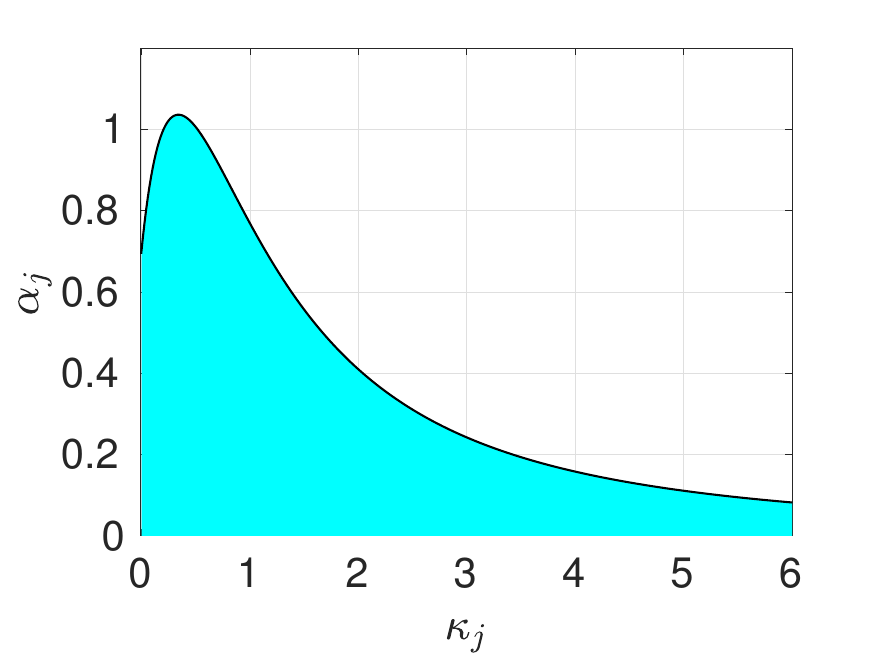}
\caption{Cyan-shaded region where limit-cycle oscillations occur in
  the $\alpha_j$ versus $\kappa_j={4\pi d_{2j} D/(D+d_{1j})}$ plane
  for the single cell problem \eqref{pert:intra} with boundary
  efflux $d_{2j}>0$ under Sel'kov kinetics
  \eqref{isolated:selkov:ode}.  Parameters are $\zeta_j=0.15$ and
  $\mu_j=2.0$. The HB threshold in \eqref{selkov:isolated_HB} for the
  linearization of the steady-state for the single cell problem is
  the boundary of this region. In the unshaded region the steady-state
  for \eqref{pert:intra} is linearly stable.}
\label{fig:selkov}
\end{figure}

By using the leading-order solution $u_{1js}={\mu_j/(1+\kappa_j)}$ to
\eqref{selkov:ss1_approx}, we can find the
region of instability of the steady-state of the {\em single cell
  problem} \eqref{pert:intra}. From \eqref{selkov:det_trace} we
conclude that the steady-state for the $j^{\mathrm{th}}$ single cell
is linearly stable, i.e.  $\mathrm{Re}(\lambda_0)<0$, only when
$\mathrm{tr}(\tilde{J}_{j})<0$. Since $\mathrm{det}(\tilde{J}_j)>0$,
instability can never occur from a zero-eigenvalue crossing and can
only occur by a Hopf bifurcation. By setting $\mathrm{tr}(\tilde{J}_{j})=0$
in \eqref{selkov:det_trace} and solving for $\alpha_j$, we obtain that
$\mathrm{Re}(\lambda_0)<0$ provided that
\begin{equation}\label{selkov:isolated_HB}
  \alpha_j > \alpha_{jH} \equiv - u_{1js}^2
  +\frac{1}{2\zeta_j}\left[- (1+\kappa_j) + \sqrt{
      (1+\kappa_j)^2 + 8\mu_j\zeta_j u_{1js} }\right] \,,
\end{equation}
where $\kappa_j={4\pi d_{2j}D/(d_{1j}+D)}$ and
$u_{1js}={\mu_j/(1+\kappa_j)}$. We identify $\alpha_{jH}$ as the Hopf
bifurcation (HB) threshold for the $j^{\mathrm{th}}$ single cell. The
corresponding HB eigenvalue is $\lambda_0=i\lambda_{0Ij}$, where 
\begin{equation}\label{selkov:lambda0}
  \lambda_{0Ij} = \sqrt{\mathrm{det}(\tilde{J}_j)} =\sqrt{
  \zeta_j (1+\kappa_j) \left(\alpha_j + u_{1js}^2\right)}\,.
\end{equation}
From the Poincar\'e-Bendixson theorem, whenever the steady-state is
unstable, the single cell exhibits limit-cycle oscillations. The
region in the $\alpha_j$ versus $\kappa_j$ plane where intracellular
oscillations occur for an single cell with boundary efflux $d_{2j}>0$
is shown in Fig.~\ref{fig:selkov} for a specific parameter
set. We find that $\alpha_{jH}$ first increases as $\kappa_j$
increases, but then decays monotonically with $\kappa_j$ beyond the
peak value.

\subsection{Two Non-Identical Cells}\label{sec:exper_two}

In this subsection we consider two cells that have the same reaction
kinetic parameters, labeled by $\alpha_c=\alpha_j$, $\mu_c=\mu_j$, and
$\zeta_c=\zeta_j$ for $j\in\{1,2\}$, but have different influx and efflux
permeabilities $d_{1j}$ and $d_{2j}$ for $j\in\{1,2\}$. Without loss of
generality we fix the inter-cell distance as $|\v{x}_1-\v{x}_2|=2$
since we can take the diffusive length scale $\sqrt{{D/\sigma}}$ in
the bulk as the parameter that measures the strength of intercellular
interactions when the cell radius $\varepsilon$ is small enough.

For this two-cell case, the steady-state $\v{u}_{1s}$ can be obtained
from the $2\times 2$ linear system \eqref{selkov:ss}, which is given explicitly
by $\v{u}_{1s}=\mu_c {\mathcal A}^{-1}\v{e}$, where
\begin{equation}\label{sec:exper_two_A}
  {\mathcal A} \equiv
  \left( {\begin{array}{cc}
            1+\kappa_1+\frac{\eps d_{11}}{D d_{21}}\kappa_1^2 R_0 &
             \frac{\eps d_{11}}{D d_{21}} \kappa_1 \kappa_2 G_{012} \\
            \frac{\eps d_{12}}{D d_{22}} \kappa_1 \kappa_2 G_{012} &
            1+\kappa_2+\frac{\eps d_{12}}{D d_{22}}\kappa_2^2 R_0 \\
  \end{array} } \right)\,,
\end{equation}  
where $\kappa_j={4\pi d_{2j} D/(d_{1j}+D)}$,
$R_0=(4\pi)^{-1} \sqrt{{\sigma/D}}$, and
$G_{012}=-(8\pi)^{-1}e^{-2\sqrt{{\sigma/D}}}$.  In terms of
$\v{u}_{1s}=(u_{11s},u_{21s})^{T}$, the HB threshold for the
$j^{\mathrm{th}}$ cell in terms of the reaction kinetic parameter
$\alpha$, as obtained by setting $\mathrm{det}\left({\mathcal M}_0\right)=0$, is
\begin{equation}\label{sec:exper_two:HB}
   \alpha_{jH} \equiv - u_{1js}^2
  +\frac{1}{2\zeta_c}\left[- (1+\kappa_j) + \sqrt{
      (1+\kappa_j)^2 + 8\mu_c\zeta_c u_{1js} }\right] \,, \quad
  j\in \lbrace{1,2\rbrace}\,.
\end{equation}

\begin{figure}[!h]
\centering
\includegraphics[width=88mm]{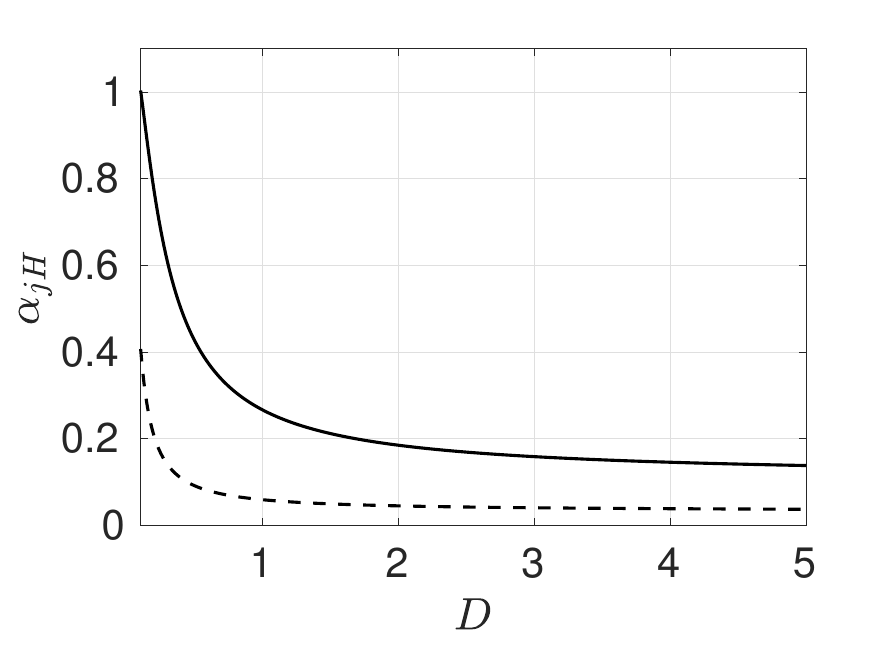}
\caption{HB thresholds $\alpha_{jH}$ versus $D$ from
  \eqref{sec:exper_two:HB} for the leading order GCEP problem
  $\mathrm{det}\left({\mathcal M}_0\right)=0$ for a two-cell
  configuration with $|\v{x}_1-\v{x}_2|=2$. For the cell at $\x_1$
  (top curve), $d_{11}=0.8$ and $d_{21}=0.4$. For the cell at $\x_2$
  (bottom curve), $d_{12}=0.4$ and $d_{22}=0.8$. Intracellular
  oscillations occur when $\alpha$ is decreased below these
  thresholds. Other parameters: $\eps=0.03$, $\sigma=0.5$,
  $\zeta_c=0.15$ and $\mu_c=2.0$. Since cell 1 has larger influx and
  smaller efflux than cell 2, we observe that there is a larger range
  of the reaction-kinetic parameter $\alpha$ where cell 1 is predicted
  to exhibit intracellular oscillations.}
\label{fig:two:selkov}
\end{figure}

To determine the effect of the permeabilities, we assume that
${d_{11}/d_{12}}>1$ and ${d_{21}/d_{22}}<1$, so that
$\kappa_1<\kappa_2$. For this scenario, the influx from the bulk to
cell 1 is larger than for cell 2, while the efflux from cell 1 is
smaller than for cell 2. As a result, we expect
$\alpha_{1H}>\alpha_{2H}$ so that cell 1 will have a larger range of
$\alpha$ where intracellular oscillations will occur. For a particular
parameter set, this is illustrated in Fig.~\ref{fig:two:selkov} where
we plot the two HB thresholds in \eqref{sec:exper_two:HB} versus
$D$. From this figure, we also observe that $\alpha_{jH}$ decreases as
the bulk diffusivity $D$ increases. We can qualitatively
  interpret this trend by recalling that Brownian motion in $\R^3$ is
  non-recurrent. When $D$ increases, the bulk-diffusing species that
  mediates intercellular signaling is more likely to undergo large
  excursions away from the two cells and escape to infinity. In
  effect, increasing $D$ reduces the capture cross-section of the
  two-cell configuration, which lowers the probability of the bulk
  signal hitting a cell.


To include the effect of the ${\mathcal O}(\eps)$ terms in the GCEP
\eqref{pert:all_main}, we write \eqref{pert:all_main} as
\begin{equation}\label{sec:twocell:gcepfull}
  \left( {\mathcal M_0}(\lambda)-4\pi \eps {\mathcal G}_{\lambda}\right)
  \v{c}=\v{0} \,, \quad \mbox{with} \quad {\mathcal M}_0(\lambda) =
\left( {\begin{array}{cc}
          m_1(\lambda,\alpha_c) & 0 \\
          0 & m_2(\lambda,\alpha_c)
        \end{array} } \right)\,,
\end{equation}
where
\begin{equation}\label{sec:pert:m0}
  \quad m_{j}(\lambda,\alpha_c)\equiv 1+\frac{D}{d_{1j}} +
  \frac{4\pi D d_{2j}}{d_{1j}} \v{e}_1^{T}(\lambda I - J_j(\alpha_c))^{-1}\v{e}_1\,,
\end{equation}
where the Jacobians of the intracellular kinetics at the steady-state are
\begin{equation}\label{sec:twocell:jac}
  J_j(\alpha_c) \equiv
  \left( {\begin{array}{cc}
  -1 + \frac{2\mu_c u_{1js}}{\alpha_c+u_{1js}^2} & \alpha_c + u_{1js}^2 \\
 -\frac{2\zeta_c\mu_c u_{1js}}{\alpha_c+u_{1js}^2} & -\zeta_c(\alpha_c+u_{1js}^2)
            \\
  \end{array} } \right)\,, \quad j\in \lbrace{1,2\rbrace}\,.
\end{equation}

As consistent with Fig.~\ref{fig:two:selkov}, we assume that as
$\alpha$ is decreased from large values the HB threshold first occurs
for cell 1. As such, we seek a critical value of $\alpha$ within an
${\mathcal O}(\eps)$ neighborhood of the leading-order threshold
$\alpha_{1H}$ for which the GCEP \eqref{sec:twocell:gcepfull} has a
Hopf bifurcation. To do so, in \eqref{sec:twocell:gcepfull} we expand
$\lambda$, $\v{c}$ and $\alpha_c$ as
\begin{equation}\label{sec:twocell:expan}
  \lambda=\lambda_{0}  + \eps \lambda_1 + \ldots \,, \quad
  \v{c}=\v{c}_0+\eps\v{c}_1+\ldots\,, \quad \alpha_c=\alpha_{1H}+\eps\alpha_{1}
  +\ldots \,,
\end{equation}
where $m_1(\lambda_{0},\alpha_{1H})=0$ with
$\lambda_0\equiv i\lambda_{0I1}$ yields the leading-order HB threshold
at which (see \eqref{selkov:lambda0})
\begin{equation}\label{sec:twocell:det_trace}
  \lambda_{0I1} = \sqrt{(1+\kappa_1) \mathrm{det} J_{1}(\alpha_{1H})} =
    \sqrt{(1+\kappa_1)\zeta_c
  \left(\alpha_{1H}+u_{11s}^2\right)} \,, \quad \mathrm{tr}(J_1)=\kappa_1\,.
\end{equation}
Upon substituting \eqref{sec:twocell:expan} into
\eqref{sec:twocell:gcepfull} and imposing $\v{c}^{T}\v{c}=1$, we
find at leading order that $\v{c}_0=\v{e}_1=(1,0)^T$. At
${\mathcal O}(\eps)$, we get $\v{c}_1^{T}\v{c}_0=0$ and obtain the
matrix problem
\begin{equation}\label{sec:twocell:matrix_1}
\left( {\begin{array}{cc}
          0 & 0 \\
          0 & m_2^{0}
        \end{array} } \right) \v{c}_1= -\lambda_1
\left( {\begin{array}{cc}
      \partial_{\lambda} m_{1}^{0} & 0 \\
          0 & \partial_{\lambda}m_{2}^{0}
        \end{array} } \right)\v{c}_0 -\alpha_1
    \left( {\begin{array}{cc}
       \partial_{\alpha_c}m_{1}^{0} & 0 \\
          0 & \partial_{\alpha_c}m_{2}^{0}
        \end{array} } \right)\v{c}_0
+4\pi  {\mathcal G}_{\lambda_0} \v{c}_{0}\,,
\end{equation}
where we have labeled $\lambda_0=i\lambda_{0I1}$,
$m_2^{0}=m_2(\lambda_0,\alpha_c)$,
$\partial_{\lambda}m_{j}^{0}\equiv\partial_{\lambda}m_{j}(\lambda_0,\alpha_c)$
and
$\partial_{\alpha_c}m_{j}^{0}\equiv\partial_{\alpha_c}m_{j}(\lambda_0,\alpha_c)$,
which are to be evaluated at the leading-order threshold
$\lambda_0=i\lambda_{0I1}$ and $\alpha_{c}=\alpha_{1H}$. The
solvability condition for \eqref{sec:twocell:matrix_1} is that the
right-hand side is orthogonal to $\v{c}_0^{T}$. In this way, we obtain
that
\begin{equation}\label{sec:twocell:pert_sol}
 \lambda_1 \partial_{\lambda} m_{1}^{0}  + \alpha_1 \partial_{\alpha_c}m_{1}^{0}
  =4\pi R_{\lambda_0}\,,
  \qquad \v{c}_1= \frac{4\pi G_{\lambda_0}(\x_1,\x_2)}{m_{2}^{0}} \,
  \left( {\begin{array}{c}
       0 \\
       1
        \end{array} } \right)\,.
\end{equation}
Since $|\x_1-\x_2|=2$, \eqref{4:ge_lambda_ex} yields 
$R_{\lambda_0}\equiv (4\pi)^{-1}\sqrt{(\sigma+i\lambda_{0I1})/D}$ and
$G_{\lambda_0}\equiv -(8\pi)^{-1}e^{-2\sqrt{(\sigma+i\lambda_{0I1})/D}}$.

Next, we calculate $\partial_{\lambda}m_{1}^{0}$ and
$\partial_{\alpha_c} m_{1}^{0}$ by recalling that if a square
matrix $A$ is invertible and depends smoothly on a parameter $p$, then
${d A^{-1}/dp}=-A^{-1} \left({dA/dp}\right) A^{-1}$. By using this
result with $A=\lambda I - J_1(\alpha_c)$ and with $p$ either
$\lambda$ or $\alpha_c$, we conclude that
\begin{subequations}\label{sec:twocell:mder}
\begin{align}
 \partial_{\lambda}m_{1}^{0} &= -\frac{4\pi D d_{21} }{d_{11}} \v{e}_1^T\left(
          \lambda_0 I - J_1(\alpha_{1H})\right)^{-2} \v{e}_1 \,,
  \label{sec:twocell:mder_1}\\
\partial_{\alpha_c} m_{1}^{0} &= \frac{4\pi D d_{21}}{d_{11}} \v{e}_1^T\left(
    \lambda_0I - J_1(\alpha_{1H})\right)^{-1} \frac{d\, J_1(\alpha_{1H})}
    {d\alpha_c}\left(\lambda_0 -
      J_1(\alpha_{1H})\right)^{-1} \v{e}_1 \,,  \label{sec:twocell:mder_2}
\end{align}
\end{subequations}
where $\lambda_0=i\lambda_{0I1}$. By differentiating
the Jacobian in \eqref{sec:twocell:jac} with $j=1$, we identify in
\eqref{sec:twocell:mder_2} that
\begin{equation}\label{sec:twocell:der_jac}
  \frac{d\, J_1(\alpha_{1H})}{d\alpha_c}\equiv
  \left( {\begin{array}{cc}
  -\frac{2\mu_c u_{11s}}{\left(\alpha_{1H}+u_{11s}^2\right)^2} & 1 \\
 \frac{2\zeta_c\mu_c u_{11s}}{\left(\alpha_{1H}+u_{11s}^2\right)^2} & -\zeta_c
            \\
  \end{array} } \right)\,.
\end{equation}

By using \eqref{sec:twocell:jac} for $J_1$, we obtain after some algebra that
\eqref{sec:twocell:mder_1} reduces to
\begin{equation}\label{sec:twocell:h1}
 \partial_{\lambda} m_{1}^{0}=-\frac{4\pi D d_{21}}{d_{11}} \left(\frac{
    \lambda_0^2 + 2\lambda_0 \, \mathrm{det}{J_1} -\mathrm{det}{J_1}\left(
      \mathrm{tr}{J_1} + 1\right)}{
 \left(\lambda_0^2 - \lambda_0 \, \mathrm{tr}{J_1} + \mathrm{det}{J_1}\right)^2}
  \right)\,.
\end{equation}
By evaluating \eqref{sec:twocell:h1} at
$\lambda_0=i\lambda_{01H}$, and using \eqref{sec:twocell:det_trace}, we
find that \eqref{sec:twocell:h1} simplifies to
\begin{equation}\label{sec:twocell:m1lam}
  \partial_{\lambda}m_{1}^{0}=\frac{8\pi D d_{21}}{d_{11}\kappa_1^2}
  \left(-1 + \frac{i \lambda_{0I1}}{1+\kappa_1}\right)^{-1} \,.
\end{equation}
In a similar way, from a lengthy but straightforward calculation using
\eqref{sec:twocell:der_jac} we find that \eqref{sec:twocell:mder_2}
simplifies to
\begin{equation}\label{sec:twocell:m1alp}
 \partial_{\alpha_c} m_{1}^{0} = - \frac{4\pi D d_{21}}{d_{11}}\left(
    \frac{2\mu_c u_{11s}} {\left(\alpha_{1H}+u_{11s}^2\right)^2}\right)
\left( \frac{\lambda_0^2 + 2\lambda_0 \, \mathrm{det}{J_1}}{
    \left(\lambda_0^2 - \lambda_0 \, \mathrm{tr}{J_1} +
      \mathrm{det}{J_1}\right)^2}\right)\,.
\end{equation}

We first show that, if we set $\alpha_1=0$ in \eqref{sec:twocell:pert_sol},
then $\mathrm{Re}(\lambda_1)<0$. When $\alpha_1=0$, we obtain from
\eqref{sec:twocell:pert_sol} and \eqref{sec:twocell:m1lam} that
\begin{equation}
  \lambda_1= \frac{d_{11}\kappa_1^2}{2D d_{21}} R_{\lambda_0}
  \left(-1 + i \frac{\lambda_{0I1}}{1+\kappa_1}\right) \,.
\end{equation}
Since
$R_{\lambda_0}\equiv (4\pi)^{-1}\sqrt{(\sigma+i\lambda_{0I1})/D}$ is
the principal branch and $\lambda_{0I1}>0$, it follows that
$\mathrm{Re}\left(R_{\lambda_0}\right)>0$ and
$\mathrm{Im} \left(R_{\lambda_0}\right)>0$. We readily conclude that
$\mathrm{Re}(\lambda_1)<0$. As a result, we have established that at
$\alpha_c=\alpha_{1H}$, the ${\mathcal O}(\eps)$ term in the GCEP
\eqref{sec:twocell:gcepfull} pushes the purely complex conjugate pair
of eigenvalues associated with the leading-order GCEP into the stable
left-half plane $\mathrm{Re}(\lambda)<0$.

This indicates that we must decrease $\alpha_c$ below $\alpha_{1H}$ 
to obtain the HB threshold for the GCEP
\eqref{sec:twocell:gcepfull}. To calculate this shift, we set
$\mathrm{Re}(\lambda_1)=0$ in \eqref{sec:twocell:pert_sol} to get
that
\begin{equation}\label{sec:twocell:alpha_1}
  \alpha_{1}=4\pi \frac{\mathrm{Re}\left[{R_{\lambda_0}/\partial_{\lambda}m_{1}^{0}}
      \right]}{\mathrm{Re}\left[{\partial_{\alpha_c}m_{1}^{0}/
    \partial_{\lambda} m_{1}^{0}}\right]}\,.
\end{equation}
Upon substituting \eqref{sec:twocell:det_trace},
\eqref{sec:twocell:h1}, \eqref{sec:twocell:m1alp} and
$R_{\lambda_0}\equiv (4\pi)^{-1}\sqrt{(\sigma+i\lambda_{0I1})/D}$ into
\eqref{sec:twocell:alpha_1}, we obtain an explicit expression for
$\alpha_1$ given by
\begin{equation}\label{sec:twocell:alpha_1final}
  \alpha_1 = \frac{\left(\alpha_{1H}+u_{11s}^2\right)^{2} d_{11}\kappa_1^2}{
    8\pi D^{3/2} d_{21}\mu_c u_{11s}} \left(\frac{1+\kappa_1 + \mathrm{det}{J_1}
    }{ 1+\kappa_1 + 2 \mathrm{det}{J_1} } \right) \mathrm{Re}\left[
      \sqrt{\sigma + i \lambda_{0I1}} \left(-1 + i \frac{\lambda_{0I1}}
        {1+\kappa_1}\right) \right] \,.
\end{equation}
It follows that $\alpha_1<0$ as expected.  We conclude that in order
to obtain a HB threshold for the $\epsilon>0$ GCEP
\eqref{sec:twocell:gcepfull}, we must decrease $\alpha_c$ by an
${\mathcal O}(\eps)$ amount from the leading-order threshold
$\alpha_{1H}$. Moreover, at the HB point
$\alpha_c\sim\alpha_{1H}+ \eps\alpha_1$, the result
\eqref{sec:twocell:pert_sol} for the eigenvector perturbation $\v{c}_1$
shows that it has a nonzero component in cell 2. From
\eqref{nstabform:c2}, this implies that cell 2 will develop small
amplitude oscillations near the steady-state that are triggered by the
Hopf bifurcation associated with cell 1.

\begin{figure}[!h]
\centering
\includegraphics[width=88mm]{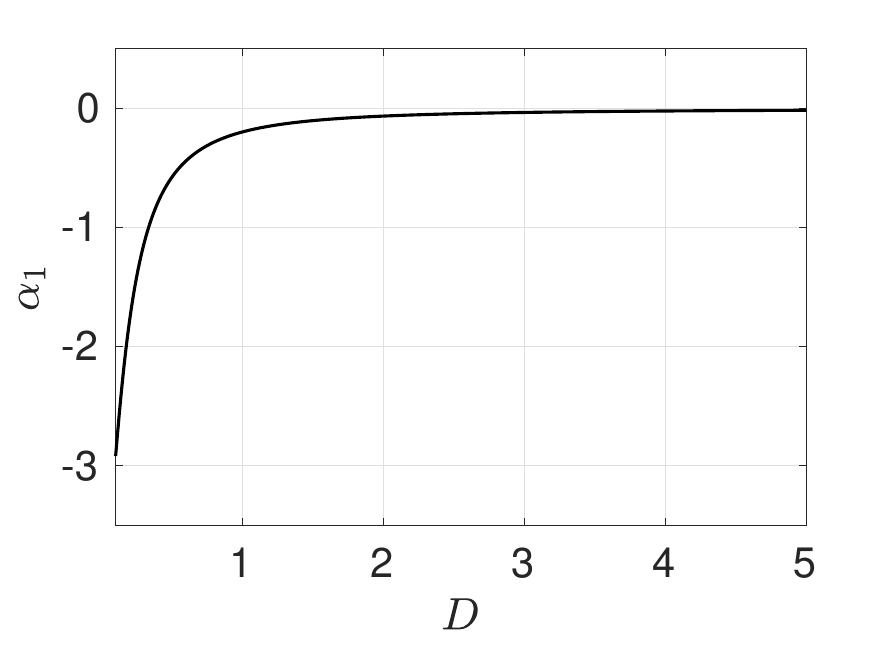}
\caption{The correction $\alpha_1$, given in \eqref{sec:twocell:alpha_1final},
  versus $D$ for the HB threshold for cell 1 for the two-cell
  configuration in Fig.~\ref{fig:two:selkov} with parameters given in
  the figure caption.}
\label{fig:twocell_alpha1}
\end{figure}

In Fig.~\ref{fig:twocell_alpha1} we plot $\alpha_1$ versus $D$ for
cell 1 in the two-cell configuration with parameter values in the
caption of Fig.~\ref{fig:two:selkov}. From this figure, we observe
that $\alpha_1<0$ but that it monotonically increases to zero as $D$
increases. This implies that, for a given cell radius $\eps$, the
leading-order HB threshold $\alpha_{1H}$ will provide an increasingly
accurate prediction of the perturbed HB threshold as $D$ increases.

\subsection{Identical Cells and where \texorpdfstring{${\mathcal G}_{\lambda}$}{G} has a constant row sum}\label{sec:exper_circ}

\begin{figure}[htbp]
  \centering
      \begin{subfigure}[b]{0.48\textwidth}
        \includegraphics[width=\textwidth]{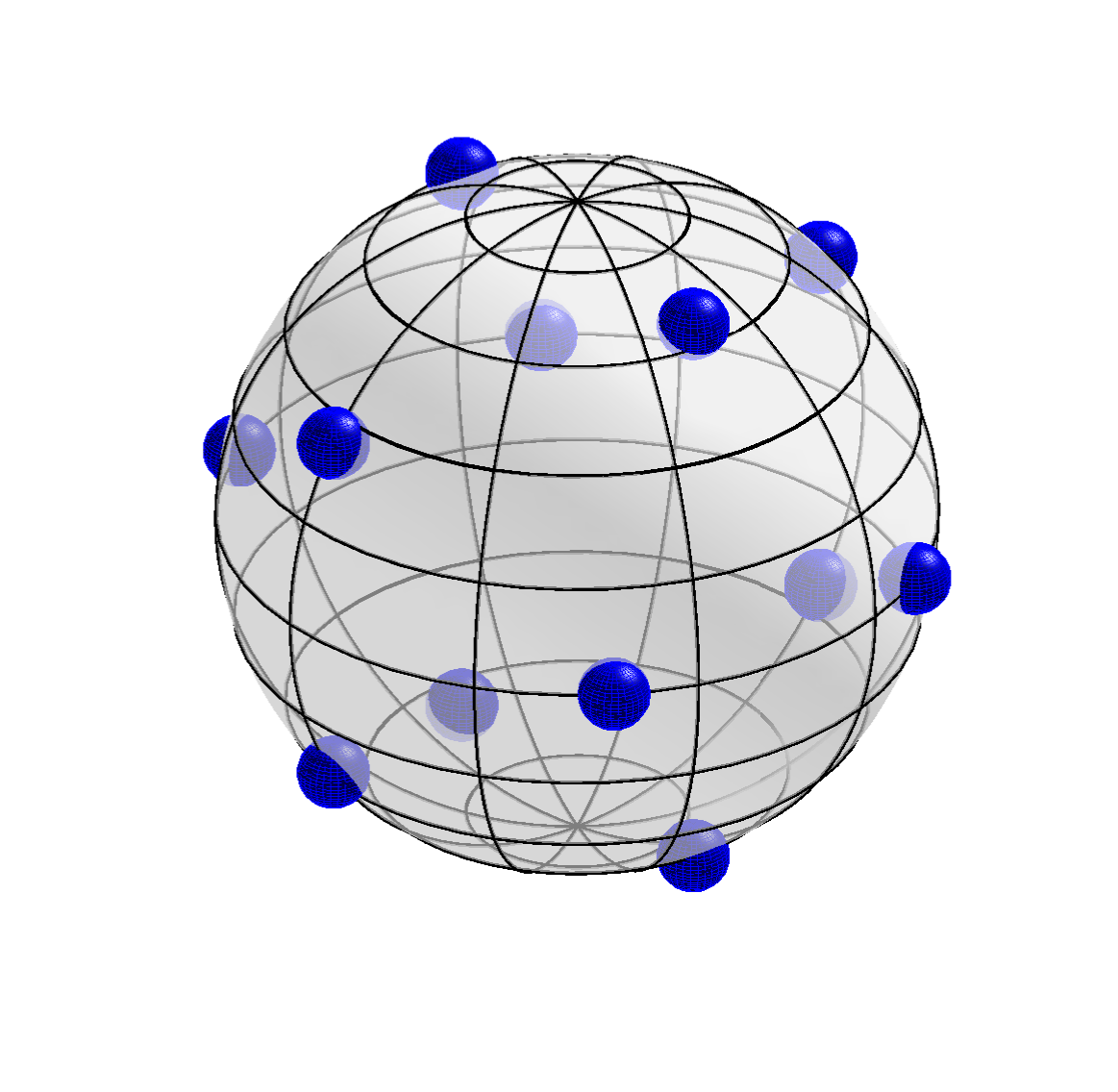}
        \caption{icosahedron cell configuration}
        \label{fig:icosa}
    \end{subfigure}
      \begin{subfigure}[b]{0.48\textwidth}
        \includegraphics[width=\textwidth]{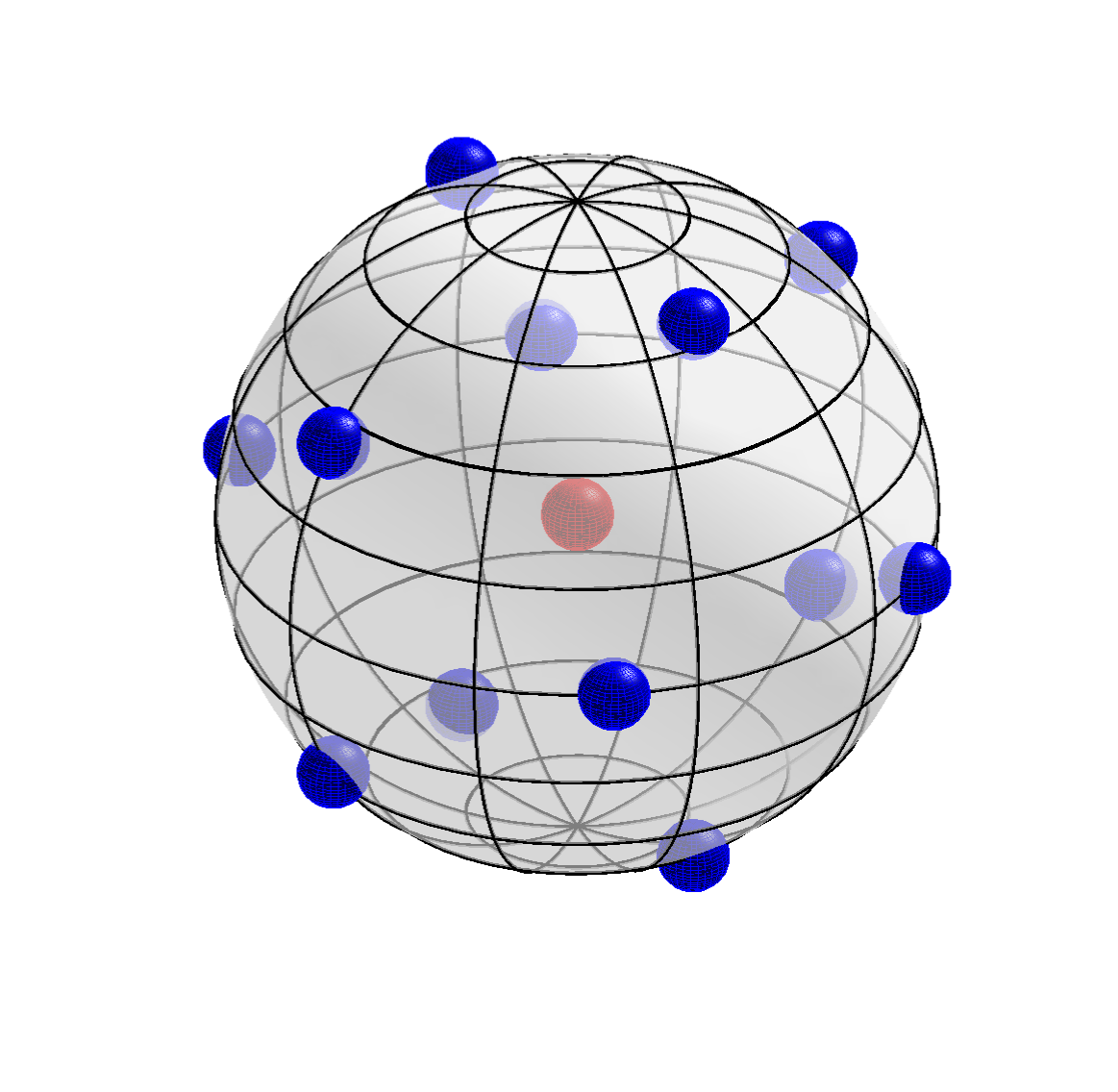}
        \caption{adding a cell in the center as a pacemaker}
        \label{fig:icosa_center}
    \end{subfigure}
    \caption{Left: Identical signaling cells in blue centered at the
      vertices of an icosahedron that lie on the unit sphere. Right:
      Same plot except but with an additional pacemaker (cell in red) centered
      at the origin that can have different permeabilities and
      reaction-kinetic parameters.}
    \label{fig:icosahedron}
\end{figure}

We now consider $N$ cells with common permeability and reaction
kinetic parameters for $j\in\lbrace{1,\ldots,N\rbrace}$, labeled by
\begin{equation}\label{num_ident:param}
  d_{1c}=d_j\,, \quad d_{2c}=d_{2j} \,, \quad \mu_c=\mu_j\,, \quad
  \zeta_c=\zeta_j\,, \quad \alpha_c=\alpha_j \,, \quad
  \kappa_c=\kappa_j\,,
  \end{equation}
  where $\kappa_c\equiv {4\pi D d_{2c}/(d_{1c}+D)}$. Moreover, we
  assume that the cell centers $\x_1,\ldots,\x_N$ are such that the
  complex-symmetric Green's matrix ${\mathcal G}_{\lambda}$ has a constant
  row sum. This matrix structure arises for any two-cell
  configuration, for cells centered at the vertices of the five
  platonic solids, for cells equally-spaced along a line of constant
  latitude on a spherical shell, among other choices. A subset of such
  matrices are those that have a circulant structure. A schematic
  illustration of an icosahedron spatial arrangement of cells centered
  on the unit sphere is shown in Fig.~\ref{fig:icosa}. The effect of
  adding a pacemaker cell to this arrangement, centered at the origin
  (see Fig.~\ref{fig:icosa_center}), will be analyzed below in \S
  \ref{sec:exper_circ_p}.

  When ${\mathcal G}_{\lambda}$ has a constant row sum, it has the
  eigenspace
  \begin{subequations}\label{num_ident:qj_def}
  \begin{equation}  
    {\mathcal G}_{\lambda}\v{e}=\gamma_N(\lambda) \v{e} \,, \qquad
    {\mathcal G}_{\lambda}\q_j=\gamma_j(\lambda) \q_j \,, \qquad
    j \in \lbrace{1,\ldots,N-1\rbrace} \,,
  \end{equation}
  where $\v{e}=(1,\ldots,1)^T$ with $^{\mathrm{T}}$ denoting
  transposition. We have the orthogonality relations
  \begin{equation}\label{num_ident:qj_ortho}
    \v{e}^{T}\q_j=0 \,, \quad j\in\{1,\ldots,N-1\} \,; \qquad
    \q_k^T\q_j=0 \,, \quad j\neq k \,, \quad j,k\in \lbrace{1,\ldots,N-1
      \rbrace}\,.
  \end{equation}
  \end{subequations}
   We refer to the eigenvector $\v{e}$ as the {\em synchronous mode}, and
   its corresponding eigenvalue $\gamma_N=\gamma_{N}(\lambda)$ is
\begin{equation}\label{num_ident:sync_eig}
 \gamma_N=R_{\lambda} + \sum\limits_{m=1}^{N-1} 
  \left( {\mathcal G}_{\lambda}\right)_{1,1+m}\,.
\end{equation}
We refer to the eigenvectors $\q_j$ for
$j\in\lbrace{1,\ldots,N-1\rbrace}$ as the {\em asynchronous
  modes}. Since ${\mathcal G}_{\lambda}$ is complex symmetric, these
modes will typically have higher multiplicity.

Since ${\mathcal G}_0\v{e}=\gamma_{N0}\v{e}$, where $\gamma_{N0}=
\gamma_{N}(0)$, the intracellular steady-states are the same for each cell.
More specifically, the solution to the steady-state problem \eqref{selkov:ss} 
is
\begin{equation}\label{num_ident:ss}
 \v{u}_{1s} = u_{c1}\v{e} \,, \quad \mbox{where} \quad
 u_{c1}\equiv \mu_{c} \left(1+\kappa_c + \frac{\eps d_{1c}}{D d_{2c}} \kappa_c^2
   \gamma_{N0}\right)^{-1} \,.
\end{equation}
From the linearization of this steady-state, we will analyze stability
thresholds associated with the GCEP \eqref{pert:all_main} where the
reaction kinetic parameter $\alpha_c$ is chosen as the primary
bifurcation parameter.  From \eqref{num_ident:ss} we observe that the
steady-state $u_{c1}$ is independent of $\alpha_c$.  For identical
cells, we write the GCEP \eqref{pert:all_main} as
\begin{equation}\label{num_ident:gcep}
  \left( m_{c} I - 4\pi \eps {\mathcal G}_{\lambda}\right)
  \v{c}=\v{0} \,, \quad m_c=m_{c}(\lambda,\alpha_c)\equiv 1+\frac{D}{d_{1c}} +
  \frac{4\pi D d_{2c}}{d_{1c}} \v{e}_1^{T}\left(\lambda I - J_c(\alpha_c)
  \right)^{-1}\v{e}_1\,,
\end{equation}
where the common Jacobian $J_c(\alpha_c)$ of the Sel'kov intracellular
kinetics is given by 
\begin{equation}\label{num_ident:jac}
  J_c(\alpha_c) \equiv
  \left( {\begin{array}{cc}
  -1 + \frac{2\mu_c u_{c1}}{\alpha_c+u_{c1}^2} & \alpha_c + u_{c1}^2 \\
 -\frac{2\zeta_c\mu_c u_{c1}}{\alpha_c+u_{c1}^2} & -\zeta_c(\alpha_c+u_{c1}^2)
            \\
  \end{array} } \right)\,.
\end{equation}
We will determine the critical values of $\alpha_c$ for which the
linearization of the steady-state undergoes a Hopf bifurcation.

To leading order, where we set $\eps=0$ in \eqref{num_ident:gcep}, we
conclude from the vanishing of the trace in \eqref{selkov:det_trace}
for the case of identical cells that
$m_{c}(\pm i \lambda_{0I},\alpha_{cH})=0$ when
\begin{equation}\label{num_ident:HB_zero}
  \alpha_{cH} \equiv - u_{c1}^2
  +\frac{1}{2\zeta_c}\left[- (1+\kappa_c) + \sqrt{
      (1+\kappa_c)^2 + 8\mu_c\zeta_c u_{c1}}\right]\,.
\end{equation}
For this value of $\alpha_{cH}$, the principal eigenvalues of the
limiting GCEP consist of an $N$-fold degenerate complex conjugate pair
of eigenvalues on the pure imaginary axis at
$\lambda=\pm i\lambda_{0I}$, given at the leading-order
threshold $\alpha_c=\alpha_{cH}$ by
\begin{equation}\label{num_ident:lam0}
  \lambda_{0I} = \sqrt{\zeta_c(1+\kappa_c) \mathrm{det}(J_c)} \,,
  \quad \mbox{where} \quad \mathrm{det}(J_c)=\alpha_{cH}+u_{c1}^2 \,, \quad
  \mathrm{tr}(J_c)=\kappa_c \,.
\end{equation}

\begin{figure}[htbp]
  \centering
      \begin{subfigure}[b]{0.48\textwidth}
      \includegraphics[width=\textwidth]{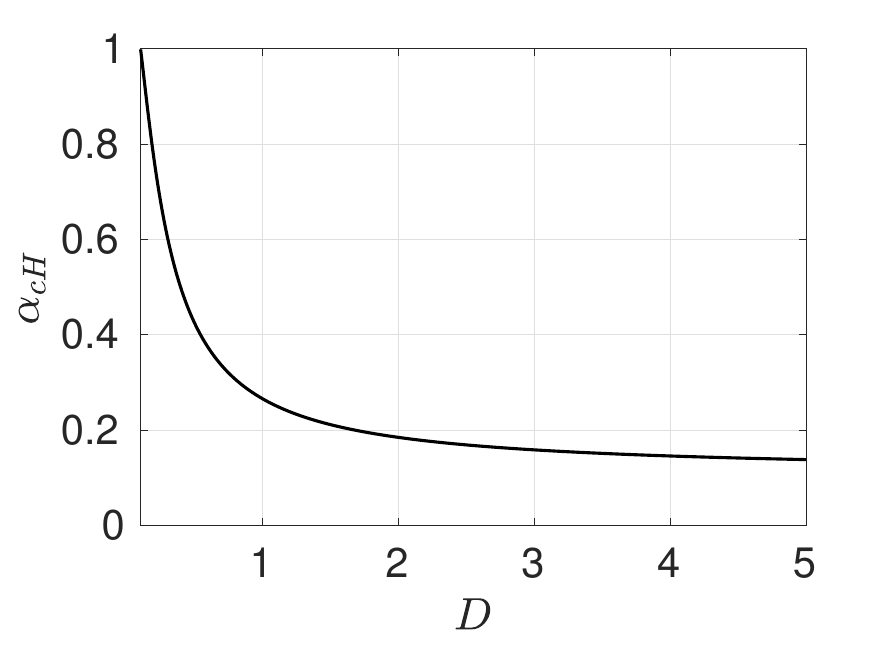}
        \caption{$\alpha_{cH}$} \label{fig:ident_alpha12}
    \end{subfigure}
      \begin{subfigure}[b]{0.48\textwidth}
      \includegraphics[width=\textwidth]{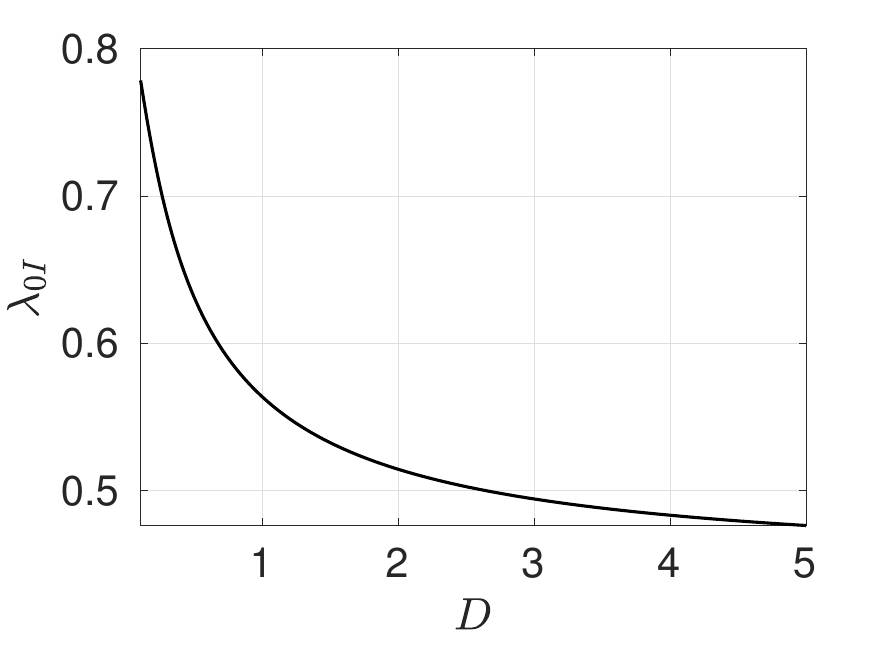}
        \caption{$\lambda_{0I}$} \label{fig:ident_lambda12}
    \end{subfigure}
    \caption{Left: The HB boundary $\alpha_{cH}$ versus $D$ in
      \eqref{num_ident:HB_zero} for an icosahedron configuration
      ($N=12$) of identical cells on the surface of the unit
      sphere. Right: the corresponding HB eigenvalue $\lambda_{0I}$
      versus $D$ given in \eqref{num_ident:lam0}. Parameters:
      $d_{1c}=0.8$, $d_{2c}=0.4$, $\sigma=0.5$, $\eps=0.03$,
      $\mu_c=2.0$ and $\zeta_c=0.15$.}
    \label{fig:ident_icosa}
\end{figure}

For a particular parameter set, in Fig.~\ref{fig:ident_alpha12} and
Fig.~\ref{fig:ident_lambda12} we plot $\alpha_{cH}$ and
$\lambda_{0I}$ versus $D$ from \eqref{num_ident:HB_zero} and
\eqref{num_ident:lam0}, respectively, when the cells are centered
at the vertices of an icosahedron. The leading-order theory
predicts that the steady-state is linearly unstable when
$\alpha_c<\alpha_{cH}$.

In terms of the spatial configuration of cells, we use the matrix
perturbation approach of \S \ref{sec:eig_pert} to analyze how this
degenerate HB threshold is unfolded by the perturbation induced by the
${\mathcal O}(\eps)$ term in the GCEP \eqref{num_ident:gcep}. Our
main result is as follows:

\begin{prop}\label{prop:stab_circ} When $N$ cells are arranged so that
  the complex symmetric matrix ${\mathcal G}_{\lambda}$ is
  diagonalizable with a constant row sum, the critical values of the
  reaction-kinetic parameter $\alpha_{c}$ for which a HB occurs for
  the linearization of the steady-state of the cell-bulk model for
  Sel'kov intracellular kinetics \eqref{isolated:selkov:ode} is 
  \begin{equation}\label{numer_ident:main}
    \alpha_c=\alpha_{cH}+\eps\alpha_{1j} + {\mathcal O}(\eps^2)\,, \quad
    j\in \lbrace{1,\ldots,N\rbrace}\,,
  \end{equation}
  where $\alpha_{cH}$ is given in \eqref{num_ident:HB_zero} and where
  \begin{equation}\label{numer_ident:alpha_1j}
    \alpha_{1j}=4\pi \frac{\mathrm{Re}\left[{\gamma_j(i\lambda_{0I})/
          \partial_{\lambda} m_{c}^{0}}\right]}{
      \mathrm{Re}\left[{\partial_{\alpha_c}m_{c}^{0}/\partial_{\lambda}m_{c}^{0}}
  \right]}\,, \quad
     j\in \lbrace{1,\ldots,N\rbrace} \,.
  \end{equation}
  In \eqref{numer_ident:alpha_1j}, $\gamma_j(\lambda)$ for
  $j\in \lbrace{1,\ldots,N\rbrace}$ are the eigenvalues of
  ${\mathcal G}_{\lambda}$ (see \eqref{num_ident:qj_def}), while
  $\partial_{\lambda}m_{c}^{0}\equiv\partial_{\lambda}{m_c(\lambda,
    \alpha_{c})}$ and
  $\partial_{\alpha_c}m_{c}^{0}\equiv\partial_{\alpha_c}{m_c(\lambda,
    \alpha_{c})}$ are the partial derivatives evaluated at the
  leading-order threshold $\lambda=\lambda_0\equiv i\lambda_{0I}$ and
  $\alpha_c=\alpha_{cH}$. The perturbed thresholds are given explicitly
  by
\begin{equation}\label{num_ident:alpha_jfinal}
  \alpha_{1j} = \frac{\left(\alpha_{cH}+u_{c1}^2\right)^{2} d_{1c}\kappa_c^2}{
    2\mu_c u_{c1} d_{2c} D} \left(\frac{1+\kappa_c + \mathrm{det}{J_c}
    }{ 1+\kappa_c + 2 \mathrm{det}{J_c} } \right) \mathrm{Re}\left[
      \gamma_{j}(i\lambda_{0I}) \left(-1 + i \frac{\lambda_{0I}}
      {1+\kappa_c}\right) \right] \,, \quad j\in\lbrace{1,\ldots,N\rbrace}\,,
\end{equation}
where $\lambda_{0I}$, $\mathrm{det}(J_c)$ and $u_{c1}$ are given in
\eqref{num_ident:lam0} and \eqref{num_ident:ss}.
\end{prop}

\begin{proof}

  We set $\v{c}=\v{q}_j$ in \eqref{num_ident:gcep} where $\v{q}_j$ is
  is an eigenvector of ${\mathcal G}_{\lambda}$ with eigenvalue
  $\gamma_j(\lambda)$ (see \eqref{num_ident:qj_def}). This yields the
  $N$ scalar problems
\begin{equation}\label{num_ident:scalar_eig}
  m_{c}(\lambda, \alpha_c)-4\pi \eps \gamma_j(\lambda)=0 \,,
  \qquad j\in \lbrace{1,\ldots,N\rbrace} \,.
\end{equation}
By expanding $\alpha_c=\alpha_{cH}+\eps \alpha_{1j}$ and
$\lambda=i\lambda_{0I}+\eps\lambda_1$, where
$m_{c}(i\lambda_{0I},\alpha_{cH})=0$, we find at
${\mathcal O}(\eps)$ that
\begin{equation}\label{num_ident:scalar_eig_per}
  \lambda_1 \partial_{\lambda}m_{c}^{0} + \alpha_{1j} 
 \partial_{\alpha_c} m_{c}^{0} = 4\pi \gamma_{j}(i\lambda_{0I}) \,, \quad
   j\in \lbrace{1,\ldots,N\rbrace} \,.
\end{equation}
Upon solving for $\lambda_1$, and setting $\mathrm{Re}(\lambda_1)=0$
for each mode, we readily obtain
\eqref{numer_ident:alpha_1j}. Finally, to obtain the explicit result
in \eqref{num_ident:alpha_jfinal}, we proceed as in
\eqref{sec:twocell:mder}--\eqref{sec:twocell:m1alp} for the two-cell
case to derive after some algebra that
\begin{equation}\label{num_ident:ratio}
  \begin{split}
    \mathrm{Re}{\left[\frac{\partial_{\alpha_c}m_{c}^{0}}{
          \partial_{\lambda}m_{c}^{0}}\right]} &=
  \frac{\mu_c u_{c1}}{\left(\alpha_{cH}+u_{c1}^2\right)^2}\left( 1 +
    \frac{\mathrm{det}(J_c)}{1+\kappa_c + \mathrm{det}(J_c)}\right)>0 \,,\\
\partial_{\lambda} m_{c}^{0} &= \frac{8\pi D d_{2c}}{d_{1c}\kappa_c^2} 
    \left(-1 + \frac{i\lambda_{0I}}{1+\kappa_c}\right)^{-1} \,.
\end{split}
\end{equation}
Upon substituting \eqref{num_ident:ratio} into \eqref{numer_ident:alpha_1j}
we obtain \eqref{num_ident:alpha_jfinal}.
\end{proof}

\begin{figure}[htbp]
  \centering
    \includegraphics[width=88mm]{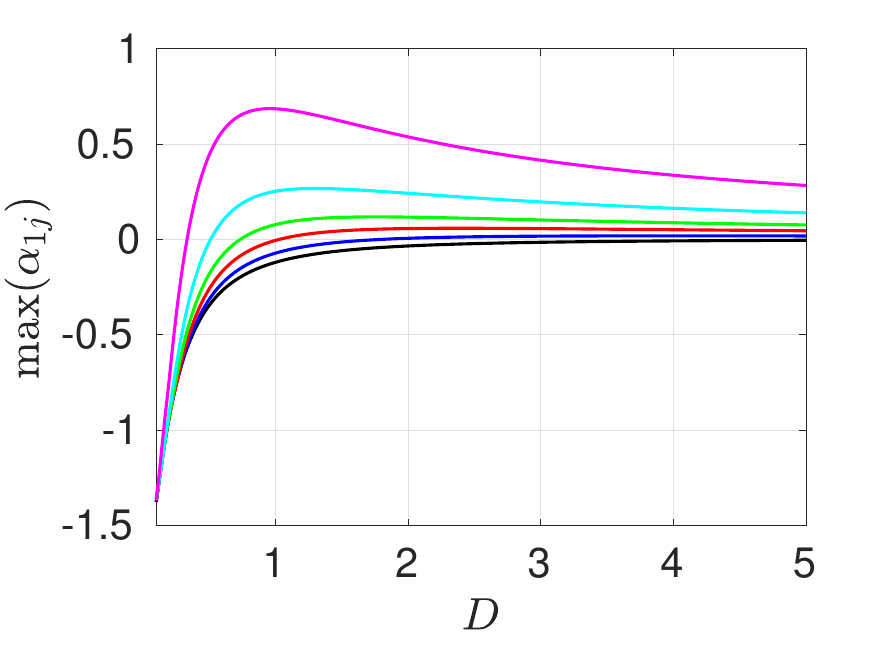}
    \caption{The correction $\mathrm{max}\left({\alpha_{1j}}\right)$,
      obtained from \eqref{ring:alpmax} and
      \eqref{num_ident:alpha_jfinal}, versus $D$ for the first HB
      threshold as $\alpha_c$ decreases below
      $\alpha_{cH}+\eps\mathrm{max}\left({\alpha_{1j}}\right)$. The cell
      configurations with vertices on the unit sphere are: $N=2$
      antipodal (black curve), $N=4$ tetrahedron (dark blue), $N=6$
      octahedron (red), $N=8$ cube-shaped (green), $N=12$ icosahedron
      (cyan), $N=20$ dodecahedron (magenta). Parameters: $d_{1c}=0.8$,
      $d_{2c}=0.4$, $\sigma=0.5$, $\eps=0.03$, $\mu_c=2.0$ and
      $\zeta_c=0.15$.}
    \label{fig:ident_alpha1}
\end{figure}

To determine the ${\mathcal O}(\eps)$ correction to the stability
threshold, we abbreviate
\begin{equation}\label{ring:alpmax}
  \mathrm{max}\left({\alpha_{1j}}\right) \equiv
  \max_{j\in\lbrace{1,\ldots,N\rbrace}} \alpha_{1j}\,.
\end{equation}
We conclude that the Hopf bifurcation {\em first occurs} when
$\alpha_{c}$ decreases below the threshold
$\alpha_{cH}+\eps\mathrm{max}(\alpha_{1j})$.  By numerically evaluating
\eqref{num_ident:alpha_jfinal} to identify the maximum value, in
Fig.~\ref{fig:ident_alpha1} we plot
$\mathrm{max}\left({\alpha_{1j}}\right)$ versus $D$ for two identical
cells and for all configurations where the identical cells are
centered at the vertices on of the five platonic solids.

For the parameter set in Fig.~\ref{fig:ident_alpha1}, and for any of
these cell configurations, our computations show that this first HB
threshold occurs for the in-phase mode $j=N$, for which the
corresponding eigenvector of the GCEP is $\v{c}=\v{e}$. As such, as
$\alpha$ is decreased from a large value, the steady-state first
becomes unstable to an in-phase HB, where all cells synchronize their
oscillations. Since intercellular communication is enhanced
  with an increasing number of cells, we observe as expected from
  Fig.~\ref{fig:ident_alpha1} that the range of $\alpha$ where
  oscillations occur increases with the number of cells for each fixed
  $D$. Moreover, for each of the five platonic solid cell
  configurations, $\mathrm{max}\left({\alpha_{1j}}\right)$ achieves a
  global maximum at some critical value of $D$ that decreases as the
  number of cells increases (this feature is not clearly visible for
  $N=2$ and $N=4$ owing to the vertical scale in
  Fig.~\ref{fig:ident_alpha1}). As a result, we conclude that there
  is some optimal value of the bulk diffusivity that provides the
  largest range in $\alpha_c$ where intracellular oscillations will
  occur. 

\subsection{The Effect of a Pacemaker Cell}\label{sec:exper_circ_p}

In what follows we consider the effect of adding a pacemaker or signaling cell
centered at the origin at $\x_{N+1}=\zero$ to the cell configurations
that were considered in \S \ref{sec:exper_circ} (see
Fig.~\ref{fig:icosa_center}). The cells at the vertices of the
platonic solids on the unit disk are taken to have identical
parameters labeled by
\begin{equation}\label{rc:ident}
  d_{1c}=d_{1j}\,, \quad d_{2c}=d_{2j}\,, \quad \alpha_c=\alpha_j\,, \quad
  \zeta_c=\zeta_j\,, \quad \mu_c=\mu_j \,,
\end{equation}
for $j\in\lbrace{1,\ldots,N\rbrace}$. The pacemaker cell is distinct in that
it can have different influx and efflux permeabilities, and possibly a
different reaction kinetic parameter $\alpha$, so that we label
\begin{equation}\label{rc:pace}
  d_{1p}=d_{1N+1}\,, \quad d_{2p}=d_{2N+1}\,, \quad \alpha_p=\alpha_{N+1}\,,
  \quad \zeta_c=\zeta_{N+1}\,, \quad \mu_c=\mu_{N+1}\,.
\end{equation}

For this configuration, in the steady-state and GCEP problems given by
\eqref{selkov:ss} and \eqref{pert:all_main}, the $(N+1)\times (N+1)$
dimensional Green's matrix, labeled by $\tilde{{\mathcal G}}_{\lambda}$, can
be partitioned into block matrix form as
\begin{equation}\label{ringc:mat_block}
  \arraycolsep=1.4pt\def\arraystretch{2.2}
    \tilde{{\mathcal G}}_{\lambda}=
    \left(
    \begin{array}{c|c}
     {\mathcal G}_{\lambda} & a_{\lambda} \v{e} \\
      \hline
    a_{\lambda} \v{e}^T & R_{\lambda}
    \end{array}
  \right) \,, \quad \mbox{where} \quad a_{\lambda}\equiv
  G_{\lambda}(\x_j,\x_{N+1})=
    -\frac{1}{4\pi} e^{-\sqrt{(\sigma+\lambda)/D}} \,, \quad
    R_{\lambda}= \frac{1}{4\pi}
        \sqrt{ \frac{\sigma+\lambda}{D}} \,.
\end{equation}
Here $\v{e}=(1,\ldots,1)^{\mathrm{T}}\in \R^N$,
${\mathcal G}_{\lambda}$ is the $N\times N$ constant row-sum Green's
matrix associated with the cell-cell interactions on the spherical
ring, and we have used $|\x_j-\x_{N+1}|=1$ for
$j\in\lbrace{1,\ldots,N\rbrace}$.
  
Next, the steady-state intracellular species $\v{u}_{1s}$ satisfying
the $(N+1)\times (N+1)$ dimensional linear system \eqref{selkov:ss}
has the form
\begin{subequations}\label{rc:ss}
  \begin{equation}\label{rc:ss1}
       \v{u}_{1s} = \left(\begin{array}{c}
      u_{c1} \v{e} \\
      u_{p1}
                      \end{array}\right)
             \,, \qquad \mbox{where} \qquad
  \left( I  +   {\mathcal K} + \frac{\eps}{D} {\mathcal P} {\mathcal K}
    \tilde{{\mathcal G}}_0 {\mathcal K}\right) \v{u}_{1s} =\mu_c
\left(\begin{array}{c}
      \v{e} \\
        1
    \end{array}\right)\,.
\end{equation}
In \eqref{rc:ss1}, the diagonal matrices ${\mathcal K}$ and
${\mathcal P}$ are defined by
\begin{equation}\label{rc:ss2}
  {\mathcal P} \equiv \mathrm{diag}\left(\frac{d_{1c}}{d_{2c}}, \ldots,
      \frac{d_{1c}}{d_{2c}},  \frac{d_{1p}}{d_{2p}} \right)\,, \qquad
        {\mathcal K}\equiv \mathrm{diag}\left( \kappa_c,\ldots,\kappa_c,
          \kappa_p \right) \,, 
\end{equation}
\end{subequations}
with $\kappa_c={4\pi D d_{2c}/\left(d_{1c}+D\right)}$ and
$\kappa_p={4\pi D d_{2p}/\left(d_{1p}+D\right)}$. By using the block
form \eqref{ringc:mat_block}, we conclude that
$u_{c1}$ and $u_{p1}$ in \eqref{rc:ss1} satisfy the $2\times 2$ linear
system 
\begin{equation}\label{ringc:equil_b}
      \left( {\begin{array}{cc}
            1+\kappa_c+\frac{\eps d_{1c}}{D d_{2c}}\kappa_c^2 \gamma_{N0} &
             \frac{\eps d_{1c}}{D d_{2c}} a_0 \kappa_c \kappa_p \\
            \frac{\eps d_{1p}}{D d_{2p}} N a_0 \kappa_p \kappa_c &
            1+\kappa_p+\frac{\eps d_{1p}}{D d_{2p}}\kappa_p^2 R_0 \\
  \end{array} } \right) \left(\begin{array}{c}
      u_{c1} \\
      u_{p1}
     \end{array}\right) = \mu_c \left(\begin{array}{c}
      1 \\
      1
     \end{array}\right) \,.
\end{equation}
Here $\gamma_{N0}$ is the eigenvalue for the synchronous mode, i.e.,
${\mathcal G}_{0}\v{e}=\gamma_{N0}\v{e}$, while
$a_0=-(4\pi)^{-1}e^{-\sqrt{\sigma/D}}$ and
$R_0=(4\pi)^{-1}\sqrt{{\sigma/D}}$ from \eqref{ringc:mat_block} and
\eqref{4:ge_lambda_ex}.

For this cell configuration, the GCEP in \eqref{5:gcep_2} is
\begin{equation}\label{ringc:gcep}
  \left( {\mathcal M_0}(\lambda)-4\pi \eps \tilde{{\mathcal G}}_{\lambda}\right)
  \v{c}=\v{0} \,, \quad \mbox{with} \quad
  {\mathcal M}_0(\lambda)=\mathrm{diag}\left(m_{c},\ldots,
    m_{c}, m_{p}\right)\,,
\end{equation}
where $\tilde{{\mathcal G}}_{\lambda}$ is given in \eqref{ringc:mat_block} and
with
\begin{subequations}\label{ringc:diag_gcep}
\begin{equation}\label{ringc:mc_mp}
  \begin{split}
 & m_c=m_{c}(\lambda,\alpha_c) \equiv 1+\frac{D}{d_{1c}} +
 \frac{4\pi D d_{2c}}{d_{1c}} \v{e}_1^{T}(\lambda I - J_c(\alpha_c))^{-1}\v{e}_1\,,
\\
  & m_p=m_{p}(\lambda,\alpha_p)\equiv 1+\frac{D}{d_{1p}} +
 \frac{4\pi D d_{2p}}{d_{1p}} \v{e}_1^{T}(\lambda I - J_p(\alpha_p))^{-1}\v{e}_1\,.
\end{split}
\end{equation}
The common Jacobian $J_c$ for the cells on the ring is given in
\eqref{num_ident:jac} with $u_{c1}$ satisfying  \eqref{ringc:equil_b},
while the Jacobian $J_p$ for the pacemaker cell evaluated at
$u_{p1}$ is
\begin{equation}
  J_p(\alpha_p) \equiv
  \left( {\begin{array}{cc}
  -1 + \frac{2\mu_c u_{p1}}{\alpha_p+u_{p1}^2} & \alpha_p + u_{p1}^2 \\
 -\frac{2\zeta_c\mu_c u_{p1}}{\alpha_p+u_{p1}^2} & -\zeta_c(\alpha_p+u_{p1}^2)
            \\
          \end{array} } \right)\,.
\end{equation}
\end{subequations}

\begin{figure}[htbp]
  \centering
    \begin{subfigure}[b]{0.48\textwidth}
      \includegraphics[width=\textwidth]{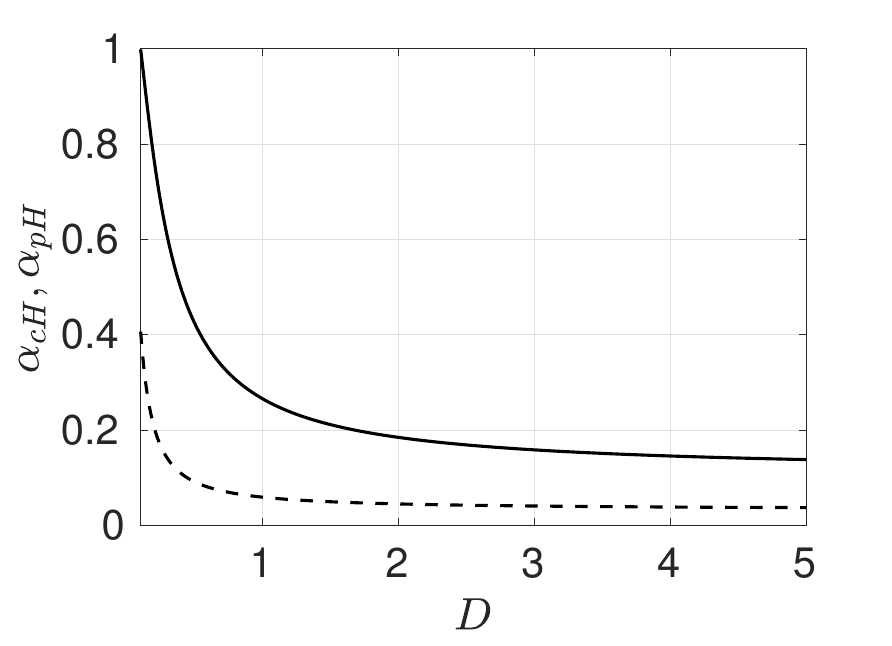}
      \caption{Case I: $\alpha_{cH}>\alpha_{pH}$}
      \label{fig:caseI:HB}
    \end{subfigure}
      \begin{subfigure}[b]{0.48\textwidth}
      \includegraphics[width=\textwidth]{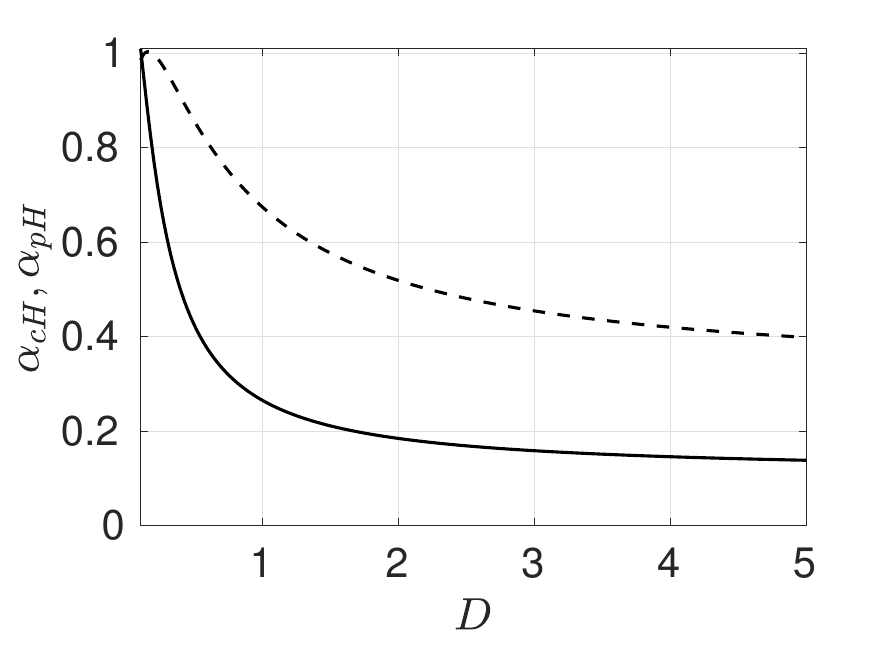}
        \caption{Case II: $\alpha_{pH}>\alpha_{cH}$} \label{fig:caseII:HB}
    \end{subfigure}
    \caption{An icosahedron arrangement of ring cells with an
      additional pacemaker cell at the center. Left: The HB thresholds
      $\alpha_{cH}$ (solid) and $\alpha_{pH}$ (dashed) versus $D$ for
      Case I with parameters $d_{1c}=0.8$, $d_{2c}=0.4$ for ring cells
      and $d_{1p}=0.4$ and $d_{2p}=0.8$ for the pacemaker cell. The ring
      cells set the HB threshold. Right: Same plot but for Case II
      with parameters $d_{1c}=0.8$, $d_{2c}=0.4$ for ring cells and
      $d_{1p}=1.2$ and $d_{2p}=0.2$ for the pacemaker cell.  The
      pacemaker cell now sets the HB threshold. Other parameters:
      $\sigma=0.5$, $\eps=0.03$, $\mu_c=2.0$ and $\zeta_c=0.15$.}
    \label{fig:caseI_II:HB}
  \end{figure}

To leading order in $\eps$, there are two possible thresholds for a
Hopf bifurcation. For {\bf Case I}, the dominant HB threshhold is for
the ring cells and occurs when $\alpha_c=\alpha_{cH}$, where
$m_c(\pm i \lambda_{0I},\alpha_{cH})=0$ and $\alpha_{cH}$ is given
in \eqref{num_ident:HB_zero}. This HB sets the stability threshold
when all roots $\lambda$ to $m_p(\lambda,\alpha_p)=0$ for
the pacemaker cell are in $\mathrm{Re}(\lambda)<0$. Alternatively, in
{\bf Case II}, the leading-order dominant HB threshold occurs for the
pacemaker cell at $\alpha_p=\alpha_{pH}$, where
$m_p(\pm i \lambda_{0I},\alpha_{pH})=0$, with $\alpha_{pH}$ given by
\begin{equation}\label{ringc:alpha_pH}
  \alpha_{pH} \equiv - u_{p1}^2
  +\frac{1}{2\zeta_c}\left[- (1+\kappa_p) + \sqrt{
      (1+\kappa_p)^2 + 8\mu_c\zeta_c u_{p1}}\right]\,.
\end{equation}
This HB arising from the pacemaker cell sets the stability threshold
provided that all the roots $\lambda$ to $m_c(\lambda,\alpha_c)=0$ for
the common ring cells are in $\mathrm{Re}(\lambda)<0$. For two
particular sets of parameters, corresponding to either Case I or II,
in Fig.~\ref{fig:caseI_II:HB} we plot the two leading-order thresholds
$\alpha_{cH}$ and $\alpha_{pH}$ versus $D$ for an icosahedreon
arrangement of cells with an additonal center cell. Case II occurs
when $d_{1p}>d_{1c}$ and $d_{2p}<d_{2c}$ so that the efflux
permeability for the pacemaker is smaller while its influx
permeability is larger than for the cells on the ring.

  For our cell configurations with a pacemaker, since there are only
  two types of eigenvectors $\v{c}_j$ for the block matrix
  \eqref{ringc:gcep} for any $\eps>0$, we can readily determine the
  perturbed HB thresholds for either Case I or II and qualitatively
  discuss the possible modes of oscillation near the onset of the
  dominant HB threshold.

Firstly, there are {\em asynchronous modes} on
the ring with an inactive, or quiescent, pacemaker cell for which
$\v{c}_j=\left(\v{q}_j,0\right)^{T}$ for $j\in\{1,\ldots,N-1\}$, where
$\v{q}_j$ is defined in \eqref{num_ident:qj_def} and satisfies
$\v{q}_j^{T}\v{e}=0$ for $j\in\{1,\ldots,N-1\}$. For these modes, we
readily find that $\lambda$ satisfies
\begin{equation}\label{ringc:asynch}
  m_c(\alpha_c,\lambda)-4\pi \eps \gamma_j(\lambda) =0 \,,
  \quad j\in\{1,\ldots,N-1\}\,,
\end{equation}
where $\gamma_j(\lambda)$ is the eigenvalue of ${\mathcal G}_{\lambda}$ as
given in \eqref{num_ident:qj_def}.  As such, for Case I where this HB
threshold sets the instability threshold, the ${\mathcal O}(\eps)$
shift in $\alpha_c$ for these asynchronous modes is given in
Proposition \ref{prop:stab_circ} with the modification that we now
seek the maximum only over $j\in \lbrace{1,\ldots,N-1\rbrace}$, so that
we define
\begin{equation}\label{ringc:alpha1_async}
  \alpha_{1\mathrm{async}} \equiv \max_{j\in\lbrace{1,\ldots,N-1\rbrace}} \alpha_{1j}\,.
\end{equation}
Moreover, in using \eqref{num_ident:alpha_jfinal} for $\alpha_{1j}$,
the steady-state $u_{c1}$ is now the solution of
\eqref{ringc:equil_b}, which has an ${\mathcal O}(\eps)$ modification
due to the pacemaker cell.

Next, we consider modes that are synchronous for the ring cells, for which
the corresponding eigenvector has the form
$\v{c}=\left(\v{e},\beta\right)^T$, where $\beta$ is to be found. From
\eqref{ringc:gcep}, we use $\v{e}^T\v{e}=N$ and
${\mathcal G}_{\lambda}\v{e}=\gamma_N\v{e}$, with
$\gamma_N= \gamma_N(\lambda)$ given in \eqref{num_ident:sync_eig}, to
obtain that
\begin{equation*}
 \left( {\begin{array}{c}
   m_c \v{e} \\
           m_p \beta
 \end{array} } \right) = 4\pi \eps  \left( {\begin{array}{c}
   {\mathcal G}_{\lambda}\v{e} +\beta a_{\lambda} \v{e} \\
     a_{\lambda} \v{e}^T \v{e} + R_{\lambda} \beta
  \end{array} } \right) =  4\pi \eps \left( {\begin{array}{c}
       (\gamma_N + \beta a_{\lambda}) \v{e} \\
        a_{\lambda} N + R_{\lambda} \beta
 \end{array} } \right)\,,
\end{equation*}
which can be written as the $2\times 2$ homogeneous matrix problem
\begin{equation}\label{ringc:amat}
  {\mathcal A} \left( {\begin{array}{c}
                         1 \\
                         \beta
                       \end{array}}\right) =
                   \left( {\begin{array}{c}
                         0 \\
                         0
                           \end{array}}\right)  \qquad
                       \mbox{where} \qquad
                       {\mathcal A}\equiv \left( {\begin{array}{cc}
                m_c-4\pi \eps \gamma_N & -4\pi \eps a_{\lambda} \\
                  -4\pi \eps a_{\lambda} N & m_p - 4\pi \eps R_{\lambda}
                                                  \end{array}}\right)\,.
\end{equation}                                              
This problem has a nontrivial solution iff $\mathrm{det}{\mathcal A}=0$,
which yields that $\lambda$ must satisfy ${\mathcal N}(\lambda)=0$,
where
\begin{equation}\label{ringc:Nlam}
  {\mathcal N}(\lambda)\equiv m_{c} m_{p} - 4\pi \eps \left(\gamma_N m_p
    +R_{\lambda} m_c\right) + 16\pi^2\eps^2 \left(\gamma_N R_{\lambda}-
    a_{\lambda}^2 N \right) \,,
\end{equation}
with $m_c=m_c(\lambda,\alpha_c)$ and $m_p=m_p(\lambda,\alpha_p)$ as defined in
\eqref{ringc:mc_mp}. At a root of ${\mathcal N}(\lambda)=0$, we obtain
that $\beta$ in \eqref{ringc:amat} can be written in the two equivalent forms
\begin{equation}\label{ringc:beta}
  \beta = \frac{m_c}{4\pi \eps a_{\lambda}} - \frac{\gamma_N}{a_{\lambda}} 
  \qquad \mbox{or} \qquad \beta= \frac{4\pi \eps a_{\lambda} N}{m_p-4\pi \eps
    R_{\lambda}} \,.
\end{equation}

We now analyze Hopf bifurcations for this {\em synchronous mode}
$\v{c}=\left(\v{e},\beta\right)^T$ for our two distinct cases. In Case
I, where the HB associated with the ring cells sets the instability
threshold, so that $m_c(i\lambda_{0I},\alpha_{cH})=0$, we expand
$\lambda=i\lambda_{0I}+\eps \lambda_1+\ldots$ and
$\alpha_c=\alpha_{cH}+\eps\alpha_{1\mathrm{sync}}+\ldots$. Upon substituting
this expansion in ${\mathcal N}(\lambda)=0$ from \eqref{ringc:Nlam} and
collecting terms of ${\mathcal O}(\eps)$, we set $\mathrm{Re}(\lambda_1)=0$
to find that
\begin{equation}\label{ringc:alpha1_sync}
  \alpha_{1\mathrm{sync}}=4\pi \frac{\mathrm{Re}\left[{\gamma_N(i\lambda_{0I})/
        \partial_{\lambda}   m_{c}^{0}}\right]}{
    \mathrm{Re}\left[{\partial_{\alpha_c}m_{c}^{0}/
        \partial_{\lambda}m_{c}^{0}}\right]}\,, \quad
     \v{c}=       \left( {\begin{array}{c}
                         \v{e} \\
                         \beta
                          \end{array}}\right) \,, \quad
  \beta = 4\pi \eps N \Big{\vert}\frac{a_{\lambda_0}}{m_p^{0}}
  \Big{\vert} e^{i\mathrm{Arg}\left({a_{\lambda_0}/m_p^{0}}\right)} + {\mathcal O}(\eps)
  \,,
\end{equation}
where $m_p^{0}\equiv m_p(i\lambda_{I0},\alpha_p)$ and
$a_{\lambda_0}\equiv a_{\lambda}(i\lambda_{I0})$. In the expression
for $\alpha_{1\mathrm{sync}}$ in \eqref{ringc:alpha1_sync}, the terms
can be evaluated using \eqref{num_ident:ratio}, to obtain
\begin{equation}\label{ringc:alpha1_sync_ex}
  \alpha_{1\mathrm{sync}} = \frac{\left(\alpha_{cH}+u_{c1}^2\right)^{2}
    d_{1c}\kappa_c^2}{
    2\mu_c u_{c1} d_{2c} D} \left(\frac{1+\kappa_c + \mathrm{det}{J_c}
    }{ 1+\kappa_c + 2 \mathrm{det}{J_c} } \right) \mathrm{Re}\left[
      \gamma_{N}(i\lambda_{0I}) \left(-1 + i \frac{\lambda_{0I}}
      {1+\kappa_c}\right) \right] \,,
\end{equation}
where $\lambda_{0I}$ and $\mathrm{det}(J_c)$ are given in
\eqref{num_ident:lam0} in terms of the solution $u_{c1}$ to
\eqref{ringc:equil_b}.  Since $\beta={\mathcal O}(\eps)$ in
\eqref{ringc:alpha1_sync} we conclude that a HB for the synchronous
mode associated with the ring cells is able to trigger small amplitude
oscillations for the pacemaker cell.  From the form of the eigenvector
$\v{c}=\left(\v{e},\beta\right)^T$, the phase shift between the
signaling gradients near the identical ring cells and the pacemaker
cell at the HB point, as characterized by \eqref{nstabform:c1}, is
determined by $\mathrm{Arg} \left({a_{\lambda_0}/m_p^{0}}\right)$.

\begin{figure}[htbp]
  \centering
    \begin{subfigure}[b]{0.48\textwidth}
      \includegraphics[width=\textwidth]{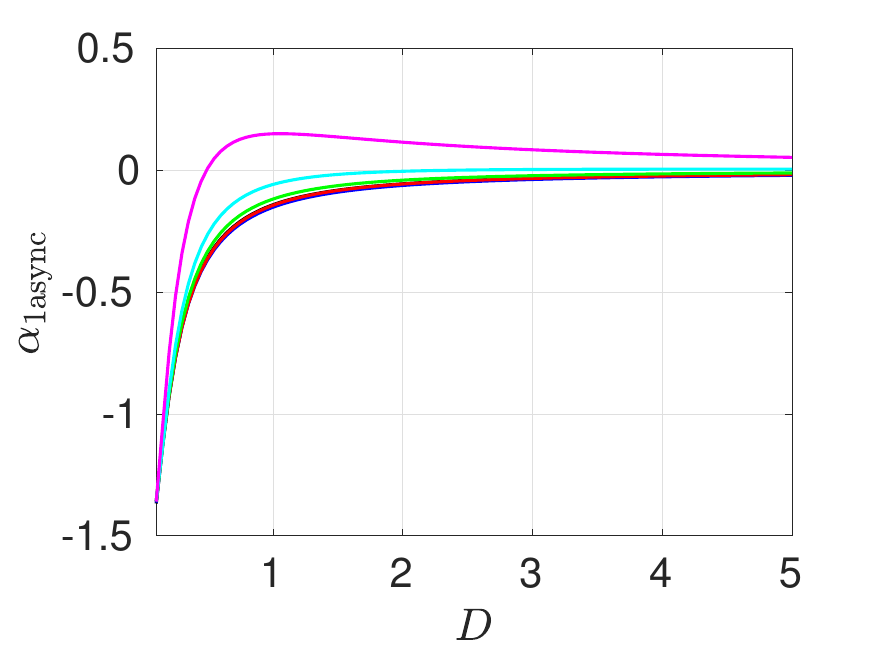}
      \caption{Case I: $\v{c}=(\v{q}_j,0)^T$ (asynchronous)}
      \label{fig:caseI_alpha_async}
    \end{subfigure}
      \begin{subfigure}[b]{0.48\textwidth}
      \includegraphics[width=\textwidth]{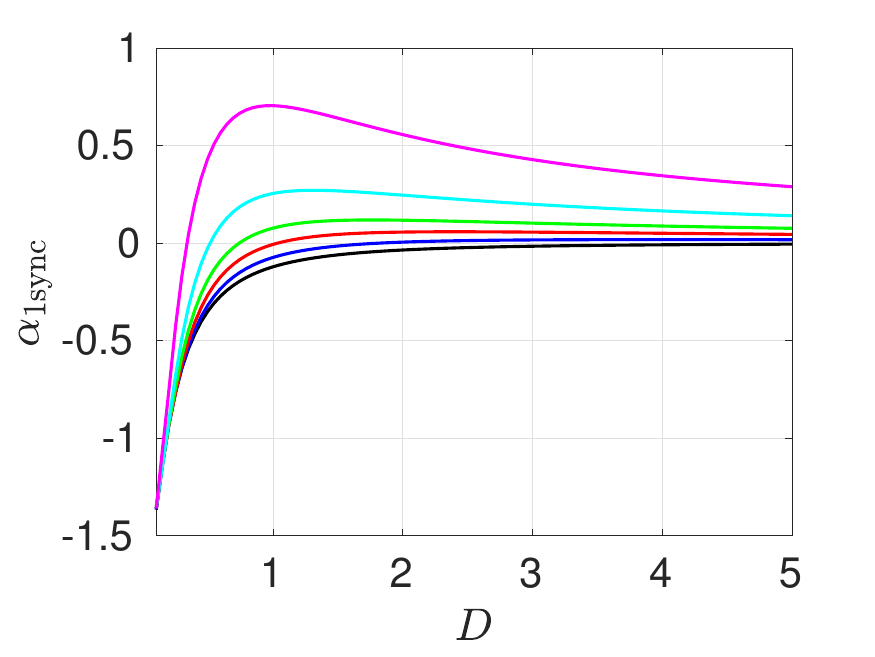}
        \caption{Case I: $\v{c}=(\v{e},\beta)^T$} \label{fig:caseI_alpha_sync}
    \end{subfigure}
    \caption{The ${\mathcal O}(\eps)$ correction to the leading-order
      HB threshold $\alpha_{cH}$ for Case I, with parameters as in
      Fig.~\ref{fig:caseI:HB}, for six different cell configurations.
      Left: $\alpha_{1\mathrm{async}}$ in \eqref{ringc:alpha1_async}
      versus $D$ for the asynchronous modes
      $\v{c}=\left(\v{q}_j,0\right)^T$, for which $\v{q}_j^{T}\v{e}=0$
      for $j\in \lbrace{1,\ldots,N-1\rbrace}$. Right:
      $\alpha_{1\mathrm{sync}}$ in \eqref{ringc:alpha1_sync_ex} versus
      $D$ for the synchronous mode
      $\v{c}=\left(\v{e},\beta\right)^T$. The six cell configurations
      are: $N=2$ antipodal (black curve), $N=4$ tetrahedron (blue),
      $N=6$ octahedron (red), $N=8$ cube-shaped (green), $N=12$
      icosahedron (cyan), and $N=20$ dodecahedron (magenta). Since
      $\alpha_{1\mathrm{sync}}>\alpha_{1\mathrm{async}}$, it is the
      synchronous mode that first develops intracellular oscillations
      through a HB as $\alpha_c$ decreases.}
      \label{fig:caseI_alpha}
  \end{figure}

In Fig.~\ref{fig:caseI_alpha_async} and
Fig.~\ref{fig:caseI_alpha_sync} we plot $\alpha_{1\mathrm{async}}$ and
$\alpha_{1\mathrm{sync}}$ versus $D$, as defined in
\eqref{ringc:alpha1_async} and \eqref{ringc:alpha1_sync_ex}, for the
Case I parameters as given in the caption of Fig.~\ref{fig:caseI:HB}
for our six cell configurations. We observe from
Fig.~\ref{fig:caseI_alpha_async} that for small $D$ small we have
$\alpha_{1\mathrm{async}}<0$ for all the cell configurations. As such,
when $D$ is small, the range of $\alpha_c$ where asynchronous
instabilities occur is ${\mathcal O}(\eps)$ smaller than that predicted by
the leading-order theory.  By comparing
Fig.~\ref{fig:caseI_alpha_async} and Fig.~\ref{fig:caseI_alpha_sync},
where we observe $\alpha_{1\mathrm{sync}}>\alpha_{1\mathrm{async}}$,
we predict that it is the synchronous mode
$\v{c}=\left(\v{e},\beta\right)^T$ that first becomes destabilized to
an oscillatory instability as $\alpha_c$ decreases. As expected, among
our six cell configurations, Fig.~\ref{fig:caseI_alpha_sync} shows
that the dodecahedron configuration, which has the largest
number of cells, has the largest range of $\alpha$ where the ring cells
have synchronized intracellular oscillations.

  Alternatively, for Case II where the HB instability threshold is set
  by the pacemaker cell, so that $m_p(i\lambda_{0I},\alpha_{pH})=0$,
  we expand $\lambda=i\lambda_{0I}+\eps \lambda_1+\ldots$ and
  $\alpha_p=\alpha_{pH}+\eps\alpha_{p1}+\ldots$ in
  ${\mathcal N}(\lambda)=0$, as given in \eqref{ringc:Nlam}, and set
  $\mathrm{Re}(\lambda_1)=0$. At ${\mathcal O}(\eps)$, we calculate that
\begin{equation}\label{ringc:alphap_dom}
  \alpha_{p1}=4\pi \frac{\mathrm{Re}\left[ {R_{\lambda_0}/
        \partial_{\lambda}   m_{p}^{0}}\right]} {
    \mathrm{Re}\left[{\partial_{\alpha_p}m_{p}^{0}/
        \partial_{\lambda}m_{p}^{0}}\right]}\,, \quad
     \v{c}=       \left( {\begin{array}{c}
                         \v{e} \\
                         \beta
                          \end{array}}\right) \,, \quad
  \beta = \frac{1}{4\pi \eps} \Big{\vert}\frac{m_c^{0}}{a_{\lambda_0}}
  \Big{\vert} e^{i\mathrm{Arg}\left( {m_{c}^{0}/a_{\lambda_0}}\right)} +
  {\mathcal O}(1)    \,,
\end{equation}
where $m_c^{0}\equiv m_c(i\lambda_{I0},\alpha_c)$ and
$R_{\lambda_0}\equiv (4\pi)^{-1}\sqrt{(\sigma+i\lambda_{0I})/D}$. In
analogy with \eqref{num_ident:ratio}, we have at the leading-order HB
threshold that
\begin{equation}\label{num_ident:ringc:ratio}
  \begin{split}
    \mathrm{Re}{\left[\frac{\partial_{\alpha_p}m_{p}^{0}}{
          \partial_{\lambda}m_{p}^{0}}\right]} &=
  \frac{\mu_c u_{p1}}{\left(\alpha_{pH}+u_{p1}^2\right)^2}\left( 1 +
    \frac{\mathrm{det}{J_p}}{1+\kappa_p + \mathrm{det}{J_p}}\right)>0 \,,\\
\partial_{\lambda} m_{p}^{0} &= \frac{8\pi D d_{2p}}{d_{1p}\kappa_p^2}
    \left(-1 + \frac{i\lambda_{0I}}{1+\kappa_p}\right)^{-1} \,,
\end{split}
\end{equation}
where $\lambda_{0I} = \sqrt{\zeta_c(1+\kappa_p) \mathrm{det}{J_p}}$
and $\mathrm{det}(J_p)=\alpha_{pH}+u_{p1}^2$. Here $u_{p1}$ is the
solution to \eqref{ringc:equil_b} and $\alpha_{pH}$ is given in
\eqref{ringc:alpha_pH}.  Finally, upon substituting
\eqref{num_ident:ringc:ratio} into \eqref{ringc:alphap_dom}, we obtain
the explicit perturbed threshold
\begin{equation}\label{ringc:alphap_dom_ex}
  \alpha_{p1} = \frac{\left(\alpha_{pH}+u_{p1}^2\right)^{2} d_{1p}\kappa_p^2}{
    8\pi D^{3/2} d_{2p}\mu_c u_{p1}} \left(\frac{1+\kappa_p + \mathrm{det}{J_p}
    }{ 1+\kappa_p + 2 \mathrm{det}{J_p} } \right) \mathrm{Re}\left[
      \sqrt{\sigma + i \lambda_{0I}} \left(-1 + i \frac{\lambda_{0I}}
        {1+\kappa_p}\right) \right] <0 \,.
\end{equation}
Since $\alpha_{p1}<0$, we conclude for Case II, where the pacemaker cell
sets the HB threshold, that the range of $\alpha$ where oscillations occur
is smaller than predicted by the leading-order theory.

\begin{figure}[htbp]
  \centering
        \includegraphics[width=88mm]{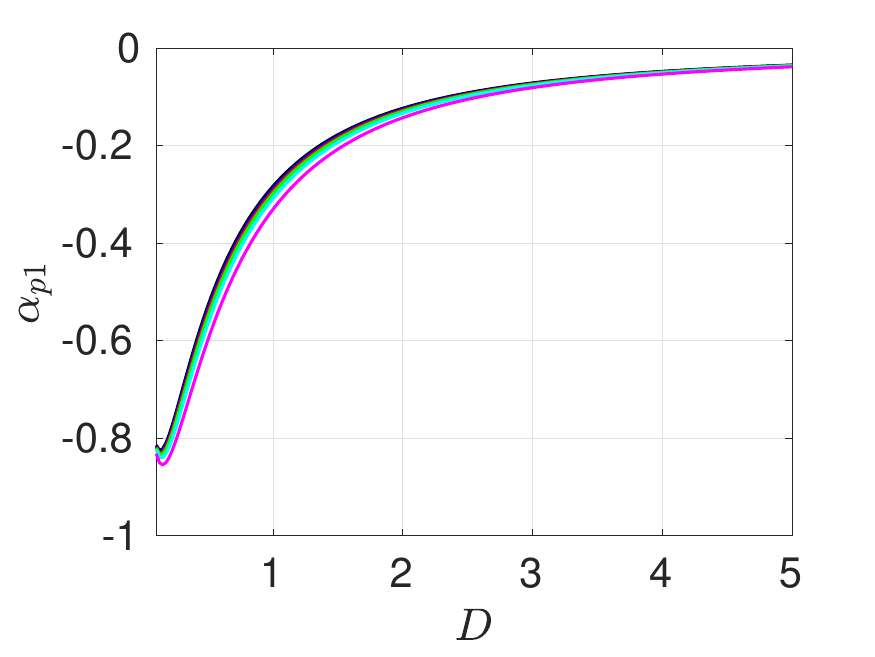}
        \caption{The ${\mathcal O}(\eps)$ correction $\alpha_{p1}$ to
          the leading-order HB threshold $\alpha_{pH}$ for Case II,
          with parameters in Fig.~\ref{fig:caseII:HB}, for the six
          different cell configurations described in the caption of
          Fig.~\ref{fig:caseI_alpha}. The corresponding eigenvector of
          the GCEP is $\v{c}=\left(\v{e}, \beta\right)^T$ with
          $\beta={\mathcal O}(\eps^{-1})$. There is only a very
          minimal effect on the perturbed threshold from the spatial
          configuration of the ring cells.}
        \label{fig:caseII_alpha}
\end{figure}

Since $|\beta|={\mathcal O} (\eps^{-1})$ in \eqref{ringc:alphap_dom}
it follows that the synchronous oscillations in the ring cells at
onset, which are triggered by the HB of the pacemaker cell, have a
much smaller amplitude as compared to that for the pacemaker cell. In
Fig.~\ref{fig:caseII_alpha} we plot $\alpha_{p1}$ for the six cell
configurations for the Case II parameters corresponding to
Fig.~\ref{fig:caseII:HB}. From this figure, we observe that the
perturbed threshold $\alpha_{p1}$ depends only very weakly on the
spatial configuration of ring cells. This weak dependence on the ring
cell arrangement arises from the ${\mathcal O}(\eps)$ influence of
$\gamma_{N0}$ in the determination of $u_{p1}$ in
\eqref{ringc:equil_b}.

The key qualitative implication of our analysis for Case II is that a
Hopf bifurcation in {\em a single pacemaker cell}, can initiate small
amplitude synchronous intracellular oscillations for the entire
collection of ring cells, which otherwise would be in an ambient state
without the pacemaker. 

\section{Time-Marching Scheme for the Integro-differential System}
\label{sec:algorithm}

In this section, and in Appendix \ref{app:numerics}, we derive a
numerical time-marching scheme to solve the asymptotic
integro-differential system \eqref{2:reduced} efficiently. In \S
\ref{sec:selkov_experiments} we apply our scheme to the case of
Sel'kov kinetics \eqref{isolated:selkov:ode} in order to confirm the
HB thresholds of \S \ref{sec:selkov}.

\begin{figure}
	\centering
	\includegraphics[width=\textwidth]{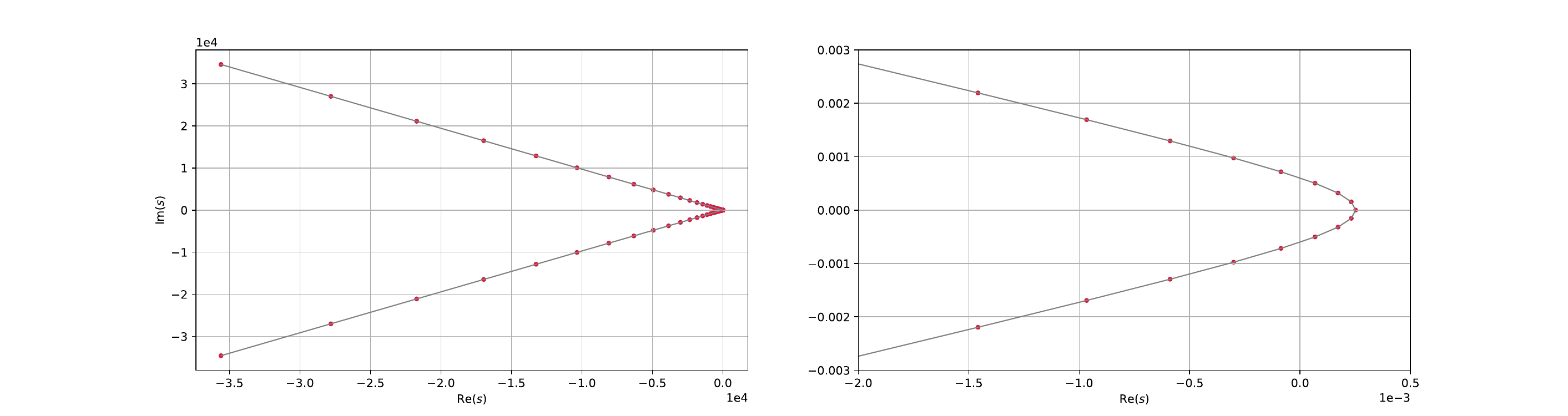}
	\caption{Left: Plot of the hyperbolically shaped contour $\Gamma$ in
        the Laplace transform plane defined by \eqref{dh:curve} with
        $\theta=0.95$, $\alpha=0.8$ and $\beta=0.7$ in
        \eqref{approx:cond}. The discretization points
        $\lbrace{s_{\ell}\rbrace}$ for the \emph{sum-of-exponentials}
        approximation of $g(t)$ and $G(a_{jk}, t)$ with $n=75$
        are the red dots.  Right: Zoomed plot showing the
        discretization near the origin.}
    \label{fig:laplace_discret}
\end{figure}

For $t>0$ we first represent the integral kernels $g(t)$
and $G(a_{jk},t)$ defined in \eqref{3:lap2} and \eqref{2:free_3d},
respectively, in the Laplace transform space along the vertical
Bromwich contour defined by $\Gamma_B\equiv \lbrace{s \, \vert
    \, \mathrm{Re}(s)=b>0\rbrace}$ as
\begin{equation*}
  g(t) = \frac{1}{2\pi i} \int_{\Gamma_B} \hat{g}(s) e^{st}\;ds\,,
  \qquad G(a_{jk}, t) = \frac{1}{2\pi i} \int_{\Gamma_B}
  \hat{G}(a_{jk}, s) e^{st}\;ds \,,
\end{equation*}
where the Laplace transforms are given by (see \eqref{3:lap2} and
\eqref{2:free_3d})
\begin{equation*}
  \hat{g}(s) \equiv \mathcal{L}[g](s) = \frac{\sqrt{s+\sigma}}{s},
  \qquad \hat{G}(a_{jk}, s) \equiv \mathcal{L}[G(a_{jk}, \cdot)](s) =
  \frac{e^{-a_{jk}\sqrt{s+\sigma}}}{4\pi D^{\frac{3}{2}} a_{jk}} \,,
\end{equation*}
with $a_{jk}\equiv {|\x_j-\x_k|/\sqrt{D}}$.  Since the only
singularities of $\hat{g}(s)$ and $\hat{G}(a_{jk},s)$ are branch cuts
along the portion of the negative real axis with $\mathrm{Re}(s)\leq
-\sigma$ and since $|\hat{g}(s)|$ and $|\hat{G}(a_{jk},s)|$ decay
sufficiently fast as $|s|\to\infty$ for $\mathrm{Re}(s)<0$, we proceed
as in \cite{greengard_2015} (see also \cite{beylkin_2005},
\cite{beylkin_2010}) by deforming the Bromwich contour $\Gamma_B$ to a
contour $\Gamma$, as depicted in Fig.~\ref{fig:laplace_discret}, where
the integrals are then discretized. This contour $\Gamma$ originates
from $\mathrm{Re}(s) \to -\infty$ in the third quadrant, intersects
the positive real axis, and returns to $\mathrm{Re}(s) \to -\infty$ in
the second quadrant.

More precisely, $\Gamma$ is chosen as the hyperbolically shaped curve
parameterized as
\begin{equation}\label{dh:curve}
  \Gamma\equiv \lbrace{\chi P(x) \,, \,\, x\in \R\rbrace} \,, \quad
  \mbox{where} \quad P(x) \equiv 1-\sin(\alpha + i\, x) \,,
\end{equation}
with $0<\alpha<{\pi/2}$, $\chi>0$, so that $P(0)>0$. The discretization
points $\{s_\ell\}$ along $\Gamma$ are chosen as $s_\ell = \chi P(x_\ell)$ with
$x_\ell= \ell h$ for $|\ell|\leq n$, so that
\begin{equation*}
  g(t) = \frac{1}{2\pi i} \int_{\Gamma_B} \hat{g}(s) e^{st}\;ds =
  \frac{\chi}{2\pi} \int_{-\infty}^{\infty} e^{\chi P(x) t}
  \hat{g}(\chi P(x)) \cos(\alpha + ix)\; dx \approx \sum_{\ell=-n}^n e_\ell
  e^{s_{\ell} t}\,,
\end{equation*}
for some coefficients $e_{\ell}$. Proceeding in a similar way for
$G(a_{jk}, t)$, we obtain their {\em sums-of-exponentials} (SOE) approximations
defined by
\begin{equation}\label{alg:soe}
  g(t) \approx g_n(t) \equiv \sum_{\ell=-n}^n e_\ell e^{s_\ell t}\,, \qquad
  G(a_{jk}, t) \approx G_n(a_{jk}, t) \equiv \sum_{\ell=-n}^n
  \omega_{jk\ell} e^{s_\ell t}\,,
\end{equation}
where the $e_{\ell}$ and $\omega_{jk\ell}$ are given by
\begin{equation}\label{alg:el_omega}
e_\ell = \frac{\chi h}{2\pi} \cos(\alpha + i\ell h)\hat{g}(s_\ell) \,, \qquad
\omega_{jk\ell} = \frac{\chi h}{2\pi} \cos(\alpha + i\ell h)
\hat{G}(a_{jk}, s_\ell) \,, \quad \mbox{for} \quad j\ne k \,.
\end{equation}
As shown in \cite{greengard_2015}, the approximation error made can be
well-controlled as characterized in the next lemma.

\begin{lemma}\label{prop:approx}(cf.~\cite{greengard_2015})
  Consider the time interval
    $0<\delta\leq t\leq T_f$, with $T_f\geq 1000\delta$, and let
    $\varepsilon_f$ with $0<\varepsilon_f<0.1$ be a prescribed
    error-tolerance. Then, for the choice of parameters $h$ and $\chi$
    defined by
    \begin{equation}\label{approx:cond}
      h=\frac{a(\theta)}{n} \,, \quad
      \chi=\frac{2\pi \beta n(1-\theta)}{T_f a(\theta)} \,, \quad
      \mbox{where} \quad
      a(\theta)\equiv \cosh^{-1}\left(\frac{2T_f}{\delta(1-\theta)
          \sin(\alpha)}\right) \,, \quad
    \end{equation}
    with $0<\alpha-\beta<\alpha+\beta<{\pi/2}$ and $0<\theta<1$, we have
    the uniform estimates
    \begin{equation}\label{approx:thm}
      \| g(t) - g_n(t)\|\leq \frac{\varepsilon_f}{\sqrt{t}} \quad
      \text{and} \quad \| G(a_{jk}, t) - G_n(a_{jk}, t)\|\leq
      \frac{\varepsilon_f}{\sqrt{t}}
      \quad \mbox{on} \quad \delta\leq t\leq T_f \,,
    \end{equation}
    when $n$ is sufficiently large of the order
    \begin{equation}\label{approx:n_cond}
   n = {\mathcal O} \left( \left(-\log\varepsilon_f + \log\log\left({T_f/\delta}
          \right) \right) \log\left({T_f/\delta}\right)\right) \,.
    \end{equation}
\end{lemma}

This result shows that the number of discretization points required
varies only logarithmically with the error tolerance and the
ratio of the bounds of the time interval.

We now illustrate how large $n$ must be to ensure specified
error-tolerances $\varepsilon_f$ for both $g_n$ and $G_n$. In order to
estimate $\sqrt{t} \, \| f(x, t) - f_a(x, t)\|$, where $f$ stands for
either $g$ (dropping the $x$-dependence) or $G$ in the following,
depending on the number $2n+1$ of terms in the
\emph{sum-of-exponentials} approximation, we proceed as in
\cite{greengard_2015}. We take a $49 \times 1000$ grid $(x_j, t_k)$
where $x_j=2^{-16+j}$ for $j\in \{1, ..., 49\}$ (not applicable for
$g(t)$), and $t_k = \delta e^{k\blacktriangle}$ with
$\blacktriangle\equiv{\log\left({T/\delta}\right)/999}$ for
$k\in\lbrace{0,\ldots, 999\rbrace}$ (equispaced on a logarithmic
scale). Note that the discretization point $x_0=0$ is a singular point
for $G(\cdot, t)$, so that we have excluded it, only keeping the 49
spatial discretization points above for Table \ref{tab:G_eps}. We then
compute $\max_{j, k} \, \sqrt{t} \, | f(x_j, t_k) - f_a(x_j, t_k)|$,
which provides the error results given in Tables \ref{tab:g_eps} and
\ref{tab:G_eps} on the three different intervals
$I_1 \equiv [\delta, T] = [10^{-3}, 1]$,
$I_2 \equiv [\delta, T] = [10^{-3}, 10^3]$ and
$I_3 \equiv [\delta, T] = [10^{-5}, 10^4]$. For two cells, it takes
less than a second to compute the coefficients $e_{\ell}$ and
$\omega_{jk\ell}$ for the interval $I_3$.

\begin{table}[h!]
  \centering
\begin{tabular}{ c c c c }
        \hline
 $\varepsilon_f$   &$n(I_1)$       &$n(I_2)$ & $n(I_3)$ \\
        \hline
$10^{-3}$   &$20  \,\, (5.657\cdot 10^{-4})$    &$43 \,\, (8.691\cdot 10^{-4})$
&$65  \,\, (9.220\cdot 10^{-4})$ \\
$10^{-6}$   &$36  \,\, (8.029\cdot 10^{-7})$    &$69 \,\, (9.414\cdot 10^{-7})$
&$101  \,\, (9.101\cdot 10^{-7})$ \\
$10^{-9}$   &$53  \,\, (7.416\cdot 10^{-10})$   &$96 \,\, (8.611\cdot 10^{-10})$
&$137 \,\, (9.690\cdot 10^{-10})$ 
\end{tabular}
\caption{The number $2n+1$ of terms needed to approximate the kernel
  $g(t) = e^{-\sigma t}/\sqrt{\pi t} + \sqrt{\sigma}
  \;\mathrm{erf}\left(\sqrt{\sigma t}\right)$ for $\sigma=1$ with
  Laplace transform $\mathcal{L}[g(t)]= \sqrt{s+\sigma}/s$ to a
  precision $\varepsilon_f$ corresponding to \eqref{approx:thm} for
  the three time intervals $I_1 = [10^{-3}, 1]$,
  $I_2 = [10^{-3}, 10^3]$ and $I_3 = [10^{-5}, 10^4]$.}
\label{tab:g_eps}
\end{table}

\begin{table}[h!]
  \centering
\begin{tabular}{ c c c c }
        \hline
 $\varepsilon_f$   &$n(I_1)$       &$n(I_2)$ & $n(I_3)$ \\
        \hline
$10^{-3}$   &$46  \,\, (7.237\cdot 10^{-4})$    &$73 \,\, (7.785\cdot 10^{-4})$
&$111  \,\, (9.797\cdot 10^{-4})$ \\
$10^{-6}$   &$62  \,\, (8.469\cdot 10^{-7})$    &$99 \,\, (7.986\cdot 10^{-7})$
&$147  \,\, (9.047\cdot 10^{-7})$ \\
$10^{-9}$   &$78  \,\, (9.350\cdot 10^{-10})$   &$125 \,\, (7.867\cdot 10^{-10})$
&$183 \,\, (8.181\cdot 10^{-10})$ 
\end{tabular}
\caption{The number $2n+1$ of terms needed to approximate the 3-D
  Green function kernel
  $G(a_{jk}, t) = \mathrm{exp}(-\sigma t -a_{jk}^2/(4t))/(4\pi D
  t)^{3/2}$ with Laplace transform
  $\mathcal{L}[G(a_{jk}, \cdot)](s) =
  \mathrm{exp}(-a_{jk}\sqrt{s+\sigma}) / (4\pi D^{\frac{3}{2}}
  a_{jk})$ to a precision $\varepsilon_f$ corresponding to
  \eqref{approx:thm} for the three time intervals
  $I_1 = [10^{-3}, 1]$, $I_2 = [10^{-3}, 10^3]$ and
  $I_3 = [10^{-5}, 10^4]$. Parameters: $D=1$ and $\sigma=1$. This
  table for $\sigma=1$ is related to Table 1 in \cite{greengard_2015}
  for $\sigma=0$, but with a few additional explicit error values.}
\label{tab:G_eps}
\end{table}

By using this SOE approximation for $g$ and $G$ in \eqref{alg:soe},
together with Duhamel's principle, in Appendix \ref{app:numerics_cj} and
\ref{app:numerics_dj} we derive an efficient  time-update
scheme for the two nonlocal terms in the integro-differential system
\eqref{2:reduced}, which we define as
\begin{equation}\label{alg:CD:def}
  C_j(t) \equiv \int_0^t u_{1j}^{\prime}(\tau) g(t-\tau)\, d\tau \,,
\qquad D_{jk}(t) \equiv \int_0^t u_{1k}(\tau) G(a_{jk},t-\tau) \,
 d\tau \,.
\end{equation}
By using these update schemes, in Appendix \ref{app:numerics_ode} 
we formulate our time-marching algorithm to compute solutions to
\eqref{2:reduced} over long-time intervals.

\section{Time-Marching Results: Sel'kov Kinetics}\label{sec:selkov_experiments}

By using the numerical time-marching scheme derived in Appendix
\ref{app:numerics} we now present numerical solutions for our
integro-differential system for various cell arrangements and
parameter values, while fixing $\eps=0.03$, $\sigma=0.5$ and the
intracellular parameters $\zeta=0.15$ and $\mu=2.0$ for each cell. We
first consider the simplest case of two diffusively coupled cells. We
then investigate polyhedral configurations in which cells are centered
at the vertices of either an icosahedron or a dodecahedron, each
touching the unit sphere. Both polyhedral arrangements
are examined with and without an additional pacemaker cell at the
center.  Regarding parameter values, we primarily vary the
intracellular kinetic parameter $\alpha$ and the bulk diffusivity $D$
of the extracellular species $U$, which interacts directly with the
intracellular species $u_{ij}$ across the boundary of the
$j^{\text{th}}$ cell. We remark that, since $\zeta=0.15\ll 1$, the
intracellular ODEs are a slow-fast system and the figures below
demonstrate that oscillatory dynamics of $u_{1j}$ resemble a
relaxation oscillation, characteristic of an excitable medium in which
$u_{2j}$ acts as the trigger that excites $u_{1j}$ before subsequently
relaxing.

\begin{figure}[ht!]
\centering
\begin{subfigure}{0.49\textwidth}
    \includegraphics[width=\textwidth]{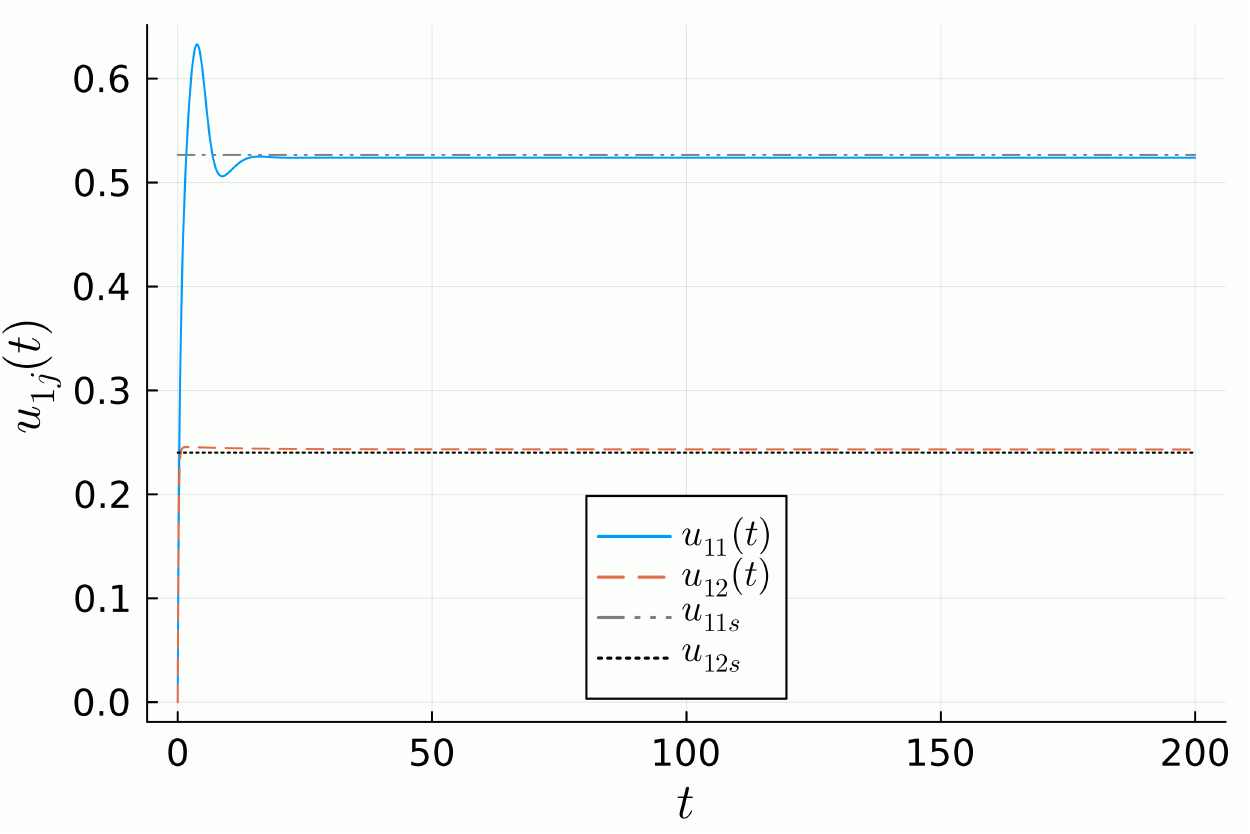}
\end{subfigure}
\hfill
\begin{subfigure}{0.49\textwidth}
    \includegraphics[width=\textwidth]{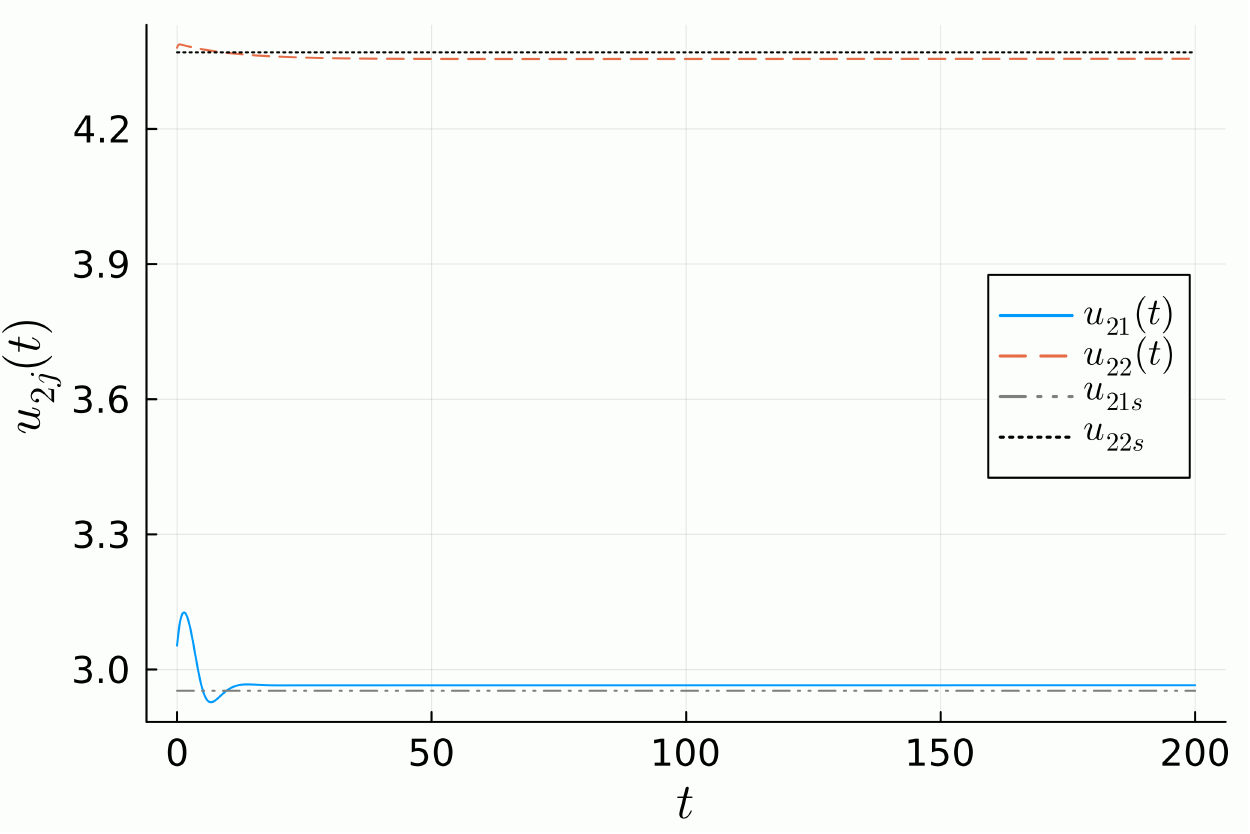}
\end{subfigure}
\begin{subfigure}{0.49\textwidth}
    \includegraphics[width=\textwidth]{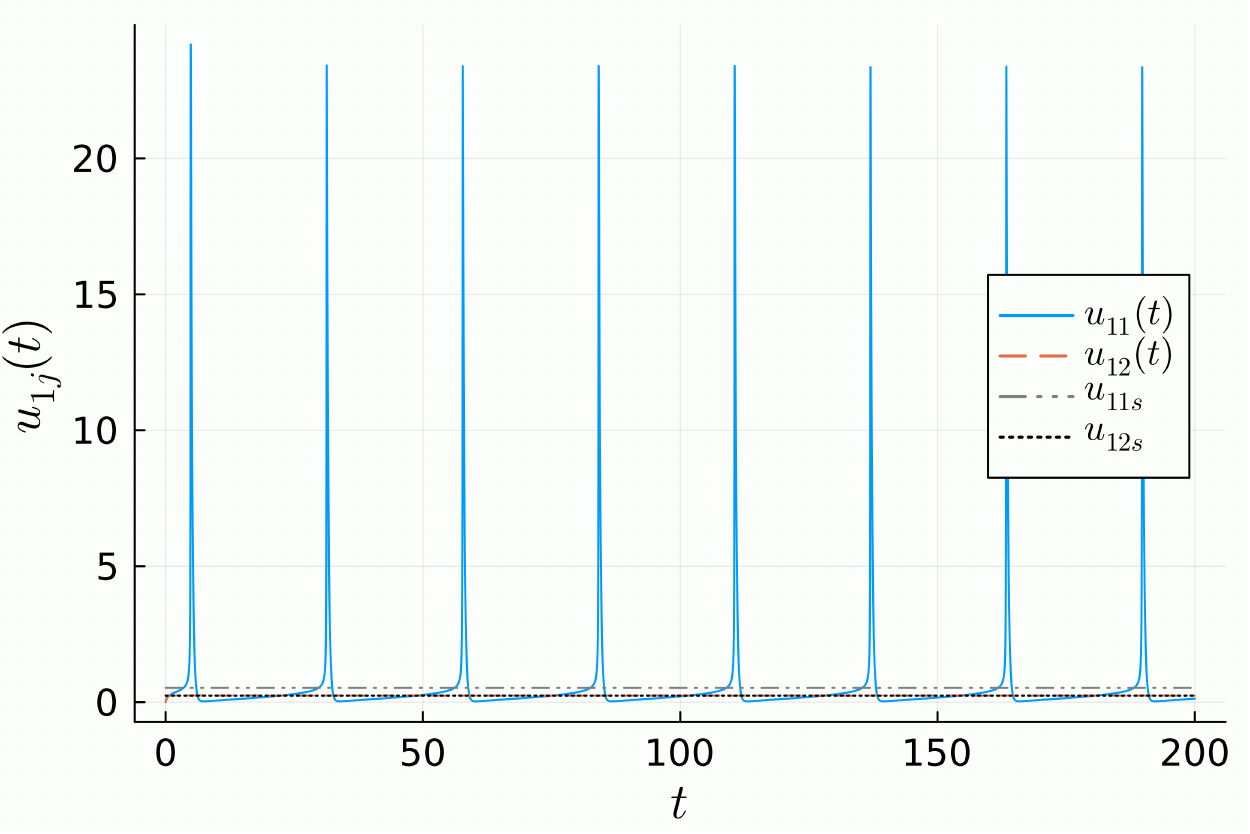}
\end{subfigure}
\hfill
\begin{subfigure}{0.49\textwidth}
    \includegraphics[width=\textwidth]{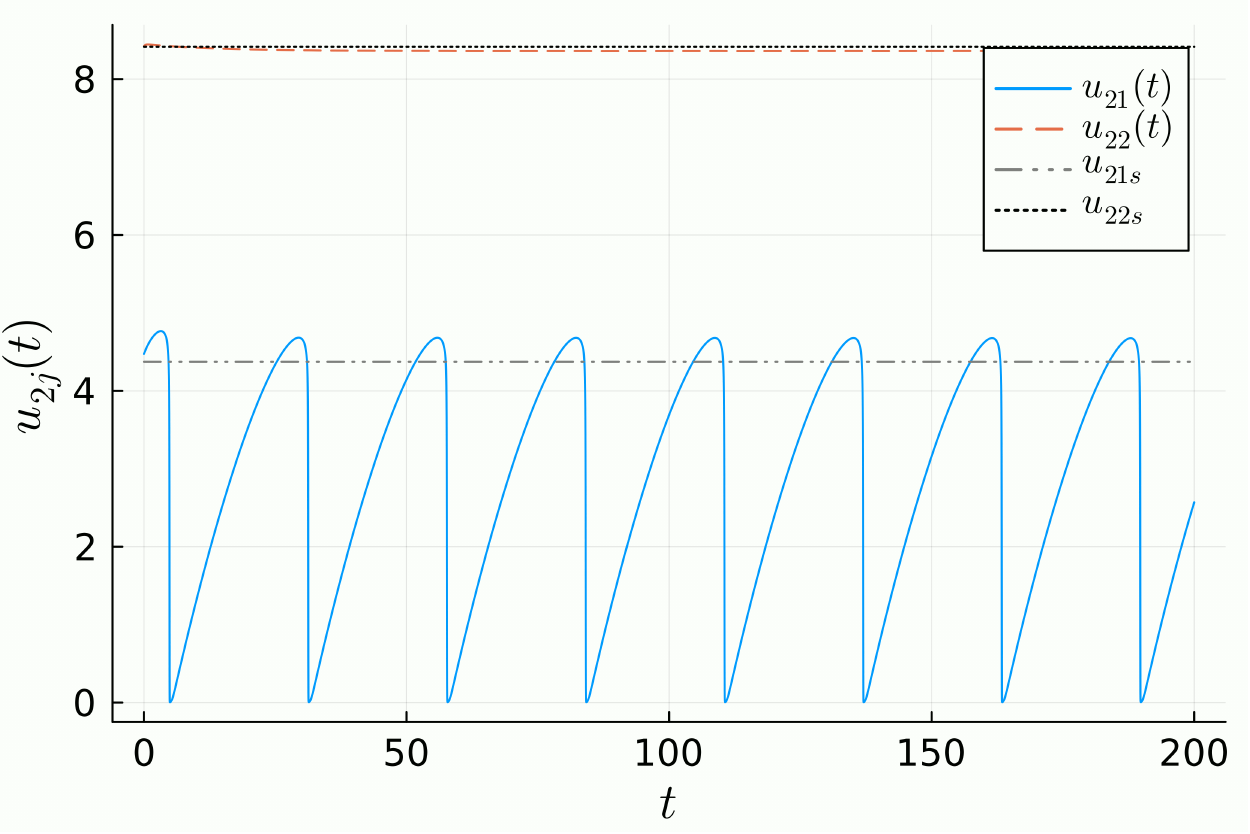}
\end{subfigure}
\caption{Dynamics for the two-cell arrangement with permeabilities
  $d_{11}=0.8$ and $d_{21}=0.4$ for cell 1 and $d_{12}=0.4$ and
  $d_{22}=0.8$ for cell 2. Top row: $\alpha = 0.4$ and $D=1$. Bottom
  row: $\alpha = 0.18$ and $D=1$. As predicted by the asymptotic
  theory (cf. Figs.~\ref{fig:two:selkov} and
  \ref{fig:twocell_alpha1}), cell 1 oscillates strongly for
  $\alpha=0.18$ while cell 2 does not. Neither cell oscillates for
  $\alpha=0.4$. All remaining parameters are as in the caption of
  Fig.~\ref{fig:two:selkov}. The initial condition is $u_{1j} = 0$
  for $j\in\{1,2\}$ and $(u_{21}, u_{22}) = (u_{21s}, u_{22s}) +
  \frac{1.1}{20} \phi_i + \frac{0.9}{20}\phi_a = (u_{21s}, u_{22s}) +
  0.1 \cdot (1, 0.1)$ where $(u_{21s}, u_{22s})$ corresponds to the
  steady-state \eqref{4:limit_sys} and $\phi_i=(1, 1)$ is the in-phase
  perturbation while $\phi_a=(1, -1)$ is the anti-phase perturbation,
  spanning the space of all possible initial perturbations. The
  numerical time stepping and plotting for each figure row took about
  30 $\rm{s}$ on a personal laptop with specifications
  13\textsuperscript{th} Gen Intel\textsuperscript{\textregistered}
  Core\textsuperscript{TM} i7-13700HX $\times$ 16, 15.3\,GB RAM, and
  LINUX MINT 22.3. The time step was taken to be $\Delta t = 0.03$
  with fivefold refinement for the nonlinear terms between the first
  two time steps. The sum-of-exponentials method parameters were
  chosen as in the caption of Fig.~\ref{fig:laplace_discret}. The
  oscillatory dynamics of $u_{1j}$ resemble a relaxation oscillation
  of an excitable medium, as $u_{1j}$ is directly excited by
  $u_{2j}$. The profile of $u_{1j}$ with its tail-shapes is very
  similar to that of the firing in biological neural systems.}
\label{fig:twononidcells}
\end{figure}

\paragraph{Two non-identical cells} To illustrate the
  linear stability theory of \S \ref{sec:selkov}, we first consider
  the case of two diffusively coupled cells where cell 1 has the
  permeabilities $d_{11}=0.8$ and $d_{21}=0.4$, while cell 2 has the
  swapped values $d_{12}=0.4$ and $d_{22}=0.8$. The remaining
  parameters are the same for both cells and are as in the caption of
  Fig.~\ref{fig:two:selkov}. Fig.~\ref{fig:two:selkov} together with
  Fig.~\ref{fig:twocell_alpha1} show that the HB threshold in $\alpha$
  for cell 1 is approximately $\alpha =
  \alpha_{1H}+\varepsilon\alpha_1 \approx 0.27 + 0.03\cdot (-0.2) =
  0.264$ when $D=1$. From the time-marching results shown in
  Fig.~\ref{fig:twononidcells} when $\alpha=0.18$, we observe that
  cell 1 oscillates exhibits relaxation oscillations, while cell 2
  does not. This is in agreement with the theoretical results obtained
  in Figs.~\ref{fig:two:selkov} and ~\ref{fig:twocell_alpha1}. When
  $\alpha=0.4$, we observe from the numerical time-marching results
  that neither cell oscillates, as predicted by our linear stability
  results encoded in Fig.~\ref{fig:two:selkov}.

\begin{figure}[ht!]
\centering
\begin{subfigure}{0.49\textwidth}
    \includegraphics[width=\textwidth]{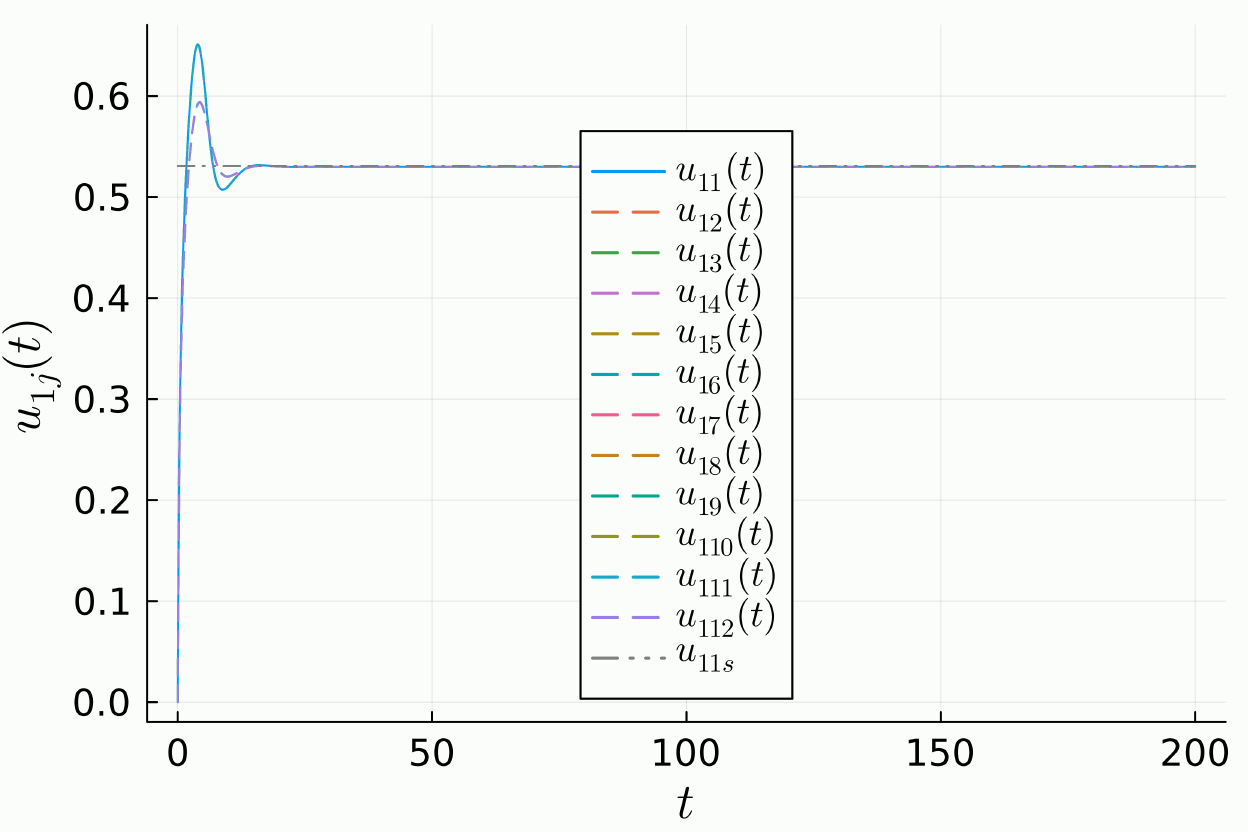}
\end{subfigure}
\hfill
\begin{subfigure}{0.49\textwidth}
    \includegraphics[width=\textwidth]{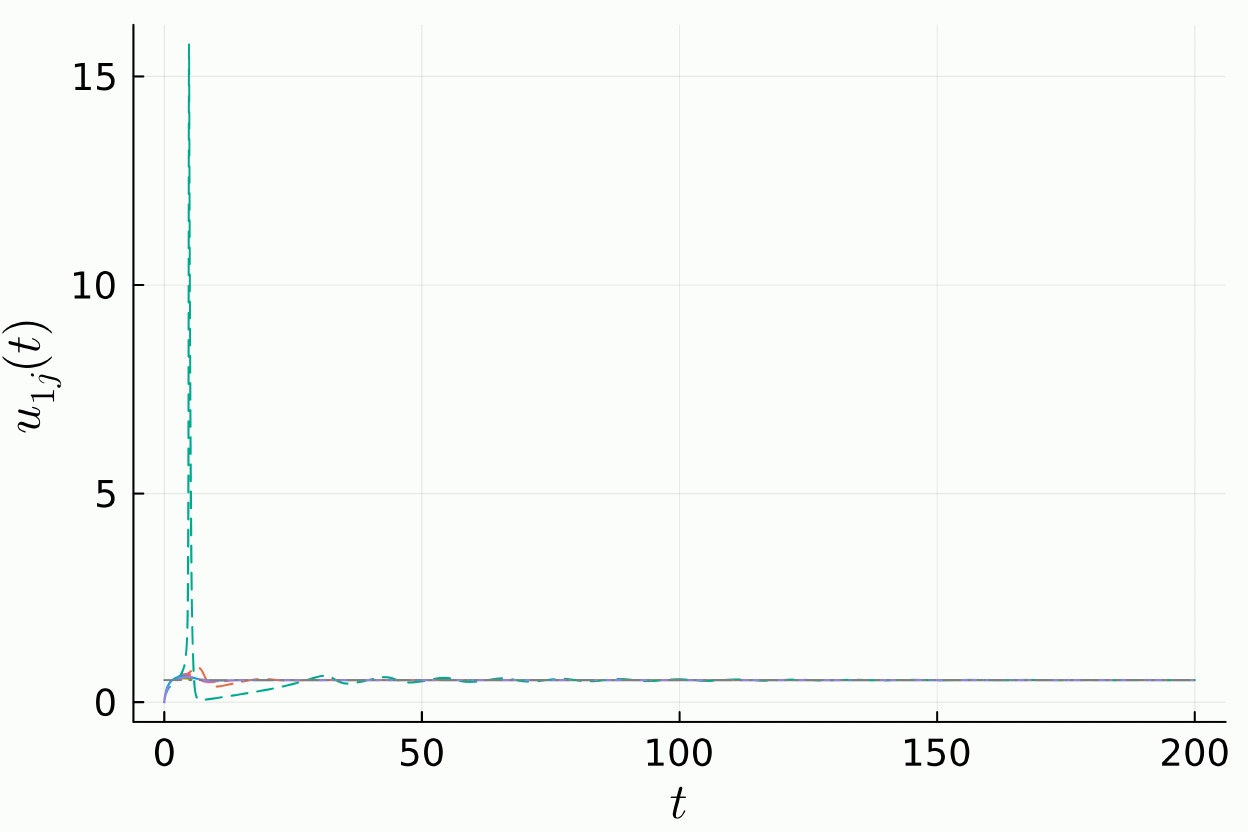}
\end{subfigure}
\begin{subfigure}{0.49\textwidth}
    \includegraphics[width=\textwidth]{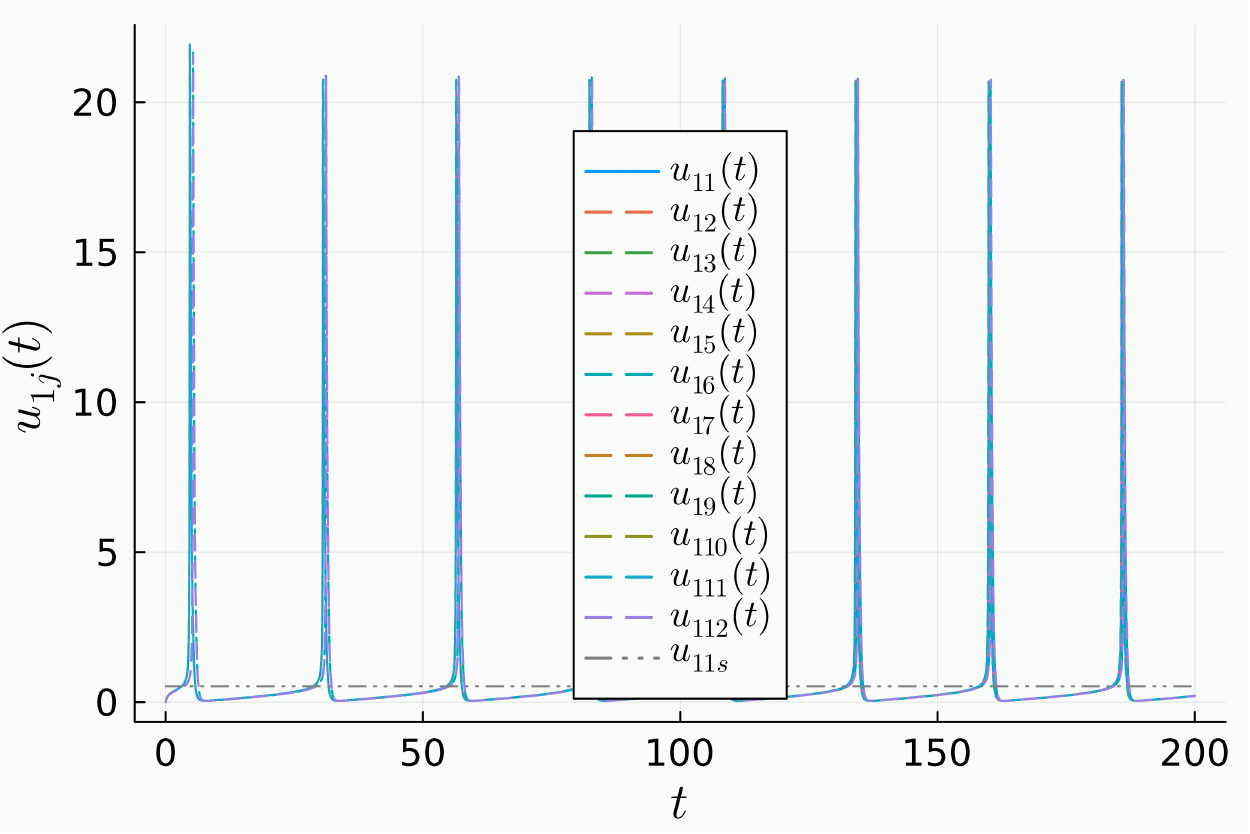}
\end{subfigure}
\hfill
\begin{subfigure}{0.49\textwidth}
    \includegraphics[width=\textwidth]{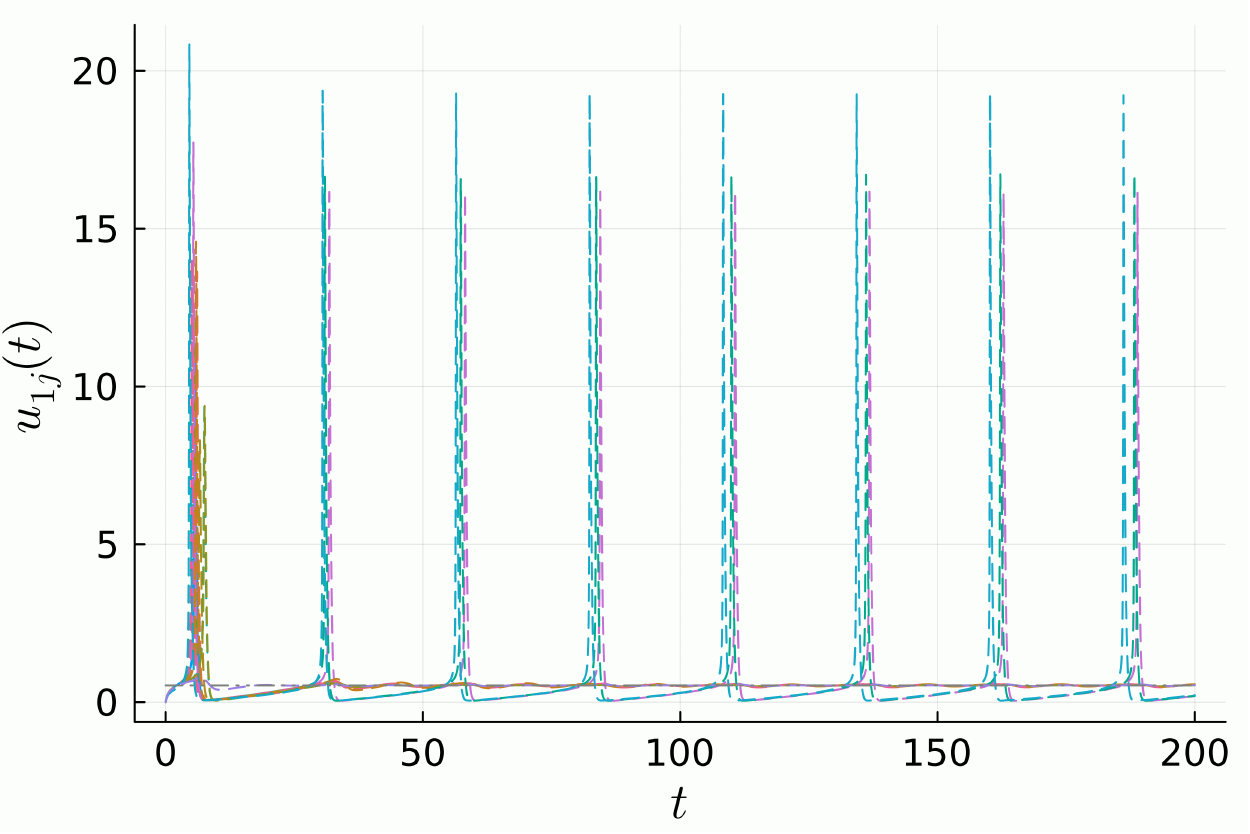}
\end{subfigure}
\begin{subfigure}{0.49\textwidth}
    \includegraphics[width=\textwidth]{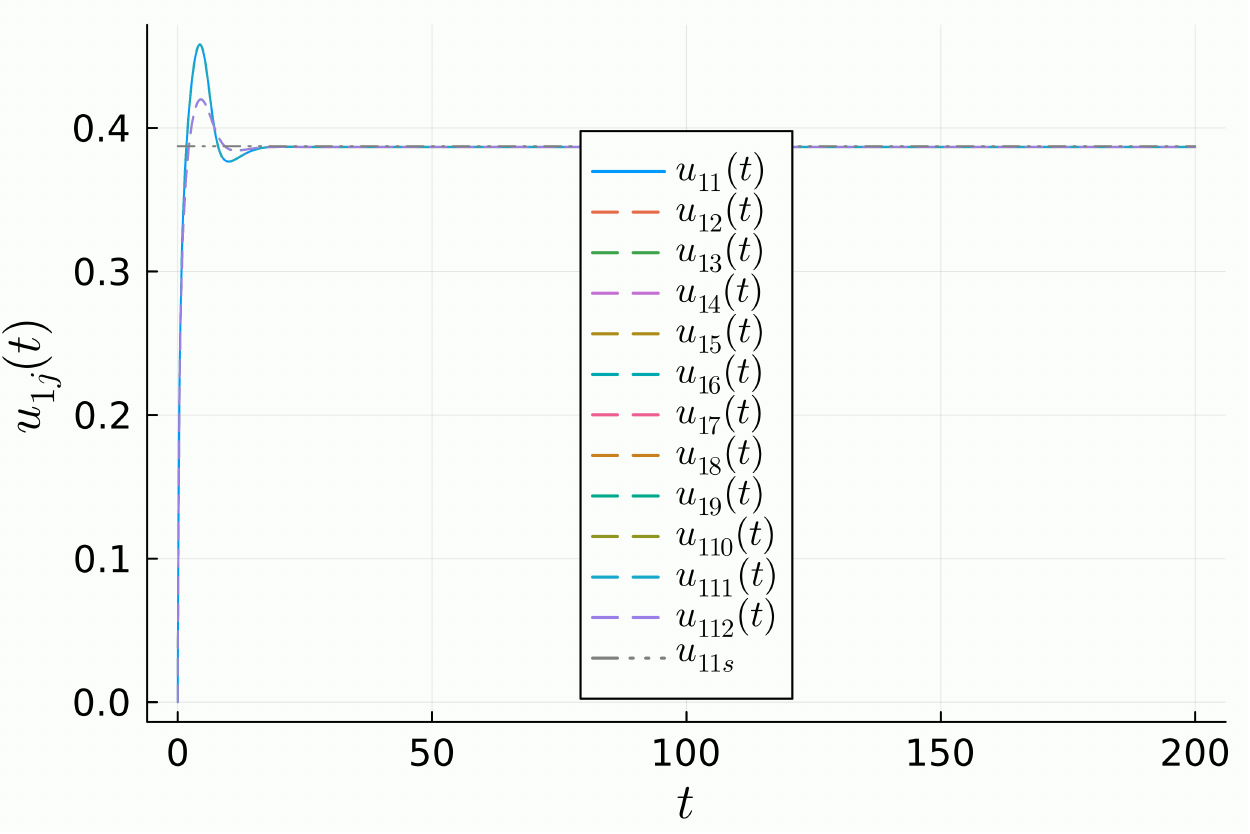}
\end{subfigure}
\hfill
\begin{subfigure}{0.49\textwidth}
    \includegraphics[width=\textwidth]{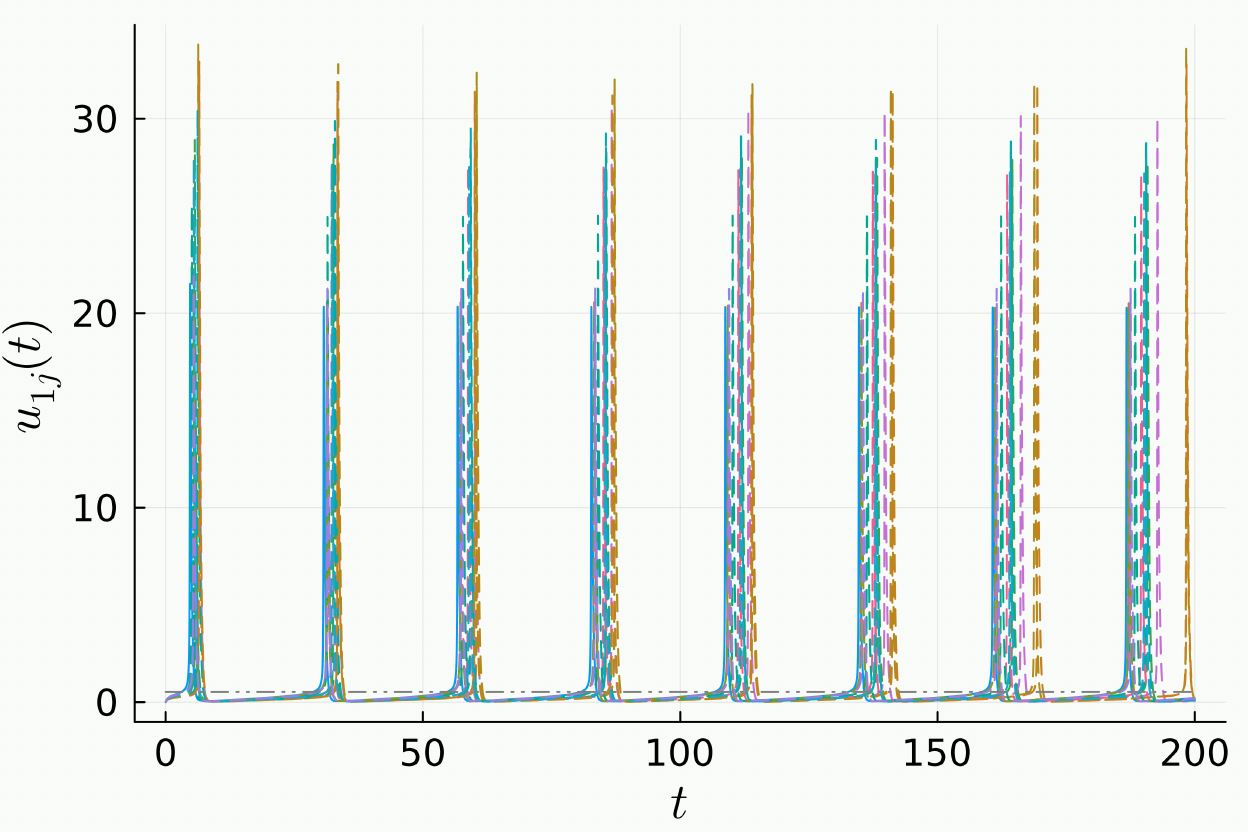}
\end{subfigure}
\caption{Dynamics of cells centered at the vertices of an
    icosahedron that is inscribed in the unit sphere each with
    permeabilities $d_{1c}=0.8$ and $d_{2c}=0.4$. Left column: When
    $\alpha = 0.4$ and $D = 1$ (top), identical cells converge to
    their coupled steady-state. For $\alpha = 0.2$ and $D = 1$,
    identical cells oscillate with increasing in-phase synchrony.
    This is consistent with the asymptotic theory leading to
    Fig.~\ref{fig:ident_icosa} together with the cyan curve of
    Fig.~\ref{fig:ident_alpha1}, for which the HB
    threshold is approximately $\alpha =
    \alpha_{cH}+\varepsilon\alpha_1 \approx 0.27 + 0.03\cdot 0.25 =
    0.2775$. Increasing $D$ to $D=4$ while fixing $\alpha=0.2$, the
    cells again converge to their coupled steady-state since
    $\alpha=0.2$ is now outside the oscillatory region of
    Fig.~\ref{fig:ident_icosa}.  Right column: Cells are slightly
    non-identical by choosing $\alpha = \alpha_j$ for each $j\in
    \lbrace{1,\ldots,N\rbrace}$ to be normally distributed according
    to $\alpha_j \sim \mathcal{N}(\alpha_{cH}+\varepsilon\alpha_1 +
    0.1, 0.05^2)$ (top; mean slightly outside the oscillatory
    parameter region), then $\alpha_j \sim
    \mathcal{N}(\alpha_{cH}+\varepsilon\alpha_1, 0.05^2)$ (center;
    mean on the bifurcation threshold), and lastly $\alpha_j \sim
    \mathcal{N}(\alpha_{cH}+\varepsilon\alpha_1 - 0.1, 0.05^2)$
    (bottom; mean inside the oscillatory parameter region). All other
    kinetic and numerical parameters are the same as for
    Fig.~\ref{fig:ident_icosa} and Fig.~\ref{fig:twononidcells},
    respectively. The initial condition is $u_{1j} = 0$ for
    $j\in\{1,...,N\}$ and ${\bf u}_2 = {\bf u}_{2s} + \frac{1.1}{20}
    \phi_i + \frac{0.9}{20}\phi_a$ where ${\bf u}_{2s}$ is the
    steady-state \eqref{4:limit_sys} and $\phi_i=(1, ..., 1)
    \in\R^{N}$ is the in-phase perturbation and $\phi_a=(1, -1, ...,
    (-1)^{N-1})$ is one of the anti-phase perturbations. The numerical
    time-stepping and plotting required $\sim 6 \text{ min}$ per
    sub-figure (hardware specifications as in the caption of
    Fig.~\ref{fig:twononidcells}).}
\label{fig:cellicosahedron}
\end{figure}

\paragraph{Cells on vertices of icosahedron}
For cells centered at the vertices of an icosahedron
  inscribed in the unit sphere with $D=1$ and permeabilities
  $d_{1c}=0.8$ and $d_{2c}=0.4$, Fig.~\ref{fig:ident_icosa}
  and the cyan curve in Fig.~\ref{fig:ident_alpha1} predict a
  perturbed HB threshold of $\alpha = \alpha_{cH} +
  \varepsilon\alpha_1 \approx 0.27 + 0.03 \cdot 0.25 = 0.2775$. This
  threshold is consistent with the numerical results in the left
  column of Fig.~\ref{fig:cellicosahedron}. More specifically, for
  $D=1$ and $\alpha = 0.4$ (top), the intracellular dynamics converge
  to a coupled steady state, whereas for $D=1$ and $\alpha = 0.2$
  (middle), they oscillate with increasing in-phase synchrony over
  time. However, increasing $D$ to $D=4$, while keeping $\alpha = 0.2$
  (bottom) restores convergence to the coupled steady state,
  consistent with the asymptotic theory underlying
  Figs.~\ref{fig:ident_icosa} and~\ref{fig:ident_alpha1}.
  Additionally, the three panels in the right column of
  Fig.~\ref{fig:cellicosahedron} show numerical solutions for slightly
  non-identical cells with $D=1$, where the parameters $\alpha_j$ for each $j\in
  \lbrace{1,\ldots,N\rbrace}$ are drawn from a normal distribution
  with standard deviation $0.05$ and a mean near the HB
  threshold as follows: {\bf Top:} $\alpha_j \sim
  \mathcal{N}(\alpha_{cH}+\varepsilon\alpha_1 + 0.1, 0.05^2)$
  (centered slightly outside the oscillatory regime); \textbf{Middle:}
  $\alpha_j \sim \mathcal{N}(\alpha_{cH}+\varepsilon\alpha_1, 0.05^2)$
  (centered directly on the bifurcation threshold); \textbf{Bottom:}
  $\alpha_j \sim \mathcal{N}(\alpha_{cH}+\varepsilon\alpha_1 - 0.1,
  0.05^2)$ (centered inside the oscillatory regime).  As expected, the
  top panel shows convergence to a steady state, whereas the middle
  and bottom panels reveal intricate synchronization dynamics. In
  particular, the collective oscillations in the bottom panel are
  reminiscent of asynchronous drift phenomena studied previously for
  weakly coupled relaxation oscillators of ODE
  systems~\cite{hoppen}. In the pure ODE setting, when relaxation
  oscillators have slightly different natural frequencies and the
  coupling strength falls below the phase-locking threshold, their
  relative phases continuously drift apart. As a result, the spatial
  location where oscillators align and spike together shifts
  continuously over time~\cite{hoppen}. Establishing whether a similar
  mechanism occurs in our integro-differential system is an
  interesting open direction.
  
\begin{figure}[htbp]
\centering
\begin{subfigure}{0.48\textwidth}
    \includegraphics[width=\textwidth]{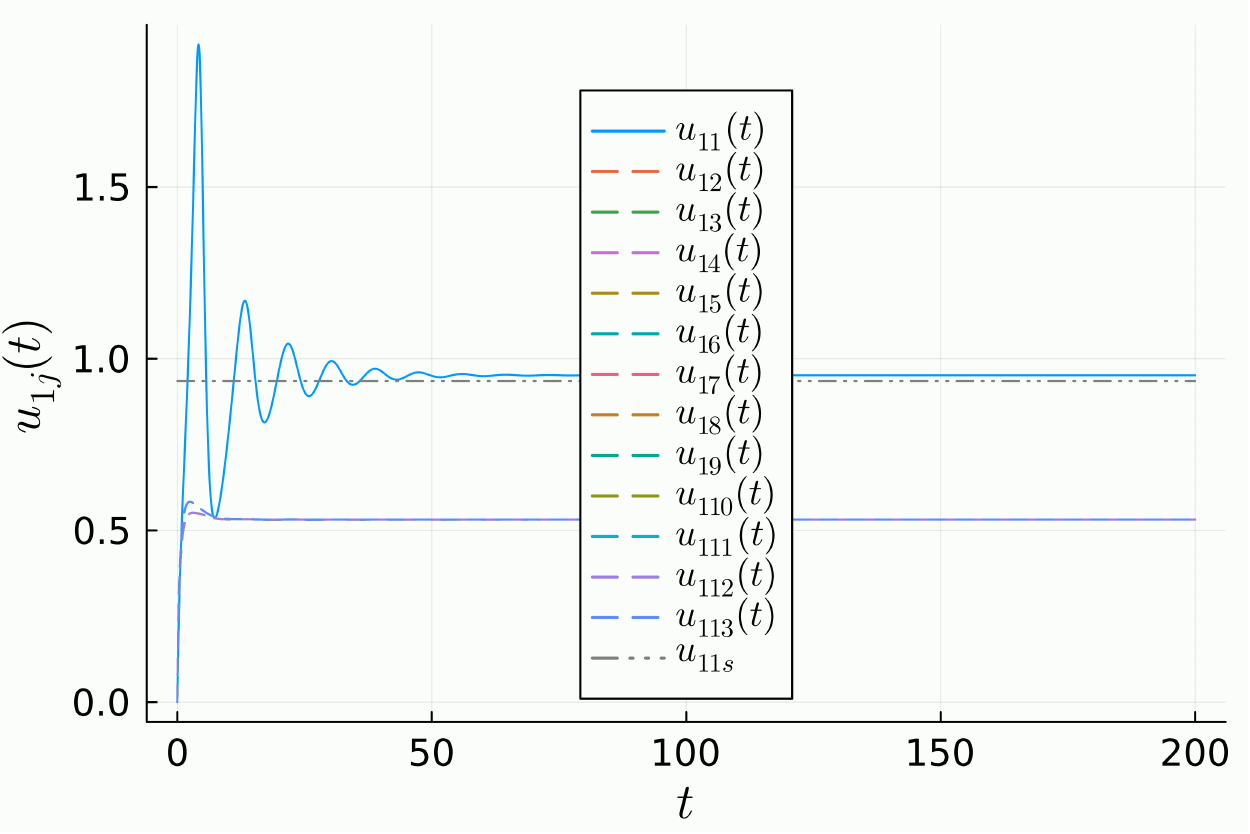}
\end{subfigure}
\hfill
\begin{subfigure}{0.48\textwidth}
    \includegraphics[width=\textwidth]{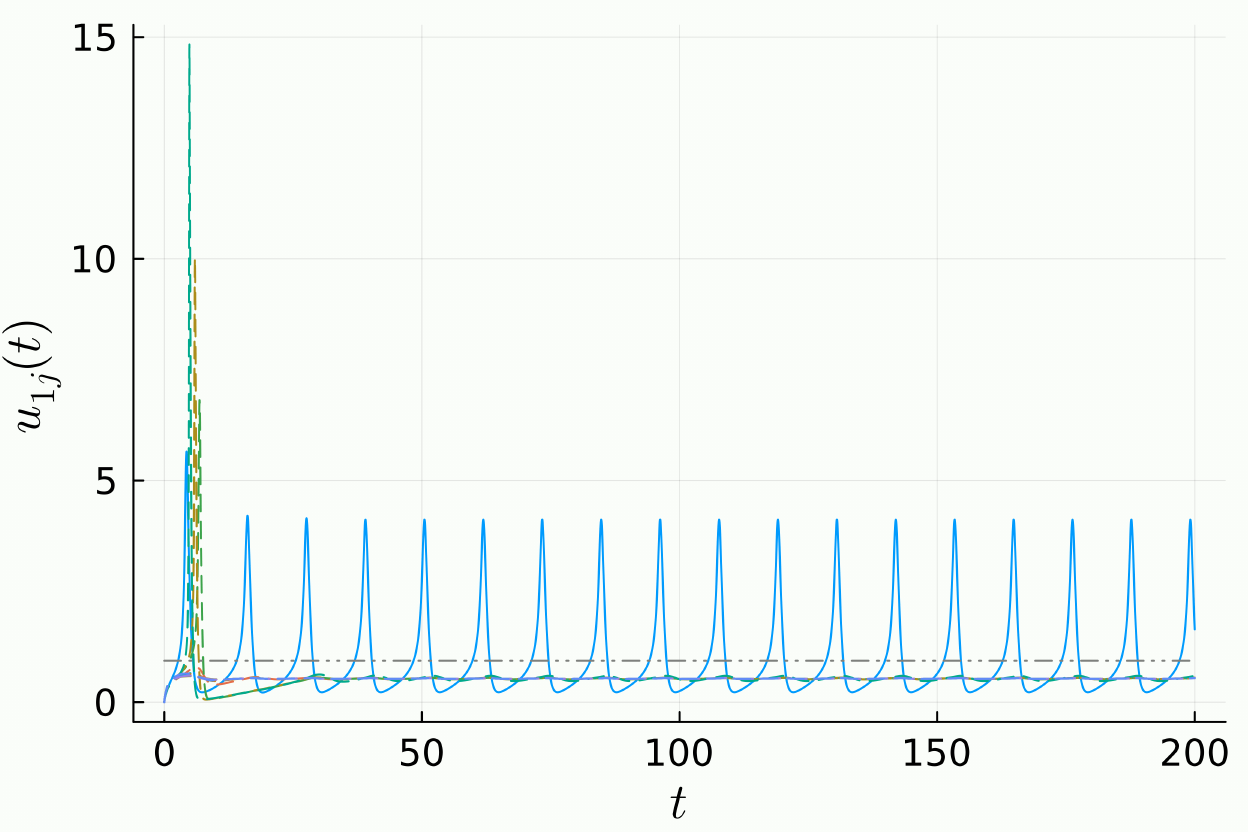}
\end{subfigure}
\begin{subfigure}{0.48\textwidth}
    \includegraphics[width=\textwidth]{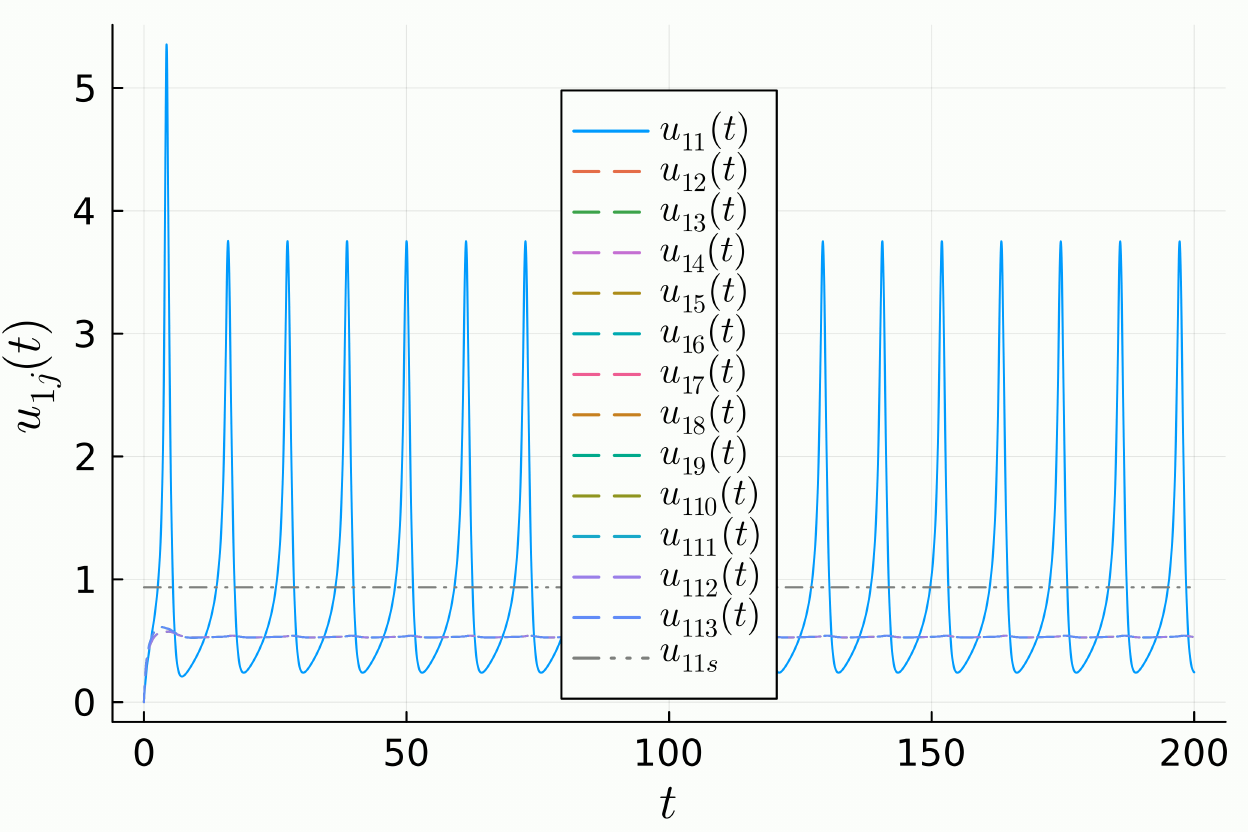}
\end{subfigure}
\hfill
\begin{subfigure}{0.48\textwidth}
    \includegraphics[width=\textwidth]{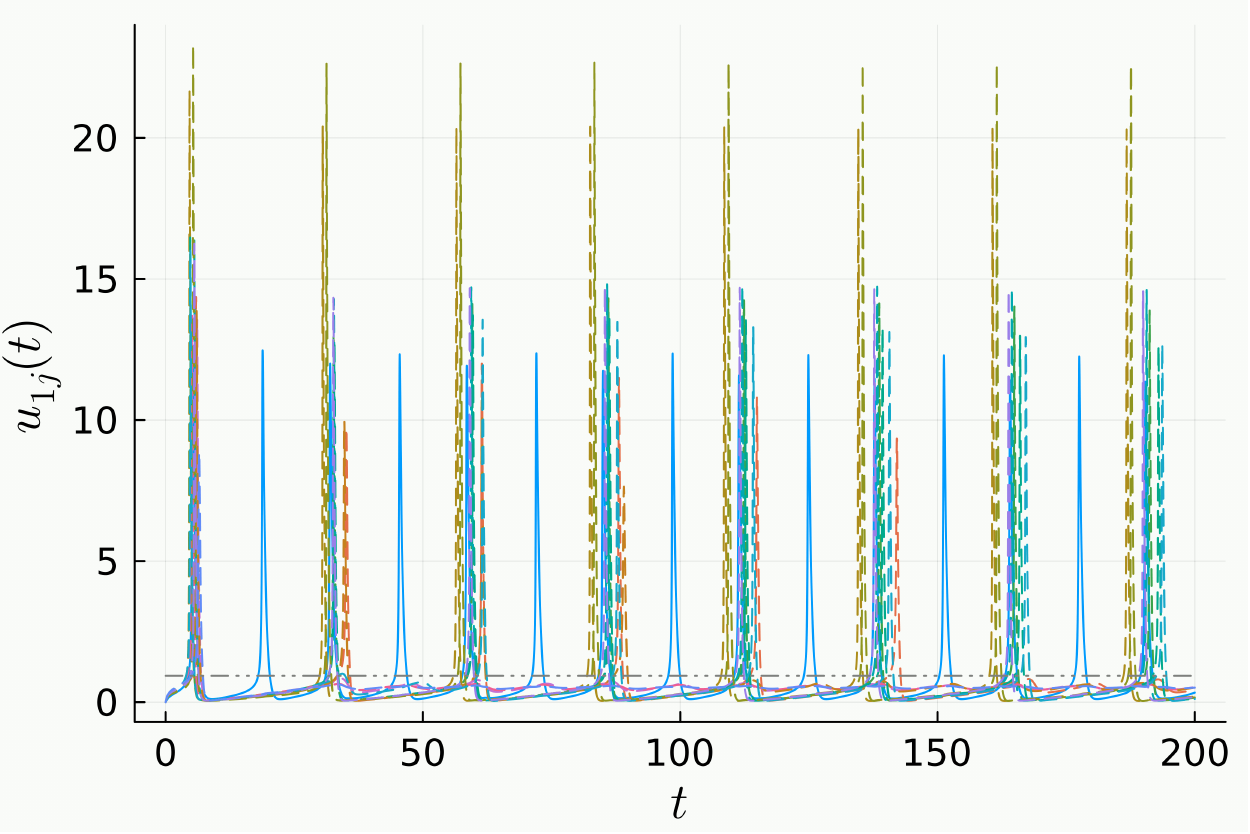}
\end{subfigure}
\begin{subfigure}{0.48\textwidth}
    \includegraphics[width=\textwidth]{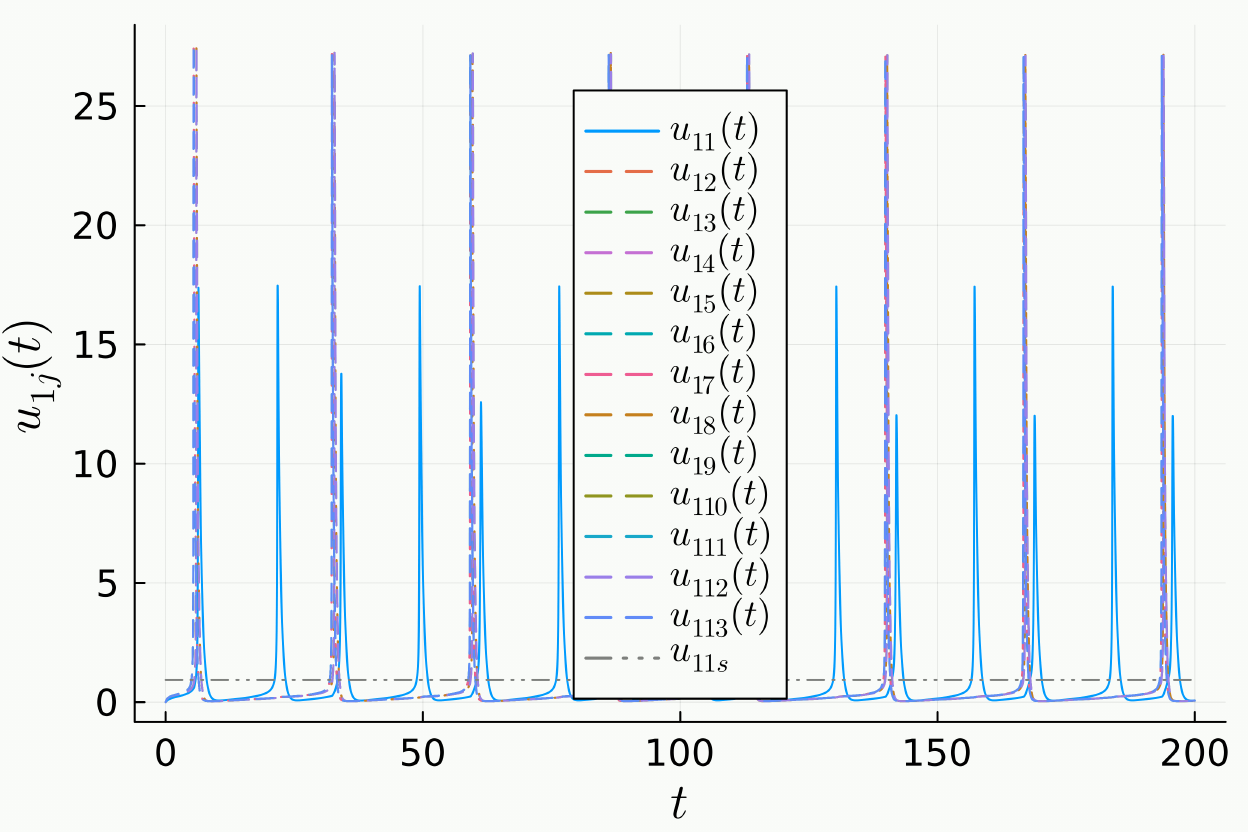}
\end{subfigure}
\hfill
\begin{subfigure}{0.48\textwidth}
    \includegraphics[width=\textwidth]{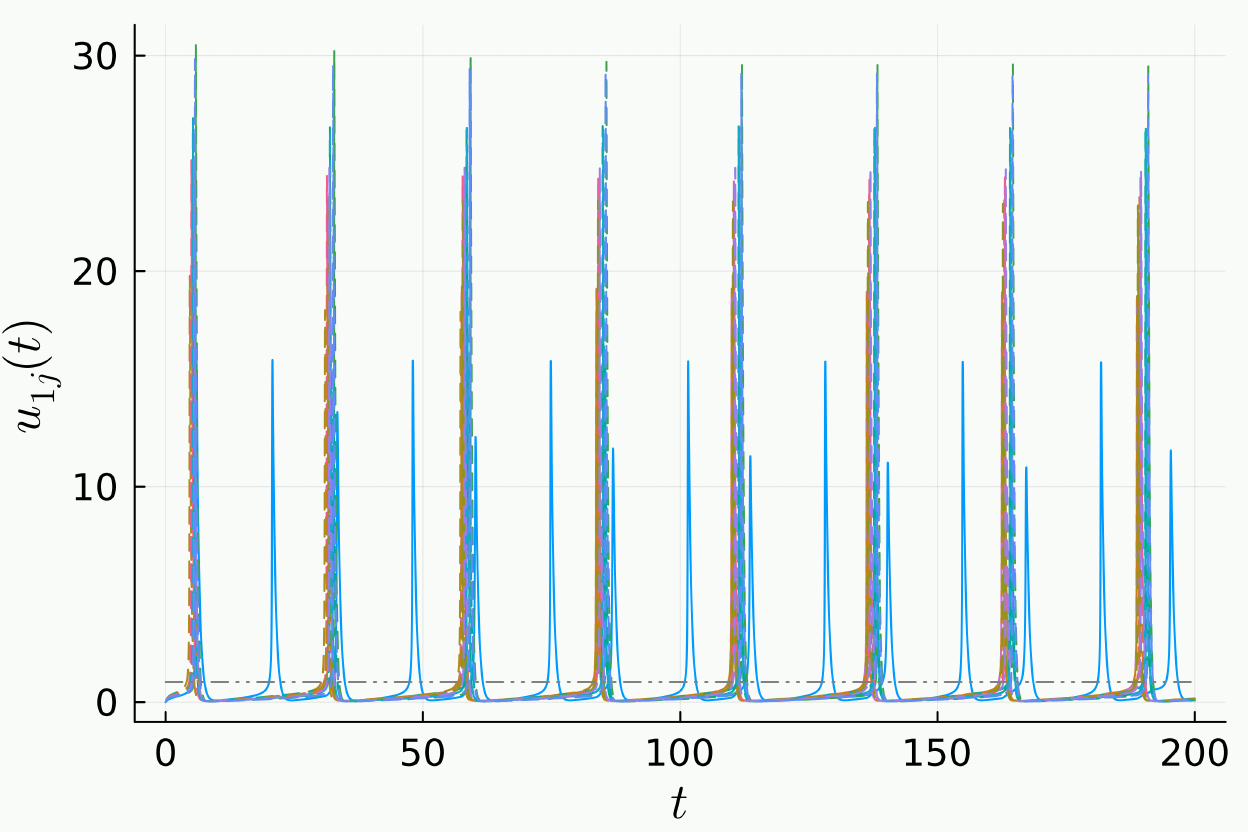}
\end{subfigure}
\caption{Dynamics of cells centered at the vertices of an icosahedron
  that is inscribed in the unit sphere with an additional cell
  centered at the origin. The cells centered at the vertices of the
  icosahedron have permeabilities $d_{1c}=0.8$ and $d_{2c}=0.4$, and
  have dynamics $u_{1j}(t)$ for $j\in \lbrace{2,\ldots,13\rbrace}$,
  while the center cell with dynamics $u_{11}(t)$ has the different
  permeabilities $d_{1p}=1.2$ and $d_{2p}=0.2$. For these parameter
  values and with $D=1$, the center cell sets the HB threshold
  (cf. Fig.~\ref{fig:caseI_II:HB} (right) with small correction of
  Fig.~\ref{fig:caseII_alpha}), and so acts as a pacemaker. Our
  asymptotic theory predicts that for $\alpha<\alpha_{pH} +
  \varepsilon\alpha_{p1} \approx 0.68 - 0.03\cdot 0.33 = 0.670$ the
  center cell oscillates and that for even lower $\alpha$-values with
  $\alpha<\alpha_{cH}+\varepsilon\alpha_{c1}\approx0.27+0.03\cdot 0.25
  = 0.2775$ the cells centered on the icosahedron vertices oscillate
  as well. Left column: When $\alpha = 0.8$ and $D = 1$ (top), the
  cells converge to their coupled steady-state. For $\alpha = 0.5$ and
  $D = 1$, only the center cell has large-amplitude
  oscillations. Decreasing to $\alpha=0.15$ while fixing $D=1$, all
  cells oscillate strongly while synchronizing with a
  phase-shift. Right column: Fixing $D=1$, cells are taken to be
  slightly non-identical by choosing $\alpha = \alpha_j$ for each
  $j\in\{1,...,N\}$ to be normally distributed according to $\alpha_j
  \sim \mathcal{N}(\alpha_{cH}+\varepsilon\alpha_{c1} + 0.1, 0.05^2)$
  (top; mean slightly above the oscillatory threshold of the
  icosahedron vertex cells), then $\alpha_j \sim
  \mathcal{N}(\alpha_{cH}+\varepsilon\alpha_{c1}, 0.05^2)$ (center;
  mean on the bifurcation threshold), and lastly $\alpha_j \sim
  \mathcal{N}(\alpha_{cH}+\varepsilon\alpha_{c1} - 0.1, 0.05^2)$
  (bottom; mean inside the oscillatory parameter region). All other
  kinetic and numerical parameters are the same as for
  Fig.~\ref{fig:caseI_II:HB} and Fig.~\ref{fig:twononidcells},
  respectively. The initial condition is as in the caption of
  Fig.~\ref{fig:icosahedron}.  The numerical time-stepping and
  plotting required $\sim 7 \text{ min}$ per sub-figure (hardware
  specifications as in the caption of Fig.~\ref{fig:twononidcells}).}
  \label{fig:cellicosahedron_center}.
\end{figure}

\paragraph{Cells on vertices of icosahedron with center cell as pacemaker}
We now place an additional cell at the center of the icosahedron
arrangement that has the modified permeabilities $d_{1p}=1.2$ and
$d_{2p}=0.2$, while the cells centered on the vertices of the
icosahedron have permeabilities $d_{1c} = 0.8$ and $d_{2c} = 0.4$. The
center cell serves as the pacemaker cell since the HB
threshold is set by it for this choice of permeabilities (see the plot
in Fig.~\ref{fig:caseI_II:HB} (right) with small correction from
Fig.~\ref{fig:caseII_alpha}). For $D=1$, the pacemaker cell oscillates
for $\alpha<\alpha_{pH} + \varepsilon\alpha_{p1} \approx
0.68-0.03\cdot 0.32 = 0.6704$ and the cells on the icosahedron
vertices should oscillate for
$\alpha<\alpha_{cH}+\varepsilon\alpha_{c1}\approx0.27+0.03\cdot 0.25 =
0.2775$, according to the asymptotic
theory. Fig.~\ref{fig:cellicosahedron_center} shows that the
corresponding numerical solutions agree well with these theoretical
predictions.

\begin{figure}[htbp]
\centering
\begin{subfigure}{0.48\textwidth}
    \includegraphics[width=\textwidth]{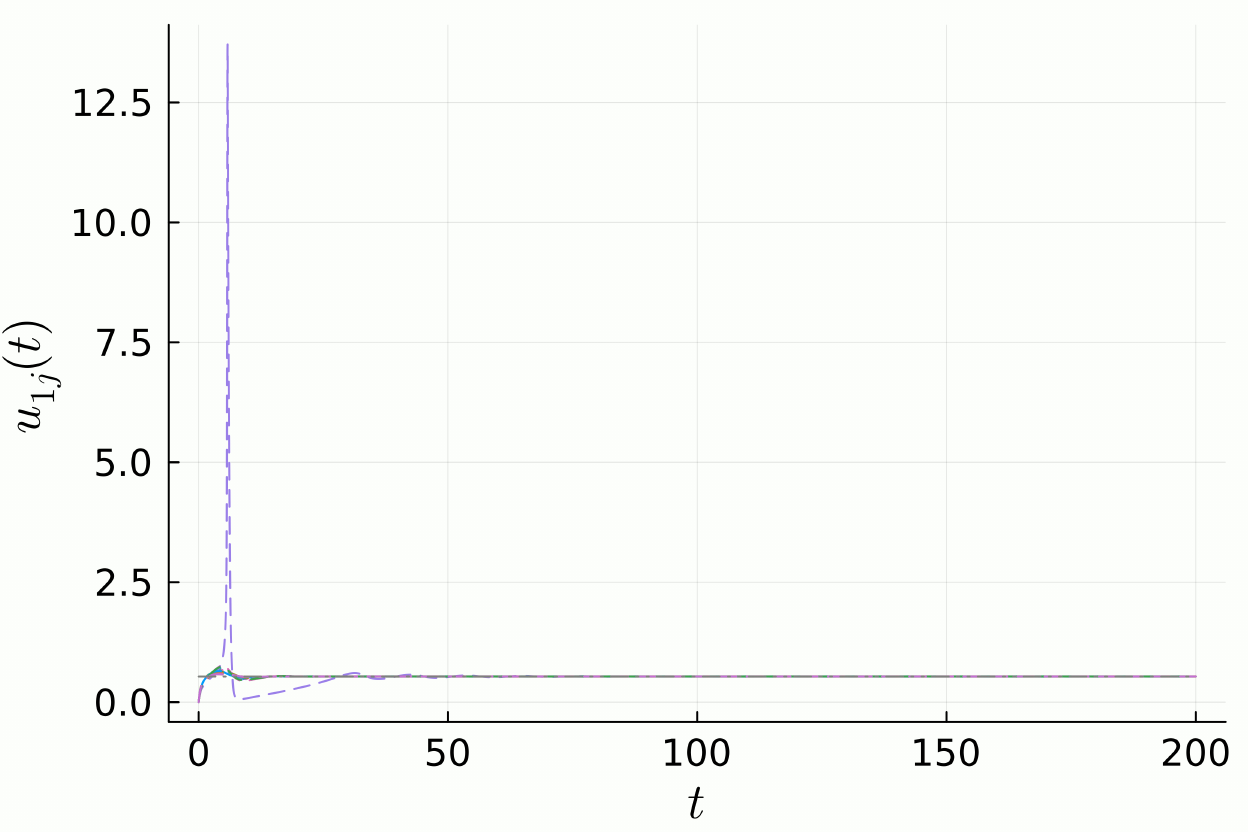}
\end{subfigure}
\hfill
\begin{subfigure}{0.48\textwidth}
    \includegraphics[width=\textwidth]{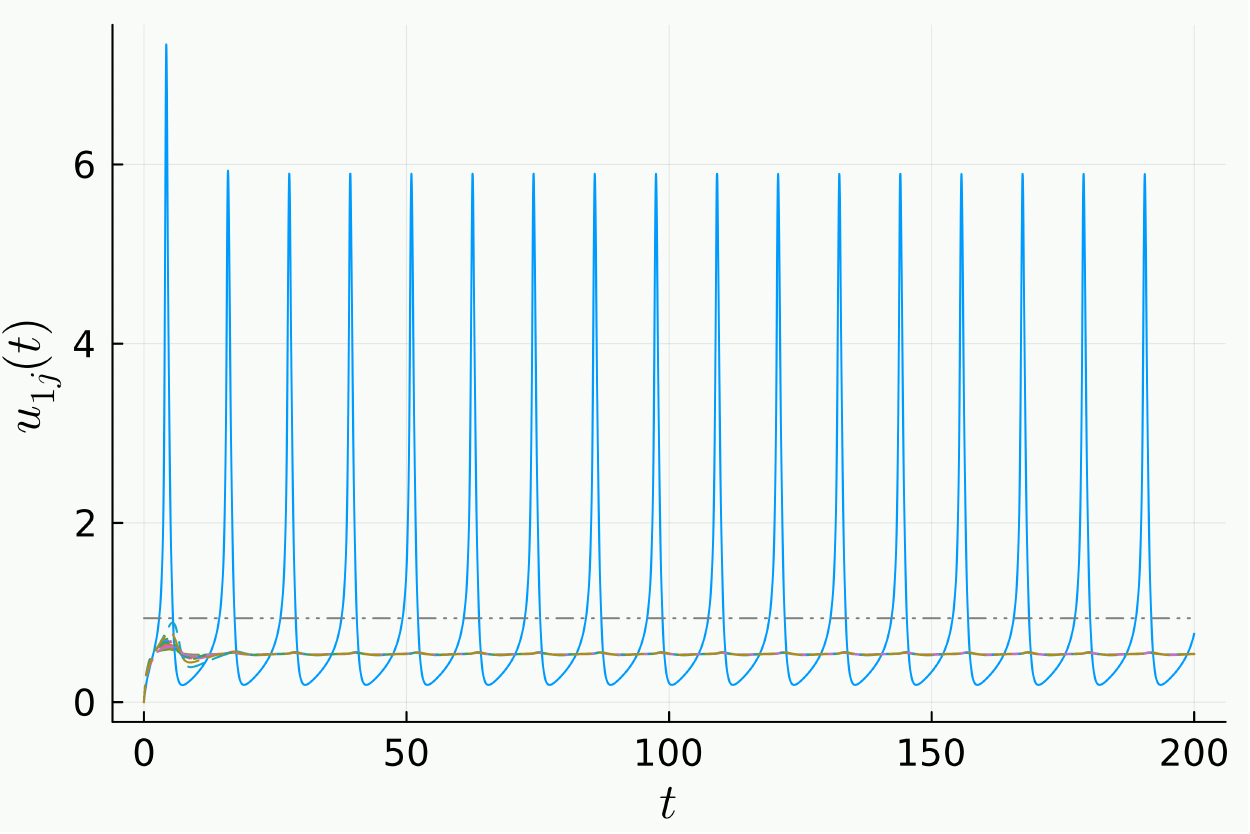}
\end{subfigure}
\begin{subfigure}{0.48\textwidth}
    \includegraphics[width=\textwidth]{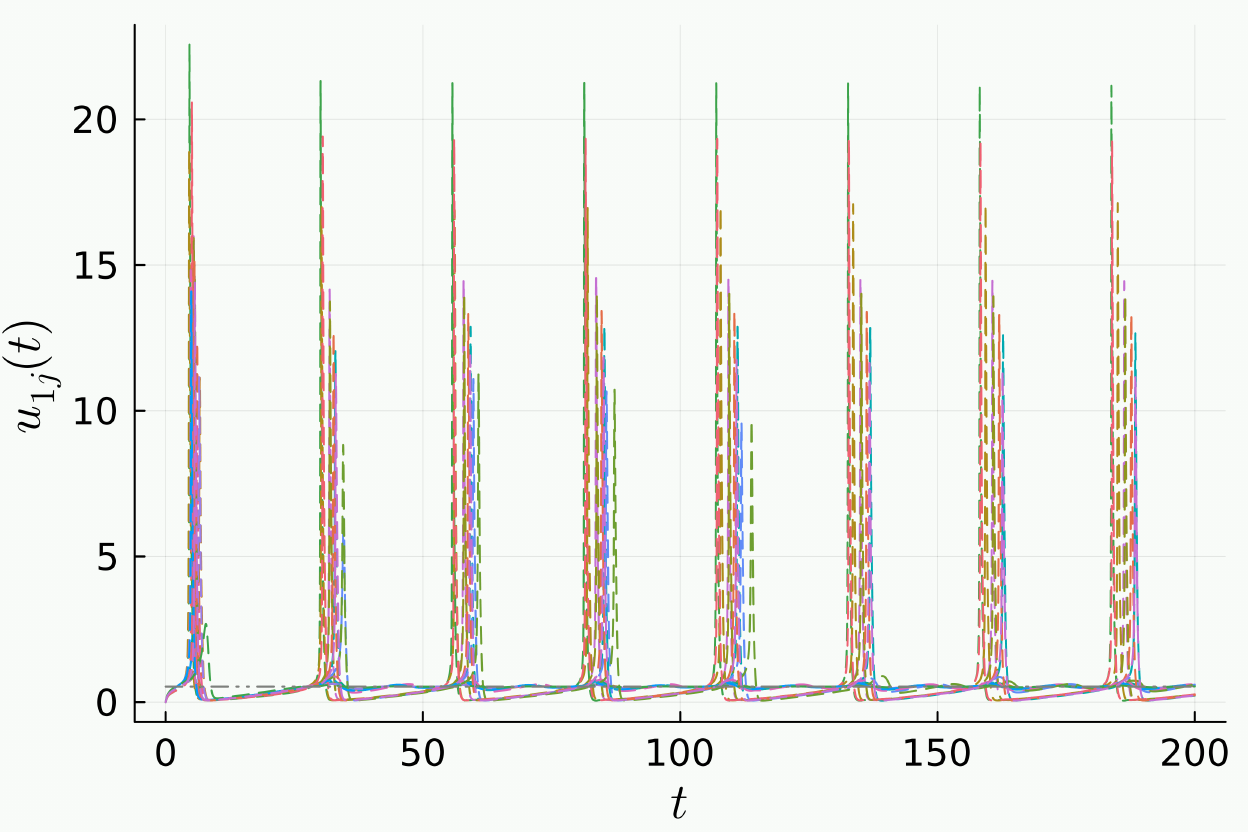}
\end{subfigure}
\hfill
\begin{subfigure}{0.48\textwidth}
    \includegraphics[width=\textwidth]{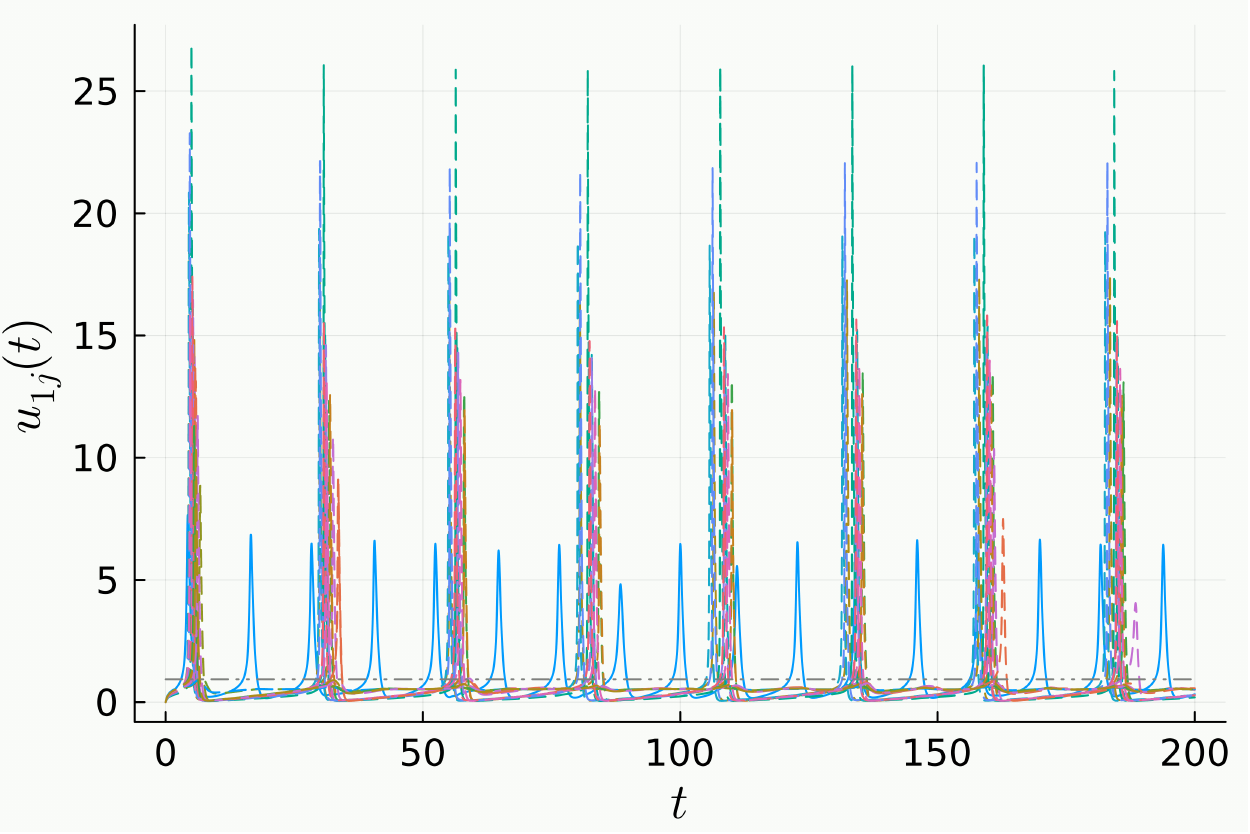}
\end{subfigure}
\begin{subfigure}{0.48\textwidth}
    \includegraphics[width=\textwidth]{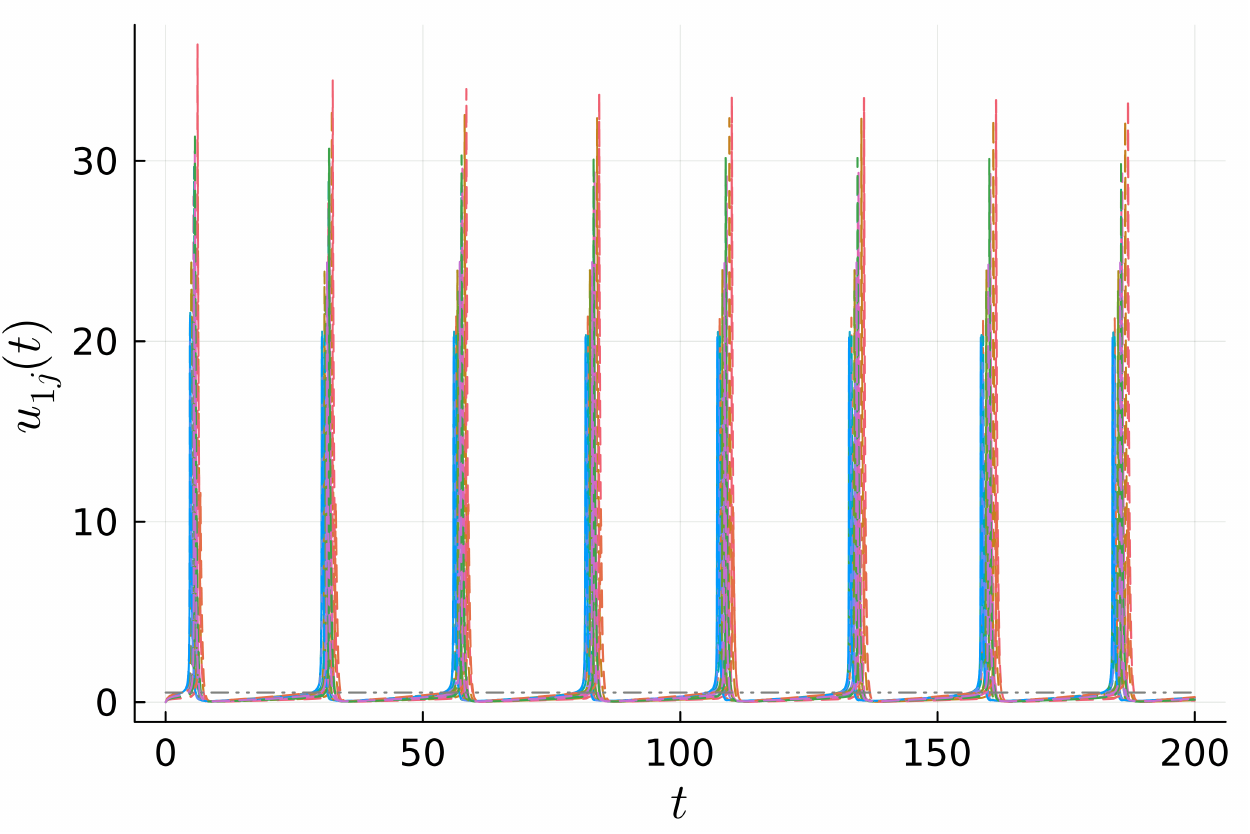}
\end{subfigure}
\hfill
\begin{subfigure}{0.48\textwidth}
    \includegraphics[width=\textwidth]{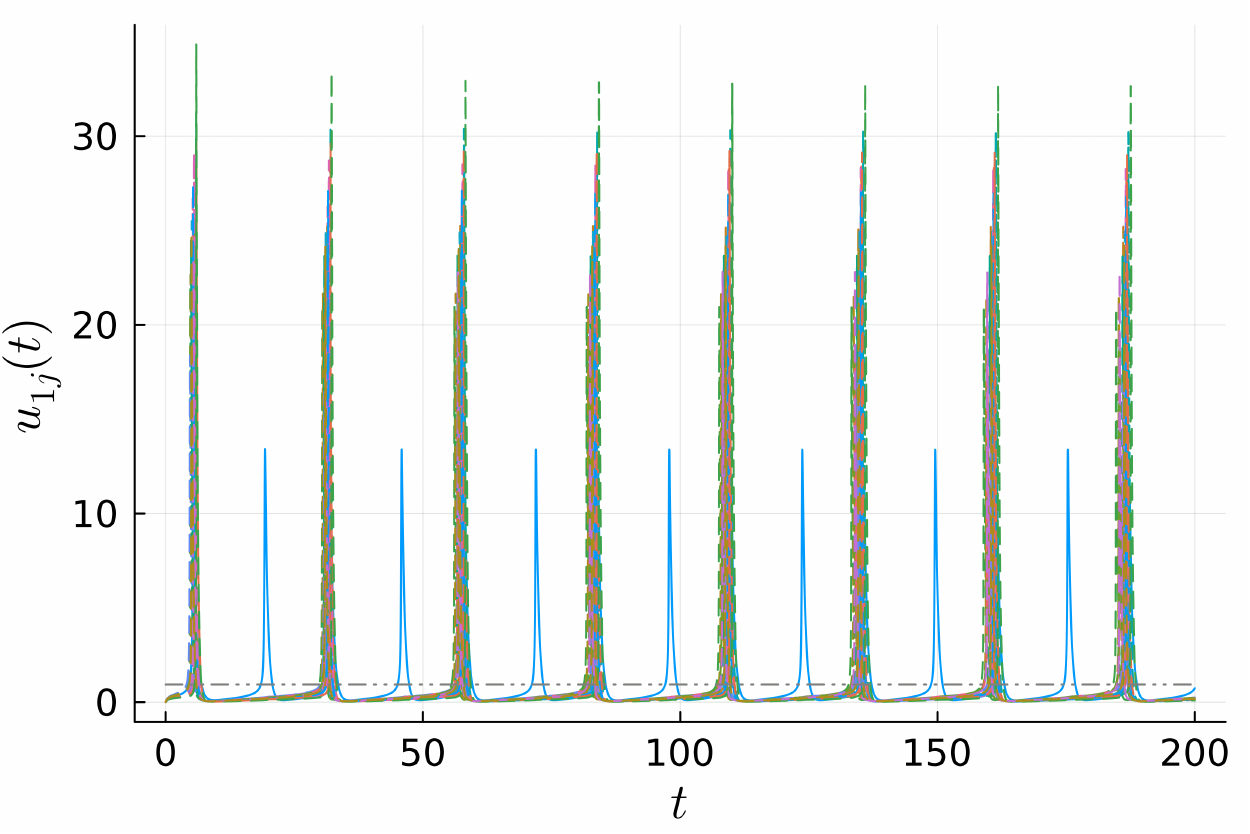}
\end{subfigure}
\caption{Dynamics of cells centered at the vertices of a
    dodecahedron that is inscribed in the unit sphere, having
    identical cell permeabilities $d_{1c}=0.8$ and $d_{2c}=0.4$, both
    without a center cell (left column) and with an additional cell
    centered at the origin (right column; center cell dynamics in
    cyan). By choosing the permeabilities of the center cell as
    $d_{1p}=1.2$ and $d_{2p}=0.2$, and fixing $D=1$, the center cell
    sets the HB threshold since from \eqref{num_ident:HB_zero} and
    \eqref{ringc:alpha_pH} we calculate $\alpha_{cH}=0.263$ and
    $\alpha_{pH}=0.592$, respectively, so that
    $\alpha_{pH}>\alpha_{cH}$.  With the pacemaker cell, our
    asymptotic theory predicts that for $\alpha<\alpha_{pH} +
    \varepsilon\alpha_{p1} \approx 0.592 - 0.03\cdot 0.33 = 0.581$ the
    center cell oscillates and when
    $\alpha<\alpha_{cH}+\varepsilon\alpha_{c1}\approx 0.263+0.03\cdot
    0.7 = 0.284$ the cells at the vertices of the dodecahedron should
    also begin to oscillate. In both columns $\alpha_j$ for each
    $j\in\lbrace{1,...,N\rbrace}$ was normalized distributed as: {\bf
      Top Row:} $\alpha_j \sim
    \mathcal{N}(\alpha_{cH}+\varepsilon\alpha_{c1} + 0.1, 0.05^2)$;
            {\bf Middle Row:} $\alpha_j \sim
            \mathcal{N}(\alpha_{cH}+\varepsilon\alpha_{c1}, 0.05^2)$;
                    {\bf Bottom Row:} $\alpha_j \sim
                    \mathcal{N}(\alpha_{cH}+\varepsilon\alpha_{c1} -
                    0.1, 0.05^2)$. The initial condition is as in the
                    caption of Fig.~\ref{fig:icosahedron}.  The
                    numerical time-stepping and plotting required
                    $\sim 15 \text{ min}$ per sub-figure for the
                    dodecahedron arrangement without center cell and
                    $\sim 17 \text{ min}$ with the pacemaker center
                    cell (hardware specifications as in the caption of
                    Fig.~\ref{fig:twononidcells}).}
\label{fig:dodecahedron_no_c_with_c}
\end{figure}

\paragraph{Cells on vertices of dodecahedron with/without center cell} 
Finally, we consider a cell configuration with cells centered
  at the vertices of a dodecahedron inscribed in the unit sphere
  having common permeabilities $d_{1c}=0.8$ and $d_{2c}=0.4$. From
  \eqref{num_ident:HB_zero}, together with the magenta curve in
  Fig.~\ref{fig:ident_alpha1}, the HB threshold when $D=1$ is
  approximately $\alpha = \alpha_{cH}+\varepsilon\alpha_1 \approx
  0.263 + 0.03\cdot 0.7 = 0.284$ in the absence of a center cell.  By
  including a pacemaker cell centered at the origin with
  permeabilities $d_{1p}=1.2$ and $d_{2p}=0.2$, we calculate from
  \eqref{ringc:alpha_pH}, together with the magenta curve in
  Fig.~\ref{fig:caseII_alpha}, that the HB threshold for the pacemaker
  cell is $\alpha = \alpha_{pH}+\varepsilon\alpha_{p1} \approx 0.592 -
  0.03\cdot 0.33 = 0.581$. As a result, the pacemaker cell controls
  the initiation of intracellular oscillations. Fixing $D=1$,
  intracellular dynamics computed from the time-marching scheme are
  shown in the right and left columns of
  Fig.~\ref{fig:dodecahedron_no_c_with_c} for the dodecahedron
  arrangement with and without the pacemaker cell, respectively. For
  these plots we again considered slightly non-identical cells where
  the reaction kinetic parameter $\alpha = \alpha_j$ for each
  $j\in\lbrace{1,...,N\rbrace}$ was normally distributed as follows:
  {\bf Top Row:} $\alpha_j \sim
  \mathcal{N}(\alpha_{cH}+\varepsilon\alpha_1 + 0.1, 0.05^2)$; {\bf
    Middle Row:} $\alpha_j\sim \mathcal{N}(\alpha_{cH}+
  \varepsilon\alpha_1, 0.05^2)$; {\bf Bottom Row:} $\alpha_j \sim
  \mathcal{N}(\alpha_{cH}+\varepsilon\alpha_1 - 0.1, 0.05^2)$.
  Without the center cell, in the top left panel of
  Fig.~\ref{fig:dodecahedron_no_c_with_c} we observe, as expected, an
  approach to the steady-state, while synchronized relaxation
  oscillations occur when $\alpha_j$ is within the oscillatory
  parameter regime (bottom left panel of
  Fig.~\ref{fig:dodecahedron_no_c_with_c}). In contrast, with the
  center cell, we observe from the top right panel that, as expected,
  only the center cell oscillates. However, when $\alpha_j$ is well
  below $\alpha_{cH}$, the cells at the vertices of the dodecahedron
  exhibit synchronized relaxation oscillations (bottom right panel).

\section{Discussion}\label{sec:discussion}

We have developed and implemented a hybrid asymptotic-numerical theory
based on the cell-bulk model \eqref{DimLess_bulk} in $\R^3$ to analyze
how collective behavior, such as synchronous oscillations, can emerge
through the long-range diffusive coupling of spatially localized
dynamically active units, referred to as ``cells''. For this
  ODE-PDE model, in the limit of small common cell radius we have
  derived a new integro-differential ODE system for the intracellular
  dynamics of a collection of cells. In this limiting system, the bulk
  diffusion field is effectively replaced with an additive time-memory
  kernel term whose strength depends on the ratio of the influx to
  efflux permeability, the bulk diffusivity $D$ and the inter-cell
  distances. An efficient time-marching scheme inspired by the {\em
  sum-of-exponentials} method of \cite{greengard_2015} and
\cite{beylkin_2005}, was developed to numerically compute solutions to
the nonlocal integro-differential system over long time intervals.  By
applying this hybrid framework to the special case of Sel'kov
intracellular kinetics, we have illustrated synchronization behavior
of intracellular oscillations for various spatial configurations of
cells. The methodology developed herein in $\R^3$ is an extension of
that developed in \cite{pelzward2D} for the 2-D case.

We now highlight the essential role of the spatial dimension in
initiating collective behavior. In the 2-D case, owing to the
$\log{r}$ behavior of the 2-D free-space Green's function, where
$r$ is the radial distance from the cell center, the signaling
molecules released from each cell do not dilute rapidly in the
bulk medium. As a result, the dynamically active units are strongly
coupled to the total cell population, rendering the 2-D cell-bulk
model highly prone to global synchronization behavior. In
\cite{pelzward2D} it was shown that the strength of cell-cell coupling
is ${\mathcal O}\left({-1/\log\eps}\right)$, which is effectively
${\mathcal O}(1)$ unless $\eps$ is extremely small. As a result, it
was observed in \cite{pelzward2D} that changing a kinetic parameter in
a single cell can trigger intracellular oscillations in the entire
cell population through a Hopf bifurcation threshold that is set by
the cell population. These intracellular oscillations become
increasingly synchronized as the bulk diffusivity increases.

In contrast, in $\R^3$ the signaling molecules released by a cell
degrade rapidly into the bulk medium owing to the sharp ${\mathcal
  O}\left({1/r} \right)$ decay of the 3-D free-space Green's
function. As a result, each single cell will undergo an oscillatory
instability only if its reaction kinetic parameters are chosen so as
to overcome the strong ${\mathcal O}\left({1/r}\right)$ diffusion loss
into the bulk medium. The inclusion of additional dynamically active
units at ${\mathcal O}(1)$ distances will only have a weak ${\mathcal
  O}(\eps)$ effect on the Hopf bifurcation thresholds associated with
an individual cell. In this way, Hopf bifurcation thresholds
triggering collective dynamics for a population of cells are
${\mathcal O}(\eps)$ close to those for individual cells. This implies
that reaction kinetic parameters need only exhibit small ${\mathcal
  O}(\eps)$ variations near the boundary between stability and
instability for each individual cell in order to initiate collective
behavior in the entire cell population that will occur over ${\mathcal
  O}(\eps^{-1})$ time scales. This qualitative feature is consistent
with the viewpoint in \cite{mora} that many biological systems operate
at their marginal stability threshold, since it is energetically
unfavorable for biological cells to strongly modify their internal
dynamics. However, from a mathematical viewpoint, in order to study
the effect of small ${\mathcal O}(\eps)$ perturbations due to distant
cells, we had to develop a novel matrix perturbation theory on the
linear stability problem for steady-state solutions in order to
calculate how the marginal stability thresholds of each cell are
modified by these small perturbations. This perturbation theory was
implemented for Sel'kov intracellular dynamics, and some of the
analytical results predicted by this analysis were confirmed by full
time-dependent simulations of the integro-differential system.

There are two key open directions that warrant further study. Firstly,
it would be interesting to analyze the role of the bulk-diffusivity
$D$ in a 3-D bounded domain with a reflecting boundary.  For an
asymptotically large $D$, it was shown in \cite{sarafa_3D} that the
cell-bulk model \eqref{DimLess_bulk} in a bounded 3-D domain can be
reduced to a limiting ODE system with global coupling. In this
well-mixed $D\gg 1$ limit, the global coupling arises from the bulk
diffusion field that, to leading order, depends only on time. As a
result, the thresholds for initiating intracellular oscillations, and
synchronization behavior, of spatially localized cells is determined
not by individual cells but instead by the entire cell population, and
these thresholds are only weakly influenced by the spatial
configuration of cells. For the unexplored regime where $D={\mathcal
  O}(1)$ in a bounded 3-D domain, we expect that intracellular
oscillations triggered by cell-cell coupling is easier to achieve than
in all of $\R^3$, but that the spatial configuration of cells will
also play a key role in any collective behavior. For either a
spherical or a finite cylindrical confining domain, the Laplace
transform of the time-dependent Green's function and its regular part
can be represented analytically as infinite series and they have
singularities along the negative real axis in the Laplace transform
space. As a result, the sum-of-exponentials numerical approach is
still applicable to this setting for computing solutions to a modified
integro-differential system that arises from the asymptotic
analysis. By allowing for a finite domain, and by replacing the
conceptual Sel'kov intracellular dynamics with dynamics relevant to
specific applications, such as Belousov-Zhabotinsky reactants
encapsulated inside microcomparments \cite{torbensen}, numerical
results from the limiting integro-differential ODE system should have
some predictive power in modeling various experimental settings
involving collective dynamics of spatially segregated units.

As a second extension of the theoretical framework, it would be
interesting to couple the dynamics of the cell locations to the
spatial gradient of intercellular signals in $\R^3$. In many
quorum-sensing applications (cf.~\cite{gregor2010}, \cite{dicty}),
cell motion due to spatial gradients of the bulk signal are believed
to be the mechanism by which cells can form spatial clusters to
strengthen short-range paracrine signaling. From a mathematical
viewpoint, clusters of cells within an ${\mathcal O}(d)$ length-scale
with $\eps\ll d \ll 1$ will interact strongly even in $\R^3$ as the
${\mathcal O}(1/r)$ spatial decay of the Green's function is no longer
relevant. More specifically, if a few dynamically active units are
clustered tightly together in $\R^3$, they can achieve a localized
quorum and synchronize independently of the rest of the cell
population, thereby creating localized pockets of activity. This
collective behavior arises from the spatial proximity of cells,
without the need for any modification of the intracellular kinetic
parameters. Distinct cell clusters that are separated by ${\mathcal
  O}(1)$ distances will interact very weakly owing to the ${\mathcal
  O}(1/r)$ decay of the Green's function.  As a result we expect that
chimera states should be readily observable in $\R^3$ whenever there
is a spatial clustering of the cell population into distinct groups of
cells. Also when spatially grouped together, the collective
  oscillation of each closely spaced group could induce in-phase or
  anti-phase oscillations of an ${\mathcal O}(1)$ distant group. A
  theoretical characterization of this anticipated phenomenon obtained
  from extending the cell-bulk model \eqref{DimLess_bulk} to allow for
  cell motion would be worthwhile.

\begin{appendices}
\renewcommand{\theequation}{\Alph{section}.\arabic{equation}}
\setcounter{equation}{0}

\section{Non-Dimensionalization of the Cell-Bulk Model}
\label{app:nondim}

In this appendix we non-dimensionalize \eqref{Dim_bulk} to obtain
\eqref{DimLess_bulk}. The dimensions, labeled by $\left[\ldots\right]$,
of the various terms in \eqref{Dim_bulk} are
\begin{equation}\label{app:dimensions}
  \begin{aligned}
   \left[{\mathcal U}\right] &= \frac{\mathrm{moles}}{\mathrm{(length)}^3}
   \,, \quad \left[\v{\mu}_j\right] = \mathrm{moles} \,, \quad
   \left[\mu_c \right] = \mathrm{moles} \,, \quad \left[D_B \right]
   =\frac{\mathrm{(length)}^2}{\mathrm{time}} \,, \quad\\ \left[k_B \right]
   &=\left[k_R\right] = \frac{1}{\mathrm{time}} \,, \quad \left[\beta_1
     \right] =\frac{\mathrm{length}}{\mathrm{time}} \,,
   \quad\ \left[\beta_2 \right]
   =\frac{1}{\mathrm{length}^{2}\times\mathrm{time}} \,.
\end{aligned}
\end{equation}
We introduce the dimensionless variables $t$, $\v{x}$, $U$ and $\v{u}$ defined
as
\begin{equation}\label{form:non_dim}
   t \equiv k_R T \,, \quad  \v{x} \equiv \frac{\v{X}}{L} \,, \quad
 U\equiv \frac{L^3}{\mu_c} {\mathcal U} \,, \quad \v{u}_j \equiv 
 \frac{\v{\mu}_j}{\mu_c} \,, \quad \eps \equiv \frac{R_0}{L}\ll 1 \,,
\end{equation}
where $L=\mathrm{min}|\v{x}_i-\v{x}_j|$ for $i\neq j$ and $R_0$ is the
common cell radius. In this way, the dimensionless bulk field
satisfies \eqref{DimLess_bulka}, while \eqref{Dim_bulkb} on the
$j^{\mathrm{th}}$ cell boundary becomes
\begin{equation}\label{app:Dim_bulkb}
    \eps D \, \partial_{n_{\x}} \, U =
    \frac{\eps \beta_{1j}}{L k_R} \, U - \frac{\eps\beta_{2j} L^2}{k_R} \,
    u_{1j}\,, \quad \x \in \partial \Omega_{\eps_j}\,.
\end{equation}
By noting that $\text{d}S_{\v{X}}=L^2\text{d}S_{\v{x}}$, the intracellular
kinetics \eqref{Dim_Intra} becomes
\begin{equation}\label{app:Dim_Intra}
  \frac{d\v{u}_j}{dt} = \v{F}_j(\v{u}_j) + \frac{\v{e}_1}{\eps}
  \int_{\partial \Omega_{\eps_j}} \left( \frac{\eps \beta_{1j}}{L k_R} -
    \frac{\eps L^2 \beta_{2j}}{k_R} \, u_{1j}\right) \, \text{d}S_{\v{x}} \,.
\end{equation}
Since
$|\partial\Omega_{\varepsilon_j}|= {\mathcal O}(\varepsilon^2)$, in
order to ensure that there is an ${\mathcal O}(1)$ efflux out of the
$j^{\mathrm{th}}$ cell into the bulk medium, we must consider the limit
where ${\beta_{2j} L^2/k_R}={\mathcal O}(\varepsilon^{-2})$. Moreover,
the feedback into the $j^{\mathrm{th}}$ cell from the bulk is
${\mathcal O}(\varepsilon)$ when the dimensionless influx parameter is
on the asymptotic range
${\beta_{1j}/(L k_R)}={\mathcal O}(\varepsilon^{-1})$. Based on these
observations, we chose the permeability scalings as in
\eqref{dim:param}. With this choice, \eqref{app:Dim_bulkb} and
\eqref{app:Dim_Intra} become \eqref{DimLess_bulkb} and
\eqref{DimLess_Intra}, respectively.

\section{Details of Transient Analysis}\label{app:laplace}

In this appendix we provide the details on the derivation of
the matching condition \eqref{3:v_final}.

To estimate \eqref{3:lap_trans} near its singularity with the largest
real part occuring at $s=0$, we let $s\to 0$ in \eqref{3:lap_trans}
and use a Taylor series approximation to obtain
\begin{equation}\label{3:small_s}
  \begin{split}    
  \hat{V}(\rho,s) &\sim \frac{d_{2j} u_{1j}(0)}{\eps \rho (d_{1j}+D)}
  \left[ \frac{1}{s} + \frac{\sqrt{s+\eps^2\sigma}}{s\sqrt{D}}
    \left( \frac{d_{1j}}{d_{1j}+D} -\rho \right) + {\mathcal O}(1) \right]\\
  & \qquad + \frac{\eps d_{2j} u_{1j}^{\prime}(0)}{\rho (d_{1j}+D)}
  \left[ \frac{1}{s^2} + \frac{\sqrt{s+\eps^2\sigma}}{s^{2} \sqrt{D}}
    \left( \frac{d_{1j}}{d_{1j}+D} -\rho \right) + {\mathcal O}(1) \right]
  \,, \quad \mbox{as} \quad s\to 0\,.
  \end{split}
\end{equation}
To invert \eqref{3:small_s}, we use the inverse Laplace transforms
${\mathcal L}^{-1}\left[{1/s^{p+1}}\right]={\tau^p/\Gamma(p+1)}$,
where $\Gamma(x)$ is the Gamma function, together with
\begin{equation}\label{3:inv_lap}
  \begin{split}
    {\mathcal L}^{-1}\left[\frac{\sqrt{s+\eps^2\sigma}}{s}\right]
    &= g_{0}(\eps^2\sigma \tau)\sqrt{\eps^2\sigma} \,, \quad
    \mbox{with} \quad g_0(\zeta)\equiv \frac{e^{-\zeta}}{\sqrt{\pi\zeta}}+
    \mathrm{erf}(\sqrt{\zeta}) \,, \\
    {\mathcal L}^{-1}\left[\frac{\sqrt{s+\eps^2\sigma}}{s^2}\right]
    &= \frac{g_{1}(\eps^2\sigma \tau)}{\sqrt{\eps^2\sigma}} \,, \quad
    \mbox{with} \quad g_1(\zeta)\equiv \left(\zeta+\frac{1}{2}\right)
    \mathrm{erf}(\sqrt{\zeta}) + \sqrt{\frac{\zeta}{\pi}} e^{-\zeta} \,,
  \end{split}
\end{equation}
where $\mathrm{erf}(z)$ is the Error function.  Upon inverting
\eqref{3:small_s}, we conclude for $\rho>1$ fixed and $\tau\to \infty$
that
\begin{equation}\label{3:v_asy_1}
  \begin{split}    
   V(\rho,\tau) &\sim \frac{d_{2j} u_{1j}(0)}{\eps \rho (d_{1j}+D)}
  \left[ 1 + \left( \frac{d_{1j}}{d_{1j}+D} -\rho \right)
    \sqrt{\frac{\eps^2\sigma}{D}} g_0(\eps^2\sigma \tau) \right]\\
  & \qquad + \frac{\eps d_{2j} u_{1j}^{\prime}(0)}{\rho (d_{1j}+D)}
  \left[ \tau + \left( \frac{d_{1j}}{d_{1j}+D} -\rho \right)
    \frac{g_1(\eps^2\sigma \tau)}{\sqrt{\eps^2 \sigma D}} \right]
  \,.
  \end{split}
\end{equation}

To match to an outer solution, we consider the intermediate time
scale regime ${\mathcal O}(1)\ll \tau\ll {\mathcal O}(\eps^{-2})$. By using
the Taylor series approximations
\begin{equation}\label{3:g0g1loc}
  g_0(\xi) \sim \frac{1}{\sqrt{\pi \xi}} + \sqrt{\frac{\xi}{\pi}} \,,
  \qquad
  g_{1}(\xi)\sim 2 \sqrt{\frac{\xi}{\pi}} + \frac{2}{3\sqrt{\pi}}\xi^{3/2}
    \,, \quad \mbox{for}\quad \xi\ll 1 \,,
\end{equation}
in \eqref{3:v_asy_1}, we obtain that
\begin{equation}\label{3:v_asy_2}
  \begin{split}    
   V(\rho,\tau) &\sim \frac{d_{2j} u_{1j}(0)}{\eps \rho (d_{1j}+D)}
   \left[ 1 + \frac{1}{\sqrt{\pi D\tau}}
     \left( \frac{d_{1j}}{d_{1j}+D} -\rho \right)
   \left(1 + \eps^2 \sigma \tau\right) \right]\\
  & \qquad + \frac{\eps d_{2j} u_{1j}^{\prime}(0)}{\rho (d_{1j}+D)}
  \left[ \tau + 2 \sqrt{\frac{\tau}{\pi D}}
    \left( \frac{d_{1j}}{d_{1j}+D} -\rho \right)
    \left(1 + \frac{\eps^2\tau \sigma}{3}\right) \right]
  \,.
  \end{split}
\end{equation}
Finally, to obtain our matching condition for the outer solution, as
derived in \S \ref{2:ode_memory}, we set $\tau={t/\eps^2}$ in
\eqref{3:v_asy_2} to obtain the result in \eqref{3:v_final}.

\section{Asymptotic Analysis for Steady-States}\label{app:strong}

In this appendix, we use a strong localized perturbation analysis
(see \cite{ward2018spots}) applied to the cell-bulk model
\eqref{DimLess_bulk} to establish that the limiting NAS \eqref{4:limit_sys}
characterizes all possible steady-states of \eqref{DimLess_bulk}.

In the limit $\eps\to 0$, the steady-state problem for the outer bulk
solution $U_s(\v{x})$ is
\begin{equation}\label{4:ssbulk}
  \begin{split}
   &  D \Delta U_{s} - \sigma U_{s} =0\,, \quad \v{x}\in \R^3\backslash
   \lbrace{\v{x}_1, \ldots,\v{x}_N\rbrace} \,; \quad U_s(\v{x})\to 0 \,\,\,
   \mbox{as}\,\,\, |\v{x}|\to \infty \,,\\
 &    U_s \sim \frac{\kappa_j u_{1js}}{4\pi D|\x-\x_j|} + B_{js}
    \,, \quad \mbox{as} \quad \v{x}\to \v{x}_j \,, \quad
    j\in\lbrace{1,\ldots,N\rbrace}\,,
  \end{split}
\end{equation}
where $\kappa_j$ is given in \eqref{2:reduced_4} and where the
behavior as $\v{x}\to \v{x}_j$ was obtained by matching to the local
or inner solution. From the steady-state of the intracellular kinetics
in \eqref{DimLess_bulk} we obtain that
\begin{equation}\label{4:ssintra}
  \v{F}_{j}\left(\v{u}_{js}\right) -\kappa_j u_{1js}\v{e}_1 +
  \eps \frac{d_{1j}}{d_{2j}} \kappa_j B_{js} \v{e}_1 =0 \,, \qquad
  j\in \lbrace{1,\ldots,N\rbrace}\,.
\end{equation}
We represent the solution to \eqref{4:ssbulk} in terms of the 
Green's function of \eqref{4:ge} as
\begin{equation}\label{4:ssbulk_sol}
  U_{s}(\v{x}) = -\sum_{k=1}^{N} \frac{\kappa_k u_{k1s}}{D}
  G_{0}(\x,\x_k) \,.
\end{equation}
By enforcing that the local behavior of $U_{s}(\v{x})$ as
$\v{x}\to \v{x}_j$ agrees with that specified in \eqref{4:ssbulk} we
readily determine $B_{js}$ from the regular part of the resulting expression
as
\begin{equation}\label{4:ss_nas}
  B_{js} =  -\frac{\kappa_j}{D} R_{0}   u_{1js} 
  - \sum_{\stackrel{k=1}{k\neq j}}^{N} \frac{\kappa_k}{D} u_{k1s}\, G_{0jk}
  \,, \qquad j\in \lbrace{1,\ldots,N\rbrace} \,,
\end{equation}
where $G_{0jk}$ and $R_0$ are given in \eqref{4:lim_sys_green}.  We
conclude that the coupled system \eqref{4:ssintra} and \eqref{4:ss_nas}
is the limiting system \eqref{4:limit_sys} obtained from the
long-time dynamics of \eqref{2:reduced}.

Finally, we combine \eqref{4:ssintra} and \eqref{4:ss_nas} to obtain a
NAS for the intracellular steady-states. To do so, we let
$\v{B}_{s}=(B_{1s},\ldots,B_{Ns})^T$ and
$\v{u}_{1s}=(u_{11s},\ldots,u_{1Ns})^T$, and write \eqref{4:ss_nas} in
matrix form as
\begin{subequations}\label{4:b0eq}
\begin{equation}\label{4:b0eq_1}
  \v{B}_{s} = -\frac{1}{D} {\mathcal G}_0 {\mathcal K} \v{u}_{1s}\,,
\end{equation}
where the steady-state Green's matrix ${\mathcal G}_0$ and the
diagonal matrix ${\mathcal K}$ are defined by
\begin{equation}\label{4:b0eq_2}
  ({\mathcal G}_{0})_{ij} = ({\mathcal G}_{0})_{ji} =G_{0ij}\,,
 \quad i\neq j \,, \quad ({\mathcal G}_{0})_{jj} = R_{0} \,; \quad
  {\mathcal K}\equiv \mathrm{diag}\left(\kappa_1,\ldots,\kappa_N\right)\,.
\end{equation}
\end{subequations}
Since $B_{js}=\v{e}_j^{T}\v{B}_{s}$, where $\v{e}_j$ is the $j^{\mathrm{th}}$ unit
basis vector, we substitute \eqref{4:b0eq_1} into
\eqref{4:ssintra} to obtain the following $m\times N$ dimensional NAS
for the steady-state intracellular species, which is equivalent to
\eqref{4:limit_sys}:
\begin{equation}\label{4:ss_final}
  \v{F}_{j}\left(\v{u}_{js}\right) -\kappa_j u_{1js}\v{e}_1 -
  \eps \frac{d_{1j}}{D d_{2j}} \kappa_j \v{e}_1 \left(
    \v{e}_j^{T} {\mathcal G}_0 {\mathcal K} \v{u}_{1s} \right) = 0\,, \qquad
    j\in \lbrace{1,\ldots,N\rbrace}\,.
\end{equation}

\section{Numerical Treatment for the Integro-Differential System}\label{app:numerics}

In this appendix, we formulate a time-marching scheme to
efficiently compute solutions to \eqref{2:reduced} over long time-intervals.

To do so, we first outline the derivation of an update scheme for the
two nonlocal terms in \eqref{alg:CD:def} that combines the SOE
approximation for their kernels together with Duhamel's principle.
Our use of Duhamel's principle in deriving an approximation scheme
relies on the following simple result.

\begin{lemma}\label{prop:duh} Let $f(t)$ be continuous and define 
  ${\mathcal F}(t)\equiv\int_{0}^{t} e^{\omega (t-\tau)} f(\tau)\,
  d\tau$.  Then, 
  ${\mathcal F}^{\prime}(t)=\omega {\mathcal F}(t) + f(t)$ with
  ${\mathcal F}(0)=0$. Moreover, we have the marching scheme
\begin{equation}\label{duh:update}
    {\mathcal F}(t+\Delta t)={\mathcal F}(t) e^{\omega \Delta t}
    +\, {\mathcal U}(t,\Delta t) \,, \quad \mbox{with} \quad
    {\mathcal U}(t,\Delta t)\equiv e^{\omega \Delta t} \int_{0}^{\Delta t}
    e^{-\omega z} f(t+z) \, dz \,.
\end{equation}
An {\em exponential time differencing} ETD2 scheme
(cf.~\cite{lopez_1}), for which the update integral ${\mathcal U}$
is exact for linear functions $f(t)$, yields, with an error
${\mathcal O}\left( (\Delta t)^3\right)$, the approximation
\begin{equation}\label{duh:etd2}
  {\mathcal F}(t+\Delta t) \approx {\mathcal F}(t) e^{\omega \Delta t}
  + f(t)\left( \frac{e^{\omega \Delta t}-1}{\omega}\right) +
  \left[ f(t+\Delta t)-f(t)\right] \left( \frac{e^{\omega \Delta t} - 1 -
      \omega \Delta t}{\omega^2\Delta t}\right) \,.
\end{equation}
\end{lemma}

\begin{proof} The proof of \eqref{duh:update} is immediate. To derive
  \eqref{duh:etd2} we substitute
  \begin{equation*}
    f(t+z) = f(t) +\frac{z}{\Delta t} \left[ f(t+\Delta t)-f(t)\right]
    +{\mathcal O}\left((\Delta t)^2\right) \,,
  \end{equation*}
  into the update integral ${\mathcal U}$ in \eqref{duh:update} and
  integrate the resulting expression explicitly.
\end{proof} 

\subsection{Update scheme for \texorpdfstring{$C_j(t)$}{C}}\label{app:numerics_cj}
Inspired by \cite{greengard_2015}, and the analysis of the
corresponding 2-D problem \cite{pelzward2D}, we begin by decomposing
$C_j(t)$ into local $C_{Lj}(t)$ and history $C_{Hj}(t)$ contributions
as
\begin{subequations}\label{alg:c_decomp}
\begin{equation}
  C_j(t) = C_{Hj}(t) + C_{Lj}(t) \,,
\end{equation}
where for some $\Delta t$, with $0<\Delta t\ll 1$, we define 
\begin{align}
   C_{Hj}(t) &\equiv \int_0^{t-\Delta t} u_{1j}^{\prime}(\tau)
                 g(t-\tau) \, d\tau\,, \label{alg:c_decomp_h} \\
  C_{Lj}(t) &\equiv \int_{t-\Delta t}^t u_{1j}^{\prime}(\tau)
  g(t-\tau)\,  d\tau =
  \int_{0}^{\Delta t} u_{1j}^{\prime}(t-z) g(z)\, dz\,.
  \label{alg:c_decomp_l}
\end{align}
\end{subequations}
Here $g(t)$ is defined in \eqref{2:reduced_3}.  The history term in
\eqref{alg:c_decomp_h} is updated using the SOE approximation for
$g(t)$ in \eqref{alg:soe} combined with Duhamel's result
\eqref{duh:etd2} for each term in the sum.  To treat the integrable
singularity in the local contribution \eqref{alg:c_decomp_l} we use
the local behavior of $g(t)$. The following lemma summarizes the
contributions of the local and history parts.

\begin{lemma} \label{lem:Cj_update}
  Fix $\Delta t\ll 1$. For $t=\Delta t$, $t=2\Delta t$, and
  $t>2\Delta t$ we have
\begin{equation} \label{eq:CLj_update}
  \begin{split}
    C_{j}(\Delta t) &\approx b_{C1} u_{1j}(\Delta t) - b_{C1} u_{1j}(0)
      \,, \\
   C_j(2\Delta t) &\approx b_1 u_{1j}^{\prime}(2\Delta t) +
   b_2 u_{1j}^{\prime}(\Delta t) +
   \left(\sum_{\ell=-n}^n e_\ell b_{C30\ell}\right) \left( u_{1j}(\Delta t)-
     u_{1j}(0) \right) \,, \\
  C_j(t+\Delta t) &\approx b_1 u_{1j}^{\prime}(t+\Delta t) + \left(b_2 +
    \sum_{\ell=-n}^n e_l b_{3\ell}\right) u_{1j}^{\prime}(t) \\
   & \qquad + \left(\sum_{\ell=-n}^n e_\ell b_{4\ell}\right)
   u_{1j}^{\prime}(t-\Delta t) + \sum_{\ell=-n}^n e_\ell e^{s_\ell \Delta t}
   H_{Cj\ell}(t)\,,
 \end{split}
 \end{equation}
 where $e_{\ell}$ is defined in \eqref{alg:el_omega} and the other
 coefficients are defined by
\begin{equation}\label{alg:C:values}
  \begin{split}
      b_1 &\equiv \frac{\sqrt{\sigma} \Delta t}{2}  \,
    \mathrm{erf}\left(\sqrt{\sigma\Delta t}\right)
    \left(-\frac{1}{4\sigma^2 (\Delta t)^2} + \frac{1}{\sigma \Delta t}
      + 1 \right)  +
    \frac{e^{-\sigma \Delta t}}{4\sigma \sqrt{\pi \Delta t}}\left( 1 +
      2\sigma \, \Delta t\right)\,,\\
    b_2 &\equiv \frac{\sqrt{\sigma}\Delta t}{2}  \,
    \mathrm{erf}\left(\sqrt{\sigma \Delta t}\right)
    \left( \frac{1}{4\sigma^2 (\Delta t)^2} + 1\right) -
    \frac{e^{-\sigma \Delta t}}{4\sigma \sqrt{\pi \Delta t}}\left( 1-
      2\sigma \, \Delta t\right)\,,\\
            b_{3\ell} &\equiv \frac{e^{s_{\ell}\Delta t}}{s_{\ell}^2 \Delta t}
         \left(e^{s_{\ell}\Delta t} -1 -s_{\ell}\Delta t\right)\,,
         \qquad b_{4\ell} \equiv \frac{e^{s_{\ell}\Delta t}}{s_{\ell}}
   \left(e^{s_{\ell}\Delta t}-1\right) - b_{3\ell}\,, \\
   b_{C0\ell} &\equiv \frac{1}{s_{\ell}\Delta t}\left(e^{2s_\ell\Delta t} - e^{s_\ell
       \Delta t}\right)\,, \qquad b_{C1} =\sqrt{\sigma} \mathrm{erf}
   \left( \sqrt{\sigma \Delta t}\right) \left( \frac{1}{2\sigma \Delta t}
     + 1 \right) + \frac{e^{-\sigma \Delta t}}{\sqrt{\pi \Delta t}} \,.
\end{split}
\end{equation}
Moreover, with $H_{Cj\ell}(\Delta t)=0$, we have the update scheme
\begin{equation}\label{alg:Hc_update}
  H_{Cj\ell}(t+\Delta t) = H_{Cj\ell}(t) e^{s_{\ell}\Delta t} + 
   \begin{cases}
b_{C0\ell} \, u_{1j}(\Delta t) - b_{C0\ell} \, u_{1j}(0)\,, & t=\Delta t\,, \\
	b_{3\ell} \, u_{1j}^{\prime}(t) + b_{4\ell} \, u_{1j}^{\prime}(t-\Delta t)\,,
        & t\geq 2\Delta t \,.
    \end{cases}
\end{equation}
\end{lemma}

\begin{proof}
  For the \textbf{first step} consider the local contribution
  $C_{Lj}(t)$ in \eqref{alg:c_decomp}.  For $z\in[0, \Delta t]$ and
  $t \geq 2\Delta t$ the ETD2 approximation \cite{lopez_1} yields
\begin{equation*}
  u_{1j}^{\prime}(t-z) = u_{1j}^{\prime}(t) +
  \frac{z}{\Delta t} \left[u_{1j}^{\prime}(t-\Delta t) -
          u_{1j}^{\prime}(t)\right] + \mathcal{O}(\Delta t^2)\,,
\end{equation*}
which, when substituted into the local integral, yields that
\begin{equation*}
  C_{Lj}(t) = b_1 u_{1j}^{\prime}(t) + b_2 u_{1j}^{\prime}(t-\Delta t) +
  \mathcal{O}\left( (\Delta t)^2\right) \,,
\end{equation*}
where $b_1$ and $b_2$ are defined by
\begin{equation}\label{alg:C:b1_b2}
  b_1 \equiv \int_{0}^{\Delta t} g(z)\, dz -\frac{1}{\Delta t}
  \int_{0}^{\Delta t} z g(z) \, dz \,, \qquad
  b_2 \equiv \frac{1}{\Delta t}\int_{0}^{\Delta t} z g(z) \, dz \,.
\end{equation}
From the detailed calculation in Appendix \ref{app:local_int}, we obtain that
\begin{equation}\label{alg:c_locint}
\begin{split}
  \int_0^{\Delta t} g(z) \, dz &= \frac{1}{2\sqrt{\sigma}} \;
  \mathrm{erf}(\sqrt{\sigma\Delta t}) + \sqrt{\sigma}\Delta t \,
  \mathrm{erf}(\sqrt{\sigma\Delta t}) + \sqrt{\frac{\Delta t}{\pi}}
                                 e^{-\sigma \Delta t}\,, \\
   \int_0^{\Delta t} z g(z) \, dz &= \frac{1}{8\sigma^{3/2}}
   \mathrm{erf}(\sqrt{\sigma \Delta t}) + \frac{\sqrt{\sigma}}{2}
   \left(\Delta t\right)^2 \mathrm{erf}(\sqrt{\sigma \Delta t}) -
   \frac{1}{4\sigma} \sqrt{\frac{\Delta t}{\pi}} e^{-\sigma \Delta t} +
   \frac{ \left(\Delta t \right)^{3/2}}{2\sqrt{\pi}} e^{-\sigma \Delta t}\,.
 \end{split}
\end{equation}
In this way, by using \eqref{alg:C:b1_b2} we
obtain the results for $b_1$ and $b_2$ in \eqref{alg:C:values}.

For $t=\Delta t$,
$C_{j}(\Delta t)=C_{Lj}(\Delta t) = \int_{0}^{\Delta t}
u_{1j}^{\prime}(\xi) g(\Delta t - \xi)\, d\xi$. By using
$u_{1j}^{\prime}(\xi)\approx {\left[u_{1j}(\Delta
    t)-u_{1j}(0)\right]/\Delta t}$ on $0\leq \xi\leq \Delta t$, we
obtain
\begin{equation*}
  C_{Lj}(\Delta t) \approx \frac{ \left[ u_{1j}(\Delta t)-u_{1j}(0)\right]}{
    \Delta t} \int_{0}^{\Delta t} g(z) \, dz \,.
\end{equation*}
Then, upon using \eqref{alg:c_locint} for
$\int_{0}^{\Delta t} g(z)\, dz$, we obtain the result in the first
line of \eqref{eq:CLj_update}.

For the \textbf{second step} we use the SOE approximation in
\eqref{alg:soe} for $g(t)$ to calculate the history term:
\begin{equation*}
  C_{Hj}(t) \approx \int_0^{t-\Delta t} u_{1j}^{\prime}(\tau) g_n(t-\tau) \; d\tau
  = \sum_{\ell=-n}^n e_{\ell} H_{Cj\ell}(t) \quad \mbox{with} \quad
  H_{Cj\ell}(t) \equiv\int_0^{t-\Delta t}
  u_{1j}^{\prime}(\tau) e^{s_{\ell}(t-\tau)}\;  d\tau\,.
\end{equation*}
Lemma \ref{prop:duh} yields the updating rule
\begin{equation*}
  H_{Cj\ell}(t+\Delta t) = H_{Cj\ell}(t) e^{s_{\ell}\Delta t} +
  \mathcal{U}_{Cj\ell}(t, \Delta t)\,, \qquad H_{Cj\ell}(\Delta t) = 0 \,,
\end{equation*}
where for $t\geq 2\Delta t$, we have
\begin{equation*}
  \mathcal{U}_{Cj\ell}(t, \Delta t) \equiv e^{2s_{\ell}\Delta t}
  \int_0^{\Delta t} e^{-s_{\ell} z}
  u_{1j}^{\prime}(t-\Delta t + z)\; dz \approx b_{3\ell} u_{1j}^{\prime}(t) +
  b_{4\ell} u_{1j}^{\prime}(t-\Delta t)\,.
\end{equation*}
Here $b_{3\ell}$ and $b_{4\ell}$, as defined in \eqref{alg:C:values}, are
obtained from an ETD2 truncation of $\mathcal{U}_{Cj\ell}$. Moreover,
for $t=\Delta t$, we obtain with integration by parts
\begin{equation*}
  \mathcal{U}_{Cj\ell}(\Delta t, \Delta t) \approx b_{C0\ell}\,
  u_{1j}(\Delta t) -  b_{C0\ell} \, u_{1j}(0)\,,
\end{equation*}
where $b_{C0\ell}$ is defined in \eqref{alg:C:values}.  By adding the
local and history contributions we obtain the update scheme in
\eqref{eq:CLj_update} and \eqref{alg:Hc_update}.
\end{proof}

\subsection{Update scheme for \texorpdfstring{$D_{jk}(t)$}{D}}
\label{app:numerics_dj}
Next, we derive the update scheme for the history-dependent
integral $D_{jk}(t)$, as defined in \eqref{alg:CD:def}, in a very similar
way. First we separate it into a local part $D_{Ljk}(t)$ and a history
part $D_{Hjk}(t)$ as
\begin{subequations}\label{alg:d_decomp}
\begin{equation}
  D_{jk}(t) = D_{Hjk}(t) + D_{Ljk}(t) \,,
\end{equation}
where for some $\Delta t$, with $0<\Delta t\ll 1$, we define
\begin{align}
&  D_{Hjk}(t) \equiv \int_0^{t-\Delta t} u_{k1}(\tau)
  G(a_{jk},t-\tau) \, d\tau\,,  \label{alg:d_decomp_h} \\
&  D_{Ljk}(t) \equiv \int_{t-\Delta t}^t u_{k1}(\tau)
        G(a_{jk},t-\tau)\, d\tau = \int_0^{\Delta t} u_{k1}(t-z)
     G(a_{jk},z)\, dz\,.  \label{alg:d_decomp_l}
\end{align}
\end{subequations}
Here $G(a_{jk},t)$ is defined in \eqref{2:reduced_3}. Once again, the
history term in \eqref{alg:d_decomp_h} is updated using the SOE
approximation for $G(a_{jk},t)$ in \eqref{alg:soe} combined with
Duhamel's result \eqref{duh:etd2}. The local contribution
\eqref{alg:d_decomp_l}, near the integrable singularity, is estimated
analytically.  The updating scheme is summarized in the following
lemma.

\begin{lemma} \label{lem:Djk_update} Fix $\Delta t\ll 1$. For
  $t=\Delta t$, $t=2\Delta t$, and $t>2\Delta t$, we have the update
  scheme
\begin{equation} \label{eq:Djk_update}
  \begin{split}
  D_{jk}(\Delta t) & \approx c_{1jk} u_{k1}(\Delta t) + c_{2jk} u_{k1}(0) \,,\\
  D_{jk}(2\Delta t) &\approx c_{1jk} u_{k1}(2\Delta t) + \left( c_{2jk} +
    \sum_{\ell=-n}^n \omega_{jk\ell} \, b_{3\ell}\right) u_{k1}(\Delta t)+
  \left(\sum_{\ell=-n}^n \omega_{jk\ell} \, b_{4\ell}\right) u_{k1}(0)\,,\\
  D_{jk}(t+\Delta t) &\approx c_{1jk} u_{k1}(t+\Delta t) + \left(c_{2jk}
    + \sum_{\ell=-n}^n \omega_{jk\ell} \, b_{3\ell}\right) u_{k1}(t) \\
  & \qquad + \left(\sum_{\ell=-n}^n\omega_{jk\ell} \, b_{4\ell}\right)
  u_{k1}(t-\Delta t) + \sum_{\ell=-n}^n \omega_{jk\ell} e^{s_{\ell}\Delta t}
  H_{Dk\ell}(t) \,,
\end{split}
\end{equation}
where $\omega_{jkl}$ was defined in \eqref{alg:el_omega} while
$b_{3\ell}$ and $b_{4\ell}$ were defined in \eqref{alg:C:values}. The
new coefficients $c_{1jk}$ and $c_{2jk}$ are
\begin{equation}\label{alg:D:values}
  \begin{split}
    c_{2jk} &\equiv \frac{1}{(4\pi D)^{3/2}} \left(
      \frac{\sqrt{\pi}}{2\Delta t \sqrt{\sigma}}\right)
    \left[ e^{-a_{jk}\sqrt{\sigma}} \mathrm{erfc}\left(
          \frac{a_{jk}}{2\sqrt{\Delta t}} -\sqrt{\sigma \Delta t}\right)
 - e^{a_{jk}\sqrt{\sigma}} \mathrm{erfc}\left(
   \frac{a_{jk}}{2\sqrt{\Delta t}} + \sqrt{\sigma \Delta t}\right)
\right] \,, \\
c_{1jk} &\equiv \frac{\sqrt{\pi}}{(4\pi D)^{3/2}a_{jk}}
\left[ e^{-a_{jk}\sqrt{\sigma}} \mathrm{erfc}\left(
          \frac{a_{jk}}{2\sqrt{\Delta t}} -\sqrt{\sigma \Delta t}\right)
 + e^{a_{jk}\sqrt{\sigma}} \mathrm{erfc}\left(
   \frac{a_{jk}}{2\sqrt{\Delta t}} + \sqrt{\sigma \Delta t}\right)
\right] - c_{2jk} \,,
\end{split}
\end{equation}
where $\mathrm{erfc}(z)$ is the complementary error function.
Moreover, for $t\ge \Delta t$, and with $H_{D\ell}(\Delta t)=0$, we
have 
\begin{equation}\label{alg:Hd_update}
  H_{Dk\ell}(t+\Delta t) = H_{Dk\ell}(t) e^{s_{\ell}\Delta t} +
 b_{3\ell}\, u_{k1}(t) + b_{4\ell} \, u_{k1}(t-\Delta t)\,.
\end{equation}
\end{lemma}

\begin{proof}
  The \textbf{first step} is to consider the local contribution
  $D_{Ljk}(t)$. For $z\in[0, \Delta t]$ and $t \geq \Delta t$ the
  ETD2  approximation \cite{lopez_1} gives 
\begin{equation*}
  u_{k1}(t-z) = u_{k1}(t) + \frac{z}{\Delta t}
  \left[ u_{k1}(t-\Delta t) - u_{k1}(t)\right] + \mathcal{O}(\Delta t^2)\,.
\end{equation*}
Upon substituting this expression into the local integral
\eqref{alg:d_decomp_l} we obtain
\begin{equation}\label{alg:D:local_e}
  D_{Ljk}(t) \approx u_{k1}(t) c_{1jk}   + u_{k1}(t-\Delta t) c_{2jk} \,,
\end{equation}
where
\begin{equation}\label{alg:D:I1_I2}
  c_{1jk}  = I_{1} - \frac{I_2}{\Delta t} \,, \quad c_{2jk}=
  \frac{I_2}{\Delta t}\,, \qquad \mbox{with} \qquad
  I_1\equiv \int_0^{\Delta t} G(a_{jk}, z)\, dz   \,, \quad
  I_2\equiv \int_0^{\Delta t} z G(a_{jk}, t)\,dz\,.
\end{equation}
As shown in Appendix \ref{app:local_int}, the two integrals
$I_1$ and $I_2$ can be evaluated analytically as
\begin{equation}\label{alg:D:Ivalues}
  \begin{split}
    I_1 &= \frac{1}{(4\pi D)^{3/2}} \left(
      \frac{\sqrt{\pi}}{a_{jk}} \right)
    \left[ e^{-a_{jk}\sqrt{\sigma}} \mathrm{erfc}\left(
          \frac{a_{jk}}{2\sqrt{\Delta t}} -\sqrt{\sigma \Delta t}\right)
 + e^{a_{jk}\sqrt{\sigma}} \mathrm{erfc}\left(
   \frac{a_{jk}}{2\sqrt{\Delta t}} + \sqrt{\sigma \Delta t}\right)
\right] \,, \\
I_2 &= \frac{1}{(4\pi D)^{3/2}} \left(\frac{\sqrt{\pi}}{2\sqrt{\sigma}}
\right)
\left[ e^{-a_{jk}\sqrt{\sigma}} \mathrm{erfc}\left(
          \frac{a_{jk}}{2\sqrt{\Delta t}} -\sqrt{\sigma \Delta t}\right)
 - e^{a_{jk}\sqrt{\sigma}} \mathrm{erfc}\left(
   \frac{a_{jk}}{2\sqrt{\Delta t}} + \sqrt{\sigma \Delta t}\right)
\right] \,,
\end{split}
\end{equation}
which yields the result in \eqref{alg:D:values} for
$c_{1jk}$ and $c_{2jk}$.

\textbf{Secondly}, to calculate the history term we use the SOE
approximation in \eqref{alg:soe} for $G(a_{jk},t)$ to get
\begin{equation}
  D_{Hjk}(t) \approx \int_0^{t-\Delta t} u_{k1}(\tau) G_n(a_{jk}, t-\tau) \, d\tau
  = \sum_{l=-n}^n \omega_{jkl} H_{Dkl}(t) \,; \qquad
 H_{Dk\ell}(t) \equiv \int_0^{t-\Delta t} u_{k1}(\tau)
  e^{s_{\ell}(t-\tau)}\; d\tau \,.
\end{equation}
Lemma \ref{prop:duh} provides the updating rule
\begin{equation*}
  H_{Dk\ell}(t+\Delta t) = H_{Dk\ell}(t) e^{s_{\ell}\Delta t} +
  \mathcal{U}_{Dk\ell}(t, \Delta t)\,, \qquad H_{Dk\ell}(\Delta t) = 0\,,
\end{equation*}
where the update integral $\mathcal{U}_{Dk\ell}(t, \Delta t) \equiv
e^{2s_{\ell}\Delta t} \int_0^{\Delta t} e^{-s_{\ell} z}u_{k1}(t-\Delta t + z)\; dz$
is estimated with an ETD2 approximation for $t\geq \Delta t$ as
\begin{equation*}
  \mathcal{U}_{Dk\ell}(t, \Delta t) \approx 
       b_{3l} \, u_{k1}(t) + b_{4l}\, u_{k1}(t-\Delta t)\,.
\end{equation*}
The constants $b_{3\ell}$ and
$b_{4\ell}$ are as in Lemma \ref{lem:Cj_update} and are given in
\eqref{alg:C:values}.  By adding the local and history contributions
we obtain the update scheme in \eqref{eq:Djk_update} and
\eqref{alg:Hd_update}.
\end{proof}

\subsection{Time-marching scheme}\label{app:numerics_ode}

In terms of our updating scheme for the memory-dependent integrals
$C_j(t)$ and $D_{jk}(t)$ given in Lemmas \ref{lem:Cj_update} and
\ref{lem:Djk_update}, we now provide the full scheme for numerically
solving our integro-differential system \eqref{2:reduced} for the case
of two intracellular species. In terms of $C_j(t)$ and $D_{jk}(t)$ of
\eqref{alg:CD:def}, and recalling
$\kappa_j\equiv {4\pi d_{2j} D/(d_{1j}+D)}$, this system for
$j\in \lbrace{1,\ldots,N\rbrace}$ is
\begin{subequations}\label{tm:ode-alg}
\begin{align}
    \frac{d\v{u}_j}{dt} &= \v{F}_j(\v{u}_j)- \kappa_j u_{1j} \v{e}_1 +
                          \frac{\eps d_{1j}}{d_{2j}} \kappa_j B_{j}(t)
  \v{e}_1 \,, \label{tm:ode} \\
    B_{j}(t) &= -\frac{\kappa_j u_{1j}(0)}{4\pi D^{3/2}} g(t) -
               \frac{\kappa_j}{4\pi D^{3/2}} C_{j}(t)
               + \sum_{\stackrel{k=1}{k\neq j}}^{N} \kappa_k
     D_{jk}(t) \,.\label{tm:alg}
\end{align}
\end{subequations}

Let $\{t_0, ..., t_q\}$ be a uniform discretization of the time domain
$[0, T]$, where $T$ is some final chosen time, and where $t_0=0$,
$t_q=T$, and $t_{i+1} - t_i = \Delta t$ for $i\in\{0, ..., q-1\}$,
with $\Delta t = {T/q}$. Below, we denote the numerical approximation
to $u_{1j}(t_i)$ by $u_{1j}^{(i)}$ and the vector of all
$u_{1j}^{(i)}$ for $j\in\{1,\ldots,N\}$ by $\v{u}_1^{(i)}$. We use a
similar notation for the reaction kinetic
functions. Further, in the description of the algorithm
  below, $[ \;\cdot_j \;]_j$ denotes a diagonal matrix of size
  $N\times N$ with diagonal entries $(\;\cdot_j\;)_j$, while
  $[ \;\cdot_{jk} \;]_{jk}$ denotes an $N\times N$ matrix, with matrix
  entries $(\;\cdot_{jk}\;)_{jk}$ for $j\neq k$ and with zero diagonal
  entries.

\textbf{Step 1:} $0 \curvearrowright \Delta t$: For this first
time-step we discretize the reaction kinetics in each cell by the
explicit Runge-Kutta RK4 method \cite{numerics} and use
\begin{equation}\label{alg:step0}
  \begin{split}
  \v{u}_{1}^{(1)} &= \v{u}_{1}^{(0)} + \frac{\Delta t}{6} (k_1^{(0)} +
  2 k_2^{(0)} + 2 k_3^{(0)} + k_4^{(0)}) + \varepsilon \left[\frac{d_{1j}}{d_{2j}}
    \kappa_j \right]_{j}\v{B}^{(1)} \Delta t\,, \\
  \v{u}_{2}^{(1)} &= \v{u}_{2}^{(0)} + \frac{\Delta t}{6} (\tilde{k}_1^{(0)}
  + 2 \tilde{k}_2^{(0)} + 2 \tilde{k}_3^{(0)} + \tilde{k}_4^{(0)}) \,,
\end{split}
\end{equation}
where the RK4 weights are
\begin{eqnarray*}
	k_1^{(0)} &=& \hat{\v{F}}_1(\v{u}_{1}^{(0)}, \v{u}_{2}^{(0)})\,, \quad
	\tilde{k}_1^{(0)}= \v{F}_2(\v{u}_{1}^{(0)}, \v{u}_{2}^{(0)})\,, \\
  k_2^{(0)} &=&  \hat{\v{F}}_1(\v{u}_{1}^{(0)} + \frac{\Delta t}{2} k_1^{(0)}\,,
                \v{u}_{2}^{(0)} + \frac{\Delta t}{2} \tilde{k}_1^{(0)})\,, \quad
  \tilde{k}_2^{(0)} =\v{F}_2(\v{u}_{1}^{(0)} + \frac{\Delta t}{2} k_1^{(0)},
                     \v{u}_{2}^{(0)} + \frac{\Delta t}{2} \tilde{k}_1^{(0)})\,, \\
  k_3^{(0)} &=&  \hat{\v{F}}_1(u_{1}^{(0)} + \frac{\Delta t}{2} k_2^{(0)},
                \v{u}_{2}^{(0)} + \frac{\Delta t}{2} \tilde{k}_2^{(0)})\,, \quad
  \tilde{k}_3^{(0)} = \v{F}_2(\v{u}_{1}^{(0)} + \frac{\Delta t}{2} k_2^{(0)}\,,
                   \v{u}_{2}^{(0)} + \frac{\Delta t}{2} \tilde{k}_2^{(0)})\,, \\
  k_4^{(0)} &=&  \hat{\v{F}}_1(\v{u}_{1}^{(0)} + \Delta t \, k_3^{(0)},
                \v{u}_{2}^{(0)} +  \Delta t \, \tilde{k}_3^{(0)})\,, \quad
  \tilde{k}_4^{(0)}  =\v{F}_2(\v{u}_{1}^{(0)} + \Delta t \, k_3^{(0)},
                        \v{u}_{2}^{(0)} + \Delta t \, \tilde{k}_3^{(0)})\,,
\end{eqnarray*}
and where $\hat{\v{F}}_1(\v{u}_1, \v{u}_2) =
\v{F}_{1}(\v{u}_1, \v{u}_2) - \left[ \kappa_j\right]_{j} \v{u}_1$.

To determine the initial values for the recursion, we evaluate the
matching condition \eqref{3:v_bj} for $B_j$ at $t=\Delta t$. In vector
form, this determines our approximation of $\v{B}^{(1)}$ for any
$\v{u}_1^{(0)}$ as
\begin{equation} \label{eq:B1added}
  \v{B}^{(1)} = -\frac{1}{4 (\pi D)^{3/2} \sqrt{\Delta t}}
  \left[\kappa_j\right]_{j} \left( \v{u}_{1}^{(0)} +
    \Delta t \left( \sigma \v{u}_{1}^{(0)} + 2 \v{u}_{1}^{\prime (0)} \right)
    \right) \,,
 \end{equation}
where the components of $\v{u}_{1}^{\prime (0)}$ are given in \eqref{3:uprime}.
The truncation error for this approximation of $B_j$ is
$\mathcal{O}(\Delta t)$, whereas the RK4 scheme gives a truncation
error of $\mathcal{O}(\Delta t^5)$ for the reaction kinetics.

\textbf{Step 2:} $\Delta t \curvearrowright 2\Delta t$: We use an RK4 scheme
for the reaction kinetics with RK4 weights as given in the first time step,
and use a lagged $\v{B}^{(1)}$ so that
\begin{equation}\label{alg:step2}
  \begin{split}
    \v{u}_{1}^{(2)} &= \v{u}_{1}^{(1)} + \frac{\Delta t}{6}
    (k_1^{(1)} + 2 k_2^{(1)} + 2 k_3^{(1)} + k_4^{(1)}) + \varepsilon
    \left[\frac{d_{1j}}{d_{2j}} \kappa_j \right]_{j} \v{B}^{(1)} \Delta t\,, \\
    \v{u}_{2}^{(2)} &= \v{u}_{2}^{(1)} + \frac{\Delta t}{6}
    (\tilde{k}_1^{(1)} + 2 \tilde{k}_2^{(1)} + 2 \tilde{k}_3^{(1)} +
    \tilde{k}_4^{(1)})\,.
  \end{split}
\end{equation}
In terms of the computed $\v{u}_{1}^{(2)}$, we determine $\v{B}^{(2)}$
from the linear system
\begin{equation}\label{alg:b2}
  \begin{split}
  \v{B}^{(2)} &= -\frac{1}{4\pi D^{\frac{3}{2}}} A^{-1}
    \left[ \kappa_j\right]_{j} \left(\varepsilon b_2
      \left[\frac{d_{1j}\kappa_j}{d_{2j}}\right]_{j} \v{B}^{(1)} +
      b_1 \hat{\v{F}}_1(\v{u}^{(2)}) + b_2 \hat{\v{F}}_1(\v{u}^{(1)})\right)
   \\
  & \qquad\qquad + A^{-1}\left(N \v{u}_1^{(2)} + M_1 \v{u}_1^{(1)} +
    M_0 \v{u}_1^{(0)}\right)\,,
  \end{split}
\end{equation}
where $\v{B}^{(1)}$ is given in
\eqref{eq:B1added}. In \eqref{alg:b2},  we have defined the matrices
$A$, $M_1$, $M_0$, and $N$ by
\begin{eqnarray*}
  &&A \equiv I + \frac{\varepsilon b_1}{4\pi D^{\frac{3}{2}}}
     \left[\frac{d_{1j}\kappa_j^2}{d_{2j}}\right]_{j}\,, \qquad
     N \equiv [\kappa_k c_{1jk}]_{jk}\,, \\
  &&M_1 \equiv \left[\kappa_k\left(c_{2jk} + \sum_{\ell=-n}^n \omega_{jk\ell}
     b_{3\ell}\right)\right]_{jk}
     -\frac{1}{4\pi D^{\frac{3}{2}}} [\kappa_j]_j
     \left(\sum_{\ell=-n}^n e_{\ell} b_{C0\ell}\right) \,, \\
&&M_0 \equiv \left[\kappa_k \sum_{\ell=-n}^n \omega_{jk\ell} b_{4\ell}\right]_{jk} +
   \frac{1}{4\pi D^{\frac{3}{2}}} [\kappa_j]_j \left(
   \sum_{\ell=-n}^n e_\ell b_{C0\ell} - g\left(2\Delta t\right)\right) \,.
\end{eqnarray*}
Here $e_{\ell}$ and $\omega_{jkl}$ are defined in
\eqref{alg:el_omega}, $b_1$, $b_{3\ell}$, $b_{4\ell}$ and $b_{C0\ell}$
are defined in \eqref{alg:C:values}, while $c_{1jk}$ and $c_{2jk}$ are
defined in \eqref{alg:D:values}.  In deriving \eqref{alg:b2} we
evaluated the derivatives appearing in $C_j(2\Delta t)$ of
\eqref{eq:CLj_update} directly from the differential equation
\eqref{tm:ode}.

\textbf{Next steps:} $t_i \curvearrowright t_{i+1}$, for $i \geq 2$: We
use the RK4 method with a lagged $\v{B}^{(i)}$, so that
\begin{equation}\label{alg:step_i}
  \begin{split}
  \v{u}_{1}^{(i+1)} &= \v{u}_{1}^{(i)} + \frac{\Delta t}{6} (k_1^{(i)} +
  2 k_2^{(i)} + 2 k_3^{(i)} + k_4^{(i)}) + \varepsilon
  \left[\frac{d_{1j}}{d_{2j}} \kappa_j \right]_{j} \v{B}^{(i)} \Delta t\,, \\
  \v{u}_{2}^{(i+1)} &= \v{u}_{2}^{(i)} + \frac{\Delta t}{6}
  (\tilde{k}_1^{(i)} + 2 \tilde{k}_2^{(i)} + 2 \tilde{k}_3^{(i)} +
  \tilde{k}_4^{(i)})  \,.
\end{split}
\end{equation}
In terms of the computed $\v{u}_{1}^{(i+1)}$ we calculate $\v{B}^{(i+1)}$ as
\begin{equation}\label{step:bi}
\begin{split}
  \v{B}^{(i+1)} &= -\frac{1}{4\pi D^{3/2}} A^{-1} \left[\kappa_j\right]_j
  \left( b_1 \hat{\v{F}}_1(\v{u}^{(i+1)}) + (b_2 + \sum_{\ell=-n}^n e_{\ell}
      b_{3\ell})\, \hat{\v{F}}_1(\v{u}^{(i)}) +
    (\sum_{\ell=-n}^n e_{\ell} b_{4\ell})  \hat{\v{F}}_1(\v{u}^{(i-1)})
    \right)  \\
  & \quad -\frac{\eps}{4\pi D^{3/2}} A^{-1} \left[\frac{d_{1j} \kappa_j^2}{
    d_{2j}} \right]_{j} \left(  (b_2 + \sum_{\ell=-n}^n e_{\ell} b_{3\ell})
    \v{B}^{(i)}  +  (\sum_{\ell=-n}^n e_{\ell} b_{4\ell}) \v{B}^{(i-1)}
    \right)     \\
     & \quad + A^{-1} \left(N \v{u}_1^{(i+1)} + M_2 \v{u}_1^{(i)} +
       M_3 \v{u}_1^{(i-1)} + \v{H}^{(i)}\right) - \frac{g(t_{i+1})}
     {4\pi D^{\frac{3}{2}}} A^{-1}[\kappa_j]_j \v{u}_1^{(0)}\,,
  \end{split}
\end{equation}
where the matrices $M_2$ and $M_3$ are defined by
\begin{equation}\label{alg:a_m_n}
  M_2 \equiv \left[\kappa_k(c_{2jk} + \sum_{\ell=-n}^n \omega_{jk\ell} b_{3\ell})
  \right]_{jk}\,,
  \quad M_3 \equiv \left[\kappa_k \sum_{\ell=-n}^n \omega_{jk\ell} b_{4\ell}
  \right]_{jk}\,.
\end{equation}
In \eqref{step:bi} the history vector
$\v{H}^{(i)}\equiv (H_{1}^{(i)},\ldots,H_{N}^{(i)})^T$ has entries
\begin{subequations}\label{alg:hji_update}
\begin{equation}
  H_{j}^{(i)} = \sum_{\ell=-n}^n \left(- \frac{1}{4\pi D^{\frac{3}{2}}} \kappa_j
    e_{\ell} e^{s_{\ell}\Delta t} H_{Cj\ell}^{(i)}
    + \kappa_k \sum_{k=1, k\neq j}^N \omega_{jk\ell} e^{s_{\ell}\Delta t}
    H_{Dk\ell}^{(i)}\right)\,,
\end{equation}
and is updated with the scheme
\begin{equation}
  \v{H}_{C\ell}^{(i+1)} = 
  \left\{
    \begin{aligned}
      & b_{C0\ell} \v{u}_1^{(1)} - b_{C0\ell} \v{u}_1^{(0)}\,, && i=1 \\
      & \v{H}_{C\ell}^{(i)} e^{s_{\ell}\Delta t} + b_{3\ell} \hat{\v{F}}_1(\v{u}^{(i)}) +
      b_{4\ell} \hat{\v{F}}_1(\v{u}^{(i-1)}) \\
      & \qquad + \varepsilon \left[\frac{d_{1j}}{d_{2j}}
        \kappa_j \right]_{j} (b_{3\ell} \v{B}^{(i)} + b_{4\ell} \v{B}^{(i-1)})\,,
      && i \geq 2 \,,
  \end{aligned}
  \right.
\end{equation}
together with
\begin{equation}
  	\v{H}_{D\ell}^{(i+1)} = 
  \begin{cases}
    b_{3\ell} \v{u}_1^{(1)} + b_{4\ell} \v{u}_1^{(0)}\,, & i=1 \\
    \v{H}_{D\ell}^{(i)} e^{s_{\ell}\Delta t} + b_{3\ell} \v{u}_1^{(i)} + b_{4\ell}
    \v{u}_1^{(i-1)}\,,    & i \geq 2\,.
  \end{cases}
\end{equation}
\end{subequations}

A rather similar scheme was successfully developed and implemented in
\cite{pelzward2D} for studying the synchronization behavior of
diffusively-coupled cells in 2-D. Both schemes are a semi-implicit
operator-splitting scheme in the sense that the intracellular reaction
kinetics $(\hat{\v{F}}_1(\v{u}^{(i)}), \v{F}_2(\v{u}^{(i)}))$ are
integrated explicitly through the RK4 method, the coupling term with
$\v{B}^{(i)}$ is added through a simple explicit Forward-Euler step,
and the coupling function $\v{B}^{(i+1)}$ is implicitly computed using
$\v{u}_1^{(i+1)}$ at the same time step $i+1$.

\subsection{Calculation of Local Contributions}\label{app:local_int}

In this technical appendix, we first summarize the derivation of the
local contribution $C_{Lj}(t)$ in \eqref{alg:c_locint}. Recalling
\eqref{3:lap2} for the two terms in $g(t)$, we integrate by parts to
get
\begin{eqnarray*}
  \int_0^{\Delta t} \frac{e^{-\sigma t}}{\sqrt{\pi t}}\; dt &=&
     \frac{1}{\sqrt{\sigma}} \; \mathrm{erf}(\sqrt{\sigma\Delta t})\,, \\
  \sqrt{\sigma} \int_0^{\Delta t} \mathrm{erf}(\sqrt{\sigma t})\; dt &=&
   \left(\sqrt{\sigma}\Delta t - \frac{1}{2\sqrt{\sigma}}\right)\;
    \mathrm{erf}(\sqrt{\sigma\Delta t}) + \sqrt{\frac{\Delta t}{\pi}}
          e^{-\sigma \Delta t}\,,
\end{eqnarray*}
which are combined to yield $\int_{0}^{\Delta t} g(t)\, dt$ in
\eqref{alg:c_locint}. Similarly, the second integral
$\int_0^{\Delta t} tg(t)\;dt$ in \eqref{alg:c_locint} is calculated
term by term by combining the two integrals
\begin{eqnarray*}
  \int_0^{\Delta t} \sqrt{\frac{t}{\pi}} e^{-\sigma t}\, dt &=&
  -\frac{\sqrt{\Delta t}}{\sigma \sqrt{\pi}} e^{-\sigma \Delta t}
       + \frac{1}{2\sigma^{3/2}} \mathrm{erf}(\sqrt{\sigma \Delta t}) \,, \\
  \sqrt{\sigma}\int_0^{\Delta t} t \; \mathrm{erf}(\sqrt{\sigma t}) \; dt &=&
\left(\frac{1}{2}\sqrt{\sigma}(\Delta t)^2 -\frac{3}{8\sigma^{3/2}}\right)\;
\mathrm{erf}(\sqrt{\sigma\Delta t}) +\frac{\sqrt{\Delta t}}{4\sigma \sqrt{\pi}}
    e^{-\sigma\Delta t} (2\sigma\Delta t + 3)\,.
\end{eqnarray*}

Next, we outline the calculation of the integrals $I_1$ and $I_2$ in
\eqref{alg:D:I1_I2}, which leads to \eqref{alg:D:values} for the local
contribution $D_{Ljk}(t)$. Upon using \eqref{2:reduced_3} for
$G(a_{jk},t)$, we have
\begin{equation}\label{appb:I1}
  I_{1}= \int_{0}^{\Delta t} G(a_{jk},z)\, dz =
  \frac{1}{(4\pi D)^{3/2}} I_{\Delta t} \qquad \mbox{where}
  \qquad I_{\Delta t} \equiv \int_{0}^{\Delta t} \frac{e^{-\sigma z -b/z}}
  {z^{3/2}} \, dz \quad \mbox{with} \quad b\equiv \frac{a_{jk}^2}{4}\,.
\end{equation}
By changing variables as $\xi=\sqrt{b/z}$, and labeling
$c=\sqrt{b/\Delta t}>0$, we can rewrite and evaluate $I_{\Delta t}$
as
\begin{equation}\label{appb:I1_1}
 I_{\Delta t} = \frac{2}{\sqrt{b}} \int_{\sqrt{b/\Delta t}}^{\infty}
 e^{-\xi^2 - \sigma b/\xi^{2}} \, d\xi = \frac{2}{\sqrt{b}}\left( \frac{\sqrt{\pi}}
   {4}\right) \left[ e^{2\sqrt{\sigma b}}
    \mathrm{erfc}\left( c + \frac{\sqrt{\sigma b}}{c}\right) +
    e^{-2\sqrt{\sigma b}}
    \mathrm{erfc}\left( c - \frac{\sqrt{ \sigma b}}{c}\right)\right] \,.
\end{equation}
By setting $b={a_{jk}^2/4}$ and $c=\sqrt{b/\Delta t}$ in \eqref{appb:I1_1}
and combining with \eqref{appb:I1}, we obtain the result in
\eqref{alg:D:Ivalues} for $I_1$.  The result for $I_2$ in
\eqref{alg:D:Ivalues} is obtained by simply differentiating $I_1$
with respect to $\sigma$ as
\begin{equation*}
  I_{2} = \int_{0}^{\Delta t} z G(a_{jk},z)\, dz = -\frac{1}{(4\pi D)^{3/2}}
  \frac{d}{d\sigma} I_{\Delta t} \,.
\end{equation*}
Finally, by using 
$\mathrm{erfc(z)}\sim \left(\sqrt{\pi} z\right)^{-1} e^{-z^2} \left( 1
  - {1/(2z^2)}\right)$ for $z\gg 1$, we estimate for $\Delta t\ll 1$
that
\begin{equation}\label{app:D:I1I2_asy}
  I_{1} \sim \frac{4 \sqrt{\Delta t}}{(4\pi D)^{3/2} a_{jk}^2}
  \mathrm{exp}\left( -\frac{a_{jk}^2}{4\Delta t} - \sigma \Delta t\right) \,,
  \qquad
  I_{2} \sim \frac{4 \left(\Delta t\right)^{3/2}}{(4\pi D)^{3/2} a_{jk}^2}
  \mathrm{exp}\left( -\frac{a_{jk}^2}{4\Delta t} - \sigma \Delta t\right) \,.
\end{equation}
This indicates that $I_1$ and $I_2$ are both very small when $\Delta t\ll 1$.
\end{appendices}

\bibliographystyle{plain}
\bibliography{Reference_new}

@book{hoppen,
  author    = {Hoppensteadt, F. C. and Izhikevich, E. M.},
  title     = {Weakly Connected Neural Networks},
  series    = {Applied Mathematical Sciences},
  volume    = {126},
  publisher = {Springer-Verlag},
  address   = {New York},
  year      = {1997},
  isbn      = {978-0-387-94948-2}
}

@article{ding,
title={Eigenvalues of rank-one updated matrices with some applications},
author= {Ding, J. and Zhu, A.},
journal={Appl. Math. Lett.},
volume={20},
number={12},
year={2007},
pages={1223-1226}
}

@article{ridgway,
  title={Quorum-sensing induced transitions between bistable steady-states for
a cell-bulk {ODE-PDE} model with {Lux} intracellular kinetics},
  author={Ridgway, W. and  Ward, M. J. and Wetton, B. T.},
  journal={J. Math. Bio.},
  volume={84},
  number={1-2},
  year={2021}
}

@article{smjw_quorum,
  title={Localized Signaling Compartments in 2-{D} Coupled by a Bulk Diffusion
  Field: Quorum Sensing and Synchronous Oscillations in the Well-Mixed Limit},
  author={Iyaniwura, S. and Ward, M. J.},
  journal={Europ. J. Appl. Math.},
  volume={32},
  number={5},
  year={2021},
  pages={1001-1031}
}

@article{xu,
	Author = {Xu, B. and Bressloff, P.},
	Journal = {SIAM J. Appl. Math.},
	Volume = {76},
	Number = {3},
	Pages = {1844--1870},
 Title = {A {PDE-DDE} model for cell polarization in fission yeast},
	Year = {2016}}

@article{gomezmarin,
  title = {Self-Sustained Spatiotemporal Oscillations Induced by Membrane-Bulk Coupling},
  author = {Gomez-Marin, A. and Garcia-Ojalvo, J. and Sancho, J. M.},
  journal = {Phys. Rev. Lett.},
  volume = {98},
  issue = {16},
  pages = {168303},
  numpages = {4},
  year = {2007},
  month = {Apr},
  publisher = {American Physical Society}, 
}

@article{turing1d,
  title={The emergence of spatial patterns for compartmental reaction kinetics coupled by two bulk diffusing species with comparable diffusivities},
  author={Pelz, M. and Ward, M. J.},
  journal={Philos. Trans. R. Soc. A},
  volume={381},
  number={2245},
  pages={20220089},
  year={2023},
  publisher={The Royal Society}
}

@misc{pelz2026b,
      title={Dynamical behavior of diffusively coupled scalar differential equations as {W}entzell boundary conditions}, 
      author={Pelz, M.},
      year={2026},
      howpublished={\url{https://arxiv.org/abs/2608.08537}},
      note={arXiv:2608.08537 [math.AP]}
}

@misc{pelz2026a,
  author        = {Pelz, M. and Scheel, A.},
  title         = {Oscillations in a scalar differential equation coupled to a diffusive field},
  year          = {2026},
  eprint        = {2604.01135},
  archivePrefix = {arXiv},
  primaryClass  = {math.AP},
  note          = {arXiv preprint arXiv:2604.01135}
}

@article{mora,
  author    = {Mora, T. and Bialek, W.},
  title     = {Are biological systems poised at criticality?},
  journal   = {J. Stat. Phys.},
  volume    = {144},
  number    = {2},
  pages     = {268--302},
  year      = {2011},
  doi       = {10.1007/s10955-011-0229-4}
}

@article{toiya,
author={Toiya, M. and Gonzalez-Ochoa, H. O. and Vanag, V. K. and Fraden, S. and
Epstein, I. R.},
title={Synchronization of chemical micro-oscillators},
journal={J. Phys. Chem. Lett.},
year={2010},
volume={1},
number={8},
pages={1241--1246},
doi={10.1021/jz100238u}
}

@article{torbensen,
author={Torbensen K. and Rossi F. and  Ristori S. and Abou-Hassan A.},
title={Chemical communication and dynamics of droplet emulsions in
networks of {B}elousov-{Z}habotinsky micro-oscillators produced by microfluidics},
journal={Lab Chip.},
year={2017},
volume={17},
number={7},
pages={1179--1189},
doi = {10.1039/c6lc01583b}
}

@article{MNS,
  title = {Oscillatory chemical reactions in closed vessels},
  author = {Merkin, J. H. and Needham, D. J. and Scott, S. K.},
  journal = {Proc. R. Soc. Lond. A Math. Phys. Sci.},
  volume = {406},
  number = {1831},
  pages = {299--323},
  year = {1986},
}

@article{pelz26,
author = {Pelz, M. and Janu\~{s}onis, S. and Handy, G.},
title = {Compartmental-Reaction Diffusion Framework for Microscale Dynamics of Extracellular Serotonin in Brain Tissue},
journal = {SIAM Journal on Life Sciences},
volume = {1},
number = {3},
pages = {509-553},
year = {2026},
doi = {10.1137/25M1827852}
}

@article{kuramoto1995,
	Author = {Kuramoto, Y.},
	Journal = {Prog. Theor. Phys.},
	Pages = {321--330},
	Title = {Scaling behavior of turbulent oscillators with
                  non-local interaction},
	Volume = {94},
	Number = {3},
	Year = {1995}}

@article{gou2d,
  title={An asymptotic analysis of a 2-{D} model of dynamically active compartments coupled by bulk diffusion},
  author={Gou, J. and Ward, M. J. },
  journal={J. Nonliner Sci.},
  volume={26},
  number={4},
  pages={979--1029},
  year={2016},
  publisher={Springer}
}

@article{show_epstein,
 title={From chemical systems to systems chemistry: Patterns in space and
 time},
 author={Showalter, K. and Epstein, I. R.},
 journal={Chaos},
 volume={25},
 pages={097613},
 year={2015} }

@article{q_survey,
  title={Pattern forming systems coupling linear bulk diffusion to dynamically active membranes or cells},
  author={Gomez, D. and Iyaniwura, S. and Paquin-Lefebvre, F. and Ward, M. J.},
  journal={Phil. Trans. Roy. Soc. A.},
  volume={379},
  pages={20200276},
  year={2021}
}

@article{sarafa_3D,
title={Oscillatory instabilities in dynamically active signalling compartments
coupled via bulk diffusion in a {3-D} spherical domain},
author= {Iyaniwura, S. and Ward, M. J.},
journal={Physica D},
volume={476},
year={2025},
number={134645},
}

@book{AS,
  title={Handbook of mathematical functions, with formulas, graphs and
                  mathematical tables},
  author={Abramowitz, M. and Stegun, I.},
  year={1974},
  publisher={Dover Publications, NY, USA}
}

@article{pelzward2D,
title={Synchronized memory-dependent intracellular oscillations for a cell-bulk {ODE-PDE} model in $\mathbf{R}^2$},
author= {Pelz, M. and Ward, M. J.},
journal={SIAM J. Appl. Dyn. Sys.},
year={2025},
volume={24},
number={2},
pages={1752--1811}
}

@article{viana2020,
 title={An integro-differential equation for dynamical systems with diffusion-mediated coupling},
 author={Aristides, R. P. and Viana, R. L.},
 journal={Nonlinear Dynamics},
 volume={100},
 number={4},
 pages={3759--3770},
 year={2020},
}

@article{greengard_2015,
 title={Efficient sum-of-exponentials approximations for the heat kernel and their applications},
 author={Jiang, S. and Greengard, L. and Wang, S.},
 journal={Adv. Comput. Math.},
 volume={41},
 pages={529--551},
 year={2015},
}

@article{beylkin_2005,
 title={On approximation of functions by exponential sums},
 author={Beylkin, G. and Monz\'{o}n, L.},
 journal={Appl. and Comput. Harm. Anal.},
 volume={19},
 number={1},
 pages={17--48},
 year={2005},
}

@article{beylkin_2010,
 title={Approximation by exponential sums revisited},
 author={Beylkin, G. and Monz\'{o}n, L.},
 journal={Appl. and Comput. Harm. Anal.},
 volume={28},
 number={2},
 pages={131--140},
 year={2010},
}

@book{numerics,
  title={An Introduction to Numerical Analysis},
  author={S\"uli, E. and Mayers, D.},
  year={2003},
  publisher={Cambridge University Press}
}

@article{lopez_1,
 author={L\'{o}pez-Fern\'{a}ndez, M. and Lubich, C. and Sch\"{a}dle, A.},
 title={Adaptive, fast, and oblivious convolution in evolution equations with memory},
 journal={SIAM J. Sci. Comput.},
 volume={30},
 number={2},
 pages={1015--1037},
 year={2008},
}

@article{karoui2022chemical,
  author    = {Karoui, H. and Patwal, P. S. and Pavan Kumar,
                  B. V. V. S. and Martin, N.},
  title     = {Chemical communication in artificial Cells: {B}asic concepts, design and challenges},
  journal   = {Front. Mol. Biosci.},
  volume    = {9},
  pages     = {880525},
  year      = {2022},
  doi       = {10.3389/fmolb.2022.880525},
  publisher = {Frontiers Media S.A.}
}

@article {gou2015,
    AUTHOR = {Gou, J. and Li, Y. X. and Nagata, W. and Ward, M. J.},
     TITLE = {Synchronized oscillatory dynamics for a 1-{D} model of
              membrane kinetics coupled by linear bulk diffusion},
   JOURNAL = {SIAM J. Appl. Dyn. Sys.},
    VOLUME = {14},
      YEAR = {2015},
    NUMBER = {4},
     PAGES = {2096--2137},
}

@article {gou2016,
    AUTHOR = {Gou, J. and Ward, M. J.},
     TITLE = {Oscillatory Dynamics for a Coupled Membrane-Bulk Diffusion Model with {F}itzhugh-{N}agumo Kinetics},
   JOURNAL = {SIAM J. Appl. Math.},
    VOLUME = {76},
      YEAR = {2016},
    NUMBER = {2},
     PAGES = {776--804},
}

@article {paquin_1d,
    AUTHOR = {Paquin-Lefebvre, F. and Nagata, W. and Ward, M. J},
     TITLE = {Pattern formation and oscillatory dynamics in a two-dimensional coupled bulk-surface reaction-diffusion system},
   JOURNAL = {SIAM J. Appl. Math.},
    VOLUME = {80},
    Number = {3},
      YEAR = {2020},
     PAGES = {1520--1545},
}

@article{gou2017,
    AUTHOR = {Gou, J. and Chiang, W.-Y. and Lai, P.-Y. and Ward,
              M. J. and Li, Y.-X.},
     TITLE = {A theory of synchrony by coupling through a diffusive chemical
              signal},
   JOURNAL = {Physica D},
    VOLUME = {339},
      YEAR = {2017},
     PAGES = {1--17},
}

@article{ward2018spots,
	Author = {Ward, M. J},
	Journal = {Nonlinearity},
	Number = {8},
	Pages = {R189},
	Publisher = {IOP Publishing},
	Title = {Spots, traps, and patches: {A}symptotic analysis of
                  localized solutions to some linear and nonlinear
                  diffusive systems},
	Volume = {31},
	Year = {2018}}

@article{dano2,
	Author = {Dan{\o}, S. and S{\o}rensen, P. G. and Hynne, F.},
        journal={Nature},
	Volume = {402},
	Number = {6759},
	Pages = {320--322},
 Title = {Sustained oscillations in living cells},
	Year = {1999}}

@article{rodriguez,
  title={The {K}uramoto model in complex networks},
  author={Rodrigues, F. A. and Peron, T. K. D. and Ji, P. and Kurths, J.},
  journal={Physics Reports},
  volume={610},
  pages={1--98},
  year={2016},
  publisher={Elsevier},
  doi={10.1016/j.physrep.2015.10.008}
}

@article{dano,
	Author = {Dan{\o}, S. and Madsen, M. F. and S{\o}rensen, P. G.},
        journal={PNAS},
	Volume = {104},
	Number = {31},
	Pages = {12732--12736},
 Title = {Quantitative characterization of cell synchronization in yeast},
	Year = {2007}}

@article{danino2010,
  author    = {Danino, T. and Mondrag{\'o}n-Palomino, O. and Tsimring, L. and Hasty, J.},
  title     = {A synchronized molecular clock of uncoupled cellular oscillators},
  journal   = {Nature},
  volume    = {463},
  number    = {7279},
  pages     = {326--330},
  year      = {2010},
  publisher = {Nature Publishing Group}
}

@article{mukwaya2021,
  author    = {Mukwaya, V. and Mann, S. and Dou, H.},
  title     = {Chemical communication at the synthetic cell/living cell interface},
  journal   = {Communications Chemistry},
  volume    = {4},
  number    = {1},
  pages     = {161},
  year      = {2021},
  publisher = {Nature Publishing Group}
}

@article{QS_whiteley,
	Author = {Whiteley, M. and Diggle, S. P. and Greenberg, E. P.},
	Journal = {Nature},
	Volume = {551},
	Pages = {7680},
 Title = {Bacterial quorum sensing: {T}he progress and promise of an emerging reserch areas},
	Year = {November 15 2017}}

@article{Selkov,
  title={Self-Oscillations in Glycolysis 1. {A} Simple Kinetic Model},
  author={Sel'kov, E. E. },
  journal={Eur. J. Biochem.},
  volume={4},
  number={1},
  pages={79--86},
  year={1968},
  publisher={Wiley Online Library}
}

@article{smjw_diff,
  title={Synchrony and oscillatory dynamics for a 2-{D} {PDE-ODE} model of
  diffusion-mediated communication between small signaling compartments},
  author={Iyaniwura, S. and Ward, M. J.},
  journal={SIAM J. Appl. Dyn. Sys.},
  volume={20},
  number={1},
  year={2021},
  pages={438--499}
}

@article{gregor2010,
  title={The onset of collective behavior in social amoebae},
  author={Gregor, T. and Fujimoto, K. and Masaki, N. and Sawai, S.},
  journal={Science},
  volume={328},
  number={5981},
  pages={1021--1025},
  year={2010},
  publisher={American Association for the Advancement of Science}
}

@article{nandy1998,
  title={Cyclic {AMP} oscillations in {D}ictyostelium {D}iscoideum:
                  {M}odels and observations},
  author={Nanjundiah, V.},
  journal={Biophys. Chem.},
  volume={72},
  number={1-2},
  pages={1--8},
  year={1998},
  publisher={Elsevier}
}

@article{de2007dynamical,
  title={Dynamical quorum sensing: {P}opulation density encoded in cellular dynamics},
  author={De Monte, S. and d'Ovidio, F. and Dan{\o}, S. and S{\o}rensen, P. G.},
  journal={PNAS},
  volume={104},
  number={47},
  pages={18377--18381},
  year={2007},
  publisher={National Acad Sciences}
}

@article{muller2006,
  title={Cell--cell communication by quorum sensing and dimension-reduction},
  author={M{\"u}ller, J. and Kuttler, C. and Hense, B. A. and Rothballer, M. and Hartmann, A.},
  journal={J. Math. Biol.},
  volume={53},
  number={4},
  pages={672--702},
  year={2006},
  publisher={Springer}
}

@article{ingalls,
  title     = {From single cells to communities: {M}athematical perspectives on bacterial quorum sensing},
  author    = {Sadr, S.  and Zargar, B. and Aucoin, M. G. and Ingalls, B. P.},
  journal   = {Computational and Structural Biotechnology Journal},
  volume    = {27},
  pages     = {4340--4351},
  year      = {2025},
  doi       = {10.1016/j.csbj.2025.09.030},
}

@article{muller2013,
  title={Approximating the dynamics of communicating cells in a diffusive medium by ODEs—homogenization with localization},
  author={M{\"u}ller, J. and Uecker, H.},
  journal={J. Math. Biol.},
  volume={67},
  number={5},
  pages={1023--1065},
  year={2013},
  publisher={Springer}
}

@article{taga2003chemical,
  title={Chemical communication among bacteria},
  author={Taga, M. E. and Bassler, B. L.},
  journal={PNAS},
  volume={100},
  number={suppl 2},
  pages={14549--14554},
  year={2003},
  publisher={National Acad Sciences}
}

@book{goldbeter1997biochemical,
  title={Biochemical oscillations and cellular rhythms: {T}he molecular bases of periodic and chaotic behaviour},
  author={Goldbeter, A.},
  year={1997},
  publisher={Cambridge University Press}
}

@article{taylor2,
 title={Insights into collective cell behavior from populations of
                  coupled chemical oscillators},
 author={Taylor, A. F. and Tinsley, M. and Showalter, K.},
 journal={Phys. Chemistry Chem Phys.},
 volume={17},
 number={31},
 pages={20047--20055},
 year={2015} }

@article{tinsley1,
 title={Dynamical quorum sensing and synchronization in
collections of excitable and oscillatory catalytic particles},
 author={Tinsley, M. R. and Taylor, A. F. and Huang, Z. and  Wang, F.  and
 Showalter, K.},
 journal={Physica D},
 volume={239},
 number={11},
 pages={785--790},
 year={2010} }

@article{tinsley2,
 title={Emergence of collective behavior in groups of excitable 
catalyst-loaded particles: {S}patiotemporal dynamical quorum sensing},
 author={Tinsley, M. R. and Taylor, A. F. and Huang, Z. and
 Showalter, K.},
 journal={Phys. Rev. Lett.},
 volume={102},
 pages={158301},
 year={2009} }

@article{dicty,
 title={Modeling oscillations and spiral waves in {\em Dictyostelium} populations},
 author={Noorbakhsh, J. and Schwab, D. J. and Sgro, A. E. and Gregor, T. and Mehta, P.},
 journal={Phys Rev. E.},
 volume={91},
 number={6},
 pages={062711},
 year={2015} }

\end{document}